\documentclass[12pt]{article}
\usepackage{graphicx} 
\usepackage[all]{xy}
\usepackage{pdflscape}
\usepackage[hypcap=false]{caption}
\usepackage{amsmath,amssymb,amsfonts,latexsym,float,graphics,epsfig} 
\usepackage{setspace}
\usepackage{xcolor}
\usepackage[top=2cm, bottom=2cm, left=2cm, right=2cm]{geometry}
\usepackage[colorlinks]{hyperref}
\usepackage{tikz}
\usetikzlibrary{cd}

\newtheorem{theorem}{Theorem}[section]

\newtheorem{remark}[theorem]{Remark}

\newcommand{\Vol}{\operatorname{Vol}}
\newcommand{\qq}{\mathbf{q}}
\newcommand{\pp}{\mathbf{p}}

\title{On Continuum Theories in a Locally Conformal Setting}
\author{}
\date{}

\begin{document}

\maketitle

\begin{center} 
 Begüm Ateşli$^{\ast,\ast\ast}$\footnote{e-mail: 
\href{mailto:begum.atesli@medipol.edu.tr}{begum.atesli@medipol.edu.tr}, \href{mailto:begumatesli@itu.edu.tr}{begumatesli@itu.edu.tr}, corresponding author,}, 
O\u{g}ul Esen$^{\dagger,\dagger\dagger}$\footnote{e-mail: \href{mailto:oesen@gtu.edu.tr}{oesen@gtu.edu.tr},}, 
Miroslav Grmela$^{\ddagger}$\footnote{e-mail: \href{mailto:miroslav.grmela@polymtl.ca} {miroslav.grmela@polymtl.ca},}, 
  Michal Pavelka$^{\ddagger\ddagger}$\footnote{e-mail: \href{mailto:pavelka@karlin.mff.cuni.cz}{pavelka@karlin.mff.cuni.cz}} 
  \\

\bigskip

$^{\ast}$School of Business \\ İstanbul Medipol University, 34815 Kavacık, İstanbul, Turkey

\bigskip

$^{\ast\ast}$Department of Mathematics Engineering \\ İstanbul Technical University,  34467 Maslak, İstanbul, Turkey

\bigskip

$^{\dagger}$Department of Mathematics, \\ Gebze Technical University, 41400 Gebze,
Kocaeli, Turkey
\bigskip

$^{\dagger\dagger}$Center for Mathematics and its Applications,  \\ 
Khazar University, Baku, AZ1096, Azerbaijan
\bigskip

$^{\ddagger}$\'{E}cole Polytechnique de Montr\'{e}al,
  C.P.6079 suc. Centre-ville, 
 \\  Montr\'{e}al, H3C 3A7,  Qu\'{e}bec, Canada
\bigskip

$^{\ddagger\ddagger}$Mathematical Institute, Faculty of Mathematics and Physics,  \\  Charles University, 
 Sokolovsk\'{a} 83, 18675 Prague, Czech Republic

\bigskip 
\begin{abstract}
We study continuum and kinetic dynamics on locally conformally symplectic (LCS) and locally conformally cosymplectic (LCCos) manifolds. The main geometric tool is the Lichnerowicz--de~Rham differential, which allows local geometric data that are compatible only up to conformal rescaling to glue into global structures.

In the LCS setting, we analyze the locally conformal Vlasov equation from the viewpoint of continuity. Its source term does not represent particle creation or annihilation: it measures the mismatch between the ordinary phase-space volume and the conformally weighted volume determined by the Lee form. Absorbing this weight into the density yields a source-free continuity equation. We also compute the entropy production of the transported density, and illustrate the framework through a thermostatted ensemble, an LCS-driven Vlasov model with a momentum-dependent gyroscopic force, and their combination.

We then propose a new definition of LCCos manifolds in which the defining one-form and two-form carry the same conformal weight. Unlike the usual convention, this allows local cosymplectic Hamiltonian systems to glue directly into global LCCos Hamiltonian dynamics, and it normalizes the evolution vector field by the conformal unit rather than by the ordinary unit. We develop the corresponding Hamiltonian, evolution and continuum dynamics, and illustrate them with a time-dependent mass example exhibiting transport in the distinguished time direction, momentum-space drift, and conformal source terms.

\medskip 

\smallskip
 

\end{abstract}
\end{center}

\tableofcontents
\onehalfspacing

\numberwithin{equation}{section}

\setlength{\parindent}{2em}
\setlength{\parskip}{3ex}
 
\section{Introduction}

Poisson geometry provides one of the fundamental geometric frameworks for
Hamiltonian dynamics \cite{Laurent13,vaisman2012lectures,Weinstein98}. Its
importance extends simultaneously to algebra, geometry, topology, classical
mechanics, and continuum and kinetic theories; see, for example,
\cite{AbrahamMarsden,Arnold,guil90,LiMa,MaRa99,morrison1982poisson,pavelkavd2016}.
At the dynamical level, the Jacobi identity encodes the compatibility of
Hamiltonian evolution with the Poisson bracket and guarantees that Hamiltonian
vector fields form an involutive family under the Jacobi--Lie bracket. In
addition, Poisson Hamiltonian dynamics preserves the Hamiltonian function
itself. When the Hamiltonian represents the total energy, this expresses the
conservative character of the underlying dynamics.

This close relation between Poisson geometry and conservative dynamics raises
a natural question: how much of the Hamiltonian and kinetic formalism survives
when the underlying geometric structure is no longer globally closed, but only
closed in a locally conformal sense? The answer developed in this paper is
twofold: integrability is retained through the underlying Jacobi structure,
whereas volume preservation is recovered in a conformally weighted sense.
The present work is devoted to this question for symplectic and cosymplectic
continuum dynamics. Our guiding idea is to replace the ordinary de~Rham
closedness condition by closedness with respect to the
Lichnerowicz--de~Rham differential \cite{GuLi84}. In this way, locally defined
geometric data may fail to glue globally in the ordinary sense, but become
compatible after conformal rescaling. Thus, locally conformal geometry
provides a natural setting in which global geometric and dynamical structures
arise from conformally compatible local Hamiltonian systems.

Locally conformally symplectic manifolds were introduced by Vaisman
\cite{Vaisman85} and have since been studied in several geometric and
dynamical directions; see, for instance,
\cite{Banyaga,Bazzoni2018,EsenLeonSarZaj1,esen2023conformal,EsLaLeMa21,HallerRybicki,OtimanStanciu2017}.
An LCS manifold carries a natural Jacobi structure, and this Jacobi structure
allows one to define a Hamiltonian dynamics governed by the Lichnerowicz--de~Rham
differential \(d_\varphi\). Although such a dynamics is not conservative in
the ordinary sense, it retains an integrability mechanism through the
underlying Jacobi bracket. This makes the LCS framework suitable not only for
particle dynamics, but also for kinetic and continuum formulations. 
 
\noindent \textbf{Goal 1. Continuity Analysis of the LCS Kinetic Dynamics.}
In the first part of the paper, we revisit LCS Hamiltonian dynamics from this
viewpoint and lift it to densities by means of a Lie--Poisson construction.
The resulting locally conformal Vlasov equation was obtained in
\cite{esen2023conformal}. In the present work, we further analyze this equation
from the viewpoint of continuity. We show that the source term in the locally
conformal Vlasov equation does not represent a genuine creation or annihilation
of particles. Instead, it measures the mismatch between the ordinary
phase-space volume and the conformally weighted volume naturally associated
with the Lee form. This also yields an equivalent formulation: after absorbing
the conformal weight into the density function, the same dynamics becomes a
source-free continuity equation with respect to the ordinary volume. Additionally, we 
analyze the entropy evolution of this density and show that the sign of the
entropy production is governed by the Lee form contribution along the
Hamiltonian flow.

\noindent \textbf{Goal 2. Three Applications of LCS Kinetic Dynamics.}
We apply the LCS kinetic framework to three continuum models. The first is
a Nos\'e--Hoover thermostatted ensemble, where the underlying thermostat dynamics is realized
as an LCS Hamiltonian system and the corresponding kinetic equation is written
as a Liouville equation for the ensemble. The second is an LCS-driven continuum model with a momentum-dependent
gyroscopic term. Here the magnetic symplectic structure underlying the
Lorentz force law is replaced by an LCS structure. The Lee form then
produces a momentum-dependent force which is orthogonal to the momentum
and therefore gyroscopic, rather than an ordinary magnetic force on the
configuration space, and the associated Vlasov-type equation describes a
continuum of particles governed by this law. Finally, we combine
these two mechanisms in a third model, where the Nos\'e--Hoover thermostatted
dynamics is coupled with the LCS-driven momentum-dependent gyroscopic contribution.

\noindent \textbf{Goal 3. A natural local-to-global formulation of LCCos manifolds.}
Cosymplectic geometry provides the standard geometric setting for Hamiltonian
systems in which a distinguished time variable is present; see, for instance,
\cite{Cape,de2017cosymplectic,LeonTuyn}, see also \cite{zajkac2023time}. Its locally conformal counterpart, namely the theory of locally conformally
cosymplectic manifolds (LCCos), has received less attention, although related
approaches exist \cite{AtEsLe23,ChLeMa91}. In the present paper, we adopt a novel formulation of
locally conformally cosymplectic geometry which is directly adapted to the
conformal gluing of local cosymplectic Hamiltonian systems.

The point of this formulation is not merely a change of conformal weights.
Rather, it is designed so that the local Hamiltonian vector fields themselves
glue directly into a global LCCos Hamiltonian vector field. In this sense, the
locally conformal deformation is placed where it naturally belongs: while the
local Hamiltonian dynamics remains compatible under gluing, the tensorial
geometric data acquire the corresponding conformal transformation laws. This
is parallel to the well-established LCS case, where locally symplectic
Hamiltonian systems glue into a global LCS Hamiltonian dynamics.

With this choice, the resulting LCCos structure admits a natural Jacobi
structure, an associated LCCos bracket, and a bracket--vector field relation
compatible with the global Hamiltonian dynamics. This compatibility is then
used to lift the LCCos Hamiltonian dynamics to densities and to derive the
corresponding continuum theory. Throughout the paper, we compare this
formulation with the alternative conformal-weight convention and indicate
explicitly how the present choice affects the Hamiltonian equations, the
evolution equations, and their kinetic lifts. 
 
\noindent \textbf{Goal 4. Continuum dynamics on LCCos manifolds.}
We develop the LCCos kinetic theory by computing the coadjoint action of the
function algebra, deriving the LCCos Vlasov equation, and rewriting it in
transport, continuity, and weighted-density forms. As in the LCS case, the
conformal source term admits two equivalent interpretations: either one works
with the conformally weighted almost cosymplectic volume, or one uses the
ordinary volume together with the weighted density. We also provide coordinate
formulations on the extended cotangent bundle and show explicitly how the Lee
form produces transport in the distinguished time direction, momentum-space
drift, and conformal source terms.

Finally, we illustrate the LCCos construction for particle and continuum dynamics of a time-dependent mass
Hamiltonian. In the pure cosymplectic case, the Hamiltonian Vlasov equation
describes free transport on fixed time-slices. In the LCCos case, however, the
Lee form couples the spatial motion to the distinguished time direction and
modifies both the particle dynamics and the corresponding Vlasov equation. We
also compare the Hamiltonian picture with the evolution dynamics, where the
normalization condition produces an intrinsic time-propagation term.

\noindent \textbf{The Locally Conformal Geometric Dynamics Program.}
The present paper belongs to a broader program devoted to locally conformal
and conformal extensions of geometric dynamics. In this direction, the
Hamilton--Jacobi formalism on locally conformally symplectic manifolds was
developed in \cite{EsenLeonSarZaj1}, while locally conformally cosymplectic
Hamiltonian dynamics and Hamilton--Jacobi theory were studied in
\cite{AtEsLe23}, following the earlier LCCos framework of \cite{ChLeMa91}.
The variational and kinetic aspects of locally conformal dynamics were then
investigated in \cite{esen2023conformal}, where the Lie--Poisson formulation
of the locally conformal Vlasov equation was first obtained.

This locally conformal program has also been developed in several related
directions. Higher-order Lagrangian systems were studied in the locally
conformal setting in \cite{Serdar}, while discrete Lagrangian
and Hamiltonian dynamics were formulated on locally conformal frameworks in
\cite{discrete-LCS}. In the field-theoretic direction, conformal
classical field theories were examined through Cauchy data spaces and
multisymplectic geometry in \cite{EsLeSaZa-Cauchy} and \cite{Esen21-HDW}. A
closely related kinetic counterpart was developed in
\cite{Contact-kinetik}, where conformal and
contact kinetic dynamics were geometrized.

We would also like to draw attention to our recent preprint \cite{coco-begum-1}, which develops LCCos particle dynamics within the same equal-weight framework. The second part of the present paper may thus be regarded as a continuum extension of that work in the LCCos setting.

\noindent \textbf{Contents.}
The paper is organized as follows. Section~\ref{sec:2} is devoted to the
locally conformally symplectic part of the theory. We first review LCS
Hamiltonian and evolution dynamics, the associated Jacobi bracket, and the
local coordinate realization on cotangent bundles. In
Subsection~\ref{Sec-Sym-Man}, we recall the LCS kinetic theory in its
Lie--Poisson formulation, analyze the source term from the viewpoint of
continuity, and derive a weighted-density formulation. Three applications are presented: the Nos\'e--Hoover thermostatted
ensemble in Subsection~\ref{sec-NH}, the LCS-driven gyroscopic model
in Subsection~\ref{sec-LCS-MV}, and their combination in
Subsection~\ref{sec:NH-LCS-EM}.

Section~\ref{sec:3} is devoted to the locally conformally cosymplectic part of
the paper. In Subsection~\ref{sec:LCCosM}, we introduce the LCCos formulation
adopted here and explain its local-to-global meaning in comparison with the
alternative conformal-weight convention. In Subsection~\ref{sec:LCCosH}, we
derive the corresponding LCCos Hamiltonian and evolution dynamics, together
with their bracket formulations and coordinate expressions on the extended
cotangent bundle. Finally, in Subsection~\ref{Sec-LCCos-Man}, we develop the
LCCos kinetic theory by computing the coadjoint action, deriving the LCCos
Vlasov equation, and rewriting it in transport, continuity, and
weighted-density forms. The time-dependent mass Hamiltonian example, including
both particle and continuum dynamics, is presented in
Subsection~\ref{sec-time-mass}.

The appendix collects the background material and convention-dependent identities used throughout the paper. In particular, Subsection~\ref{sec-LdR} recalls the Lichnerowicz--de~Rham differential, the conformal unit, and the basic local-to-global interpretation of locally conformal structures. The appendix also fixes the sign and contraction conventions for the Schouten--Nijenhuis bracket, bivector-induced sharp maps, and Jacobi--Lie brackets, in order to clarify the Poisson and Jacobi identities used in the main text.

\section{Conformal Analysis of Symplectic Continuum Dynamics}\label{sec:2}

\subsection{Locally Conformal Symplectic Dynamics}\label{sec:lcs} 

\noindent \textbf{(Almost) Symplectic Dynamics.}
A two-form \(\Omega\) on a \(2n\)-dimensional manifold \(M\) is called almost
symplectic \cite{Arnold,holm2009geometric,LiMa-book,MaRa99} if it is
non-degenerate, that is, if
\begin{equation}
    \operatorname{Vol}_{\Omega}
    =
    \Omega^n
\end{equation}
is nowhere vanishing. We denote an almost symplectic manifold by the pair
\((M,\Omega)\).
The non-degeneracy of \(\Omega\) induces an isomorphism, called the flat map,
from the space \(\mathfrak{X}(M)\) of vector fields to the space
\(\Gamma^1(M)\) of one-form sections:
\begin{equation}\label{flat}
\flat:\mathfrak{X}(M)
\longrightarrow \Gamma^{1}(M),
\qquad
X\mapsto \iota_X\Omega=\Omega(X,\bullet).
\end{equation}
Here, \(\iota_X\) denotes the interior product, or contraction, with the vector
field \(X\). The inverse of \(\flat\) is the sharp map
\begin{equation}\label{sharp}
\sharp:\Gamma^{1}(M)\longrightarrow \mathfrak{X}(M),
\qquad
\sharp=\flat^{-1}.
\end{equation}

For a Hamiltonian function \(H\) on an almost symplectic manifold, the almost
symplectic Hamiltonian dynamics is determined by Hamilton's equation
\begin{equation}\label{HamEq}
    \iota_{X_H}\Omega=dH,
    \qquad \text{or, equivalently,} \qquad
    X_H=\sharp(dH).
\end{equation}
The dynamics is governed by the Hamiltonian vector field \(X_H\). Since the
almost symplectic two-form is non-degenerate, the sharp map is an isomorphism,
and hence \(X_H\) is uniquely determined by the one-form \(dH\).

We define the bivector field \(\Lambda\) associated with the almost symplectic
form \(\Omega\) by using the sharp map in \eqref{sharp}. For two one-forms
\(\alpha\) and \(\beta\), we set
\begin{equation}\label{Lambda}
    \Lambda(\alpha,\beta)
    =
    \Omega(\sharp\alpha,\sharp\beta).
\end{equation}
A direct computation shows that the sharp map \(\sharp_\Lambda\) induced by
\(\Lambda\), defined by right contraction, coincides with the sharp map
\(\sharp\) in \eqref{sharp}. Namely,
\begin{equation}
    \sharp_\Lambda\alpha
    =
    \iota_\alpha\Lambda
    =
    \Lambda(\bullet,\alpha)
    =
    \sharp\alpha .
\end{equation}
Referring either to the almost symplectic two-form \(\Omega\) or, equivalently,
to the bivector field \(\Lambda\), one defines a bracket on smooth functions by
\begin{equation}\label{bra-sym}
    \{F,H\}
    =
    \Omega(X_F,X_H)
    =
    \Lambda(dF,dH).
\end{equation}
This bracket is skew-symmetric and satisfies the Leibniz identity. In general,
however, it does not satisfy the Jacobi identity, and therefore it defines only
an almost Poisson bracket.

\noindent \textbf{Symplectic Dynamics.}
An almost symplectic two-form is called symplectic if it is closed, that is, if 
    $d\Omega=0$. 
In this case, the bivector field \(\Lambda\) defined in \eqref{Lambda}
commutes with itself under the Schouten--Nijenhuis bracket:
\begin{equation}
    [\Lambda,\Lambda]_{SN}=0.
\end{equation}
Equivalently, the Poisson condition \eqref{Poisson-ident-app} is satisfied.
Hence, the bracket \eqref{bra-sym} becomes a Poisson bracket by satisfying the
Jacobi identity.

Accordingly, the family of Hamiltonian vector fields is closed under the
Jacobi--Lie bracket. More precisely, with the above sign convention, one has
\begin{equation}
    [X_F,X_H]
    =
    -X_{\{F,H\}}.
\end{equation}
Consequently, the characteristic distribution generated by Hamiltonian vector
fields is involutive, and hence integrable.
In the symplectic case, the Hamiltonian flow preserves both the symplectic
two-form and the Hamiltonian function:
\begin{equation}
    \mathcal{L}_{X_H}\Omega=0,
    \qquad
    \mathcal{L}_{X_H}H=0.
\end{equation}
Thus the symplectic structure is invariant along the motion, and, when \(H\)
is interpreted as the total energy, the second identity expresses conservation
of energy.

\noindent \textbf{Locally Conformally Symplectic Dynamics.}
An almost symplectic two-form \(\Omega\) on a manifold \(M\) is called
locally conformally symplectic, abbreviated as LCS, \cite{Vaisman85} if there
exists a closed one-form \(\varphi\), called the Lee one-form, such that
\begin{equation}\label{Omega-varphi}
    d\Omega=\varphi\wedge\Omega,
    \qquad \text{or, equivalently,} \qquad
    d_{\varphi}\Omega=0 .
\end{equation}
Here \(d_{\varphi}\Omega=d\Omega-\varphi\wedge\Omega\) denotes the
Lichnerowicz--de~Rham differential defined in \eqref{LdR-Diff-}. We denote an
LCS manifold by \((M,\Omega\,;\varphi)\).

On an LCS manifold, one may still use the almost symplectic Hamiltonian
equation \eqref{HamEq}. In the present context, we call this the LCS evolution
dynamics. Thus, for a Hamiltonian function \(H\), the corresponding evolution
vector field \(X_H\) is defined by
\begin{equation}\label{evo-dyn-LCS}
    \iota_{X_H}\Omega=dH .
\end{equation}
In this case, the Hamiltonian function is preserved along the flow, whereas
the LCS two-form is not preserved in general. Indeed,
\begin{equation}
    \mathcal{L}_{X_H}H=0,
    \qquad
    \mathcal{L}_{X_H}\Omega
    =
    \varphi(X_H)\Omega-\varphi\wedge dH .
\end{equation}
However, the evolution vector fields are not closed, in general,
under the Jacobi--Lie bracket. This failure of closure should be
distinguished from the integrability of the distribution they
generate. Indeed, since $\Omega$ is non-degenerate, the vectors
$X_H(x)=\sharp_x(dH_x)$, as $H$ varies, span $T_xM$ at every point
$x\in M$. Thus the generated distribution is the full tangent
bundle and is integrable, even though the family of evolution
vector fields need not form a Lie algebra.

In the LCS setting, one may instead consider what we call the LCS Hamiltonian
dynamics. The LCS Hamiltonian vector field \(\mathcal{E}_H\) is defined by
\begin{equation}\label{LCS-Ham-global}
    \iota_{\mathcal{E}_H}\Omega=d_{\varphi}H .
\end{equation}
In this case, neither the Hamiltonian function nor the LCS two-form is
preserved in general. More precisely,
\begin{equation}\label{cons-LCS-Ham}
    \mathcal{L}_{\mathcal{E}_H}H
    =
    H\,\varphi(\mathcal{E}_H),
    \qquad
    \mathcal{L}_{\mathcal{E}_H}\Omega
    =
    \langle\varphi,\mathcal{E}_H\rangle\Omega .
\end{equation}
Accordingly, the divergence of the LCS Hamiltonian vector field with respect
to the LCS volume form is
\begin{equation}\label{LCS-Ham-divergence}
    \mathcal{L}_{\mathcal{E}_H}\mathrm{Vol}_\Omega
    =
    n\langle\varphi,\mathcal{E}_H\rangle \mathrm{Vol}_\Omega,
    \qquad
    \operatorname{div}_{\Omega}(\mathcal{E}_H)
    =
    n\langle\varphi,\mathcal{E}_H\rangle .
\end{equation}

\noindent \textbf{Jacobi Bracket in LCS Geometry.}
An LCS manifold \((M,\Omega\,;\varphi)\) carries a natural Jacobi bracket
\begin{equation}\label{LCS-Bracket}
    \{F,H\}^{\mathrm{LCS}}
    =
    \Lambda(dF,dH)
    +
    F\,\mathcal{Z}_{\varphi}(H)
    -
    H\,\mathcal{Z}_{\varphi}(F),
\end{equation}
where \(\mathcal{Z}_{\varphi}\) is the Lee vector field associated with
\(\varphi\), defined by
\begin{equation}
    \mathcal{Z}_{\varphi}=\sharp\varphi .
\end{equation}
The bracket \eqref{LCS-Bracket} is skew-symmetric and satisfies the
first-order differential-operator identity \eqref{Jacobi-bracket-Leibniz}.
Since it also satisfies the Jacobi identity, it defines a Jacobi bracket in
the sense of \eqref{Jacobi-bracket}. With the present sign convention, the
associated Jacobi pair is \((\Lambda,-\mathcal{Z}_{\varphi})\). Equivalently,
this pair satisfies
\begin{equation}\label{ident-Jac-LCS}
    [\Lambda,\Lambda]
    =
    -2\,\mathcal{Z}_{\varphi}\wedge\Lambda,
    \qquad
    [\mathcal{Z}_{\varphi},\Lambda]=0 .
\end{equation}
The relation between the LCS bracket and the LCS Hamiltonian vector field is
\begin{equation}\label{LCS-bracket-vector-relation}
    \{F,H\}^{\mathrm{LCS}}
    =
    \mathcal{E}_H(F)
    -
    F\langle\varphi,\mathcal{E}_H\rangle .
\end{equation}

As a consequence, the family of LCS Hamiltonian vector fields is closed under
the Jacobi--Lie bracket. More precisely, in light of \eqref{Jac-integ}, one has
\begin{equation}\label{sennur}
    [\mathcal{E}_F,\mathcal{E}_H]
    =
    -\mathcal{E}_{\{F,H\}^{\mathrm{LCS}}}.
\end{equation}
Thus, although the evolution vector fields on an LCS manifold
need not be closed under the Jacobi--Lie bracket, the twisted
Hamiltonian vector fields do form a Lie algebra through the
underlying Jacobi structure. The local-to-global
construction presented below shows that the twisted Hamiltonian equation
\eqref{LCS-Ham-global} and the bracket \eqref{LCS-Bracket} are not external
choices; rather, they are naturally dictated by the geometry of locally
conformally symplectic manifolds. We denote the space of LCS Hamiltonian vector
fields by
\begin{equation}\label{LCS-Ham-vfield}
    \mathfrak{X}_{\mathrm{LCS-ham}}(M)
    =
    \bigl\{
    \mathcal{E}_H\in\mathfrak{X}(M)
    :
    \iota_{\mathcal{E}_H}\Omega=d_\varphi H
    \bigr\}.
\end{equation}

For later use, we record the following integration identity. Assuming that the
boundary terms vanish, equations \eqref{LCS-Ham-divergence} and
\eqref{LCS-bracket-vector-relation} imply \cite{esen2023conformal}
\begin{equation}\label{LC-assoc}
\int_M
K\{H,F\}^{\mathrm{LCS}}\,\mathrm{Vol}_\Omega
=
\int_M
\left(
\{K,H\}^{\mathrm{LCS}}
+
(n+2)\langle\varphi,\mathcal{E}_H\rangle K
\right)
F\,\mathrm{Vol}_\Omega .
\end{equation}
In the symplectic limit \(\varphi=0\), one has
\(\mathcal{E}_H=X_H\), \(\mathrm{div}_{\Omega}(X_H)=0\), and the LCS bracket
reduces to the canonical Poisson bracket. Hence \eqref{LC-assoc} becomes the
usual skew-adjointness relation
\[
\int_M
K\{H,F\}^{(S)}\,\mathrm{Vol}_\Omega
=
\int_M
\{K,H\}^{(S)}F\,\mathrm{Vol}_\Omega .
\]

\noindent \textbf{From Local to Global.}
Let \(M\) be an even-dimensional manifold with an open cover
\(\{U_\alpha\}\). Assume that each \(U_\alpha\) carries a symplectic form
\(\Omega_\alpha\), and that there exist smooth functions
\(\sigma_\alpha\) such that
\begin{equation}
    e^{\sigma_\alpha}\Omega_\alpha
    =
    e^{\sigma_\beta}\Omega_\beta
    \qquad
    \text{on } U_\alpha\cap U_\beta .
\end{equation}
Then the rescaled forms glue to a global two-form \(\Omega\), locally given by
\begin{equation}\label{bb-1}
    \Omega\vert_\alpha
    =
    e^{\sigma_\alpha}\Omega_\alpha .
\end{equation}
We assume that $d\sigma_\alpha=d\sigma_\beta$ on overlaps.
These local one-forms then define a global closed Lee form
$\varphi$ by
\begin{equation}\label{bb-2}
    \varphi\vert_\alpha=d\sigma_\alpha .
\end{equation}
Since \(d\Omega_\alpha=0\), one obtains locally
\(d\Omega\vert_\alpha=\varphi\vert_\alpha\wedge\Omega\vert_\alpha\), hence
\(d\Omega=\varphi\wedge\Omega\), or equivalently \(d_\varphi\Omega=0\). Thus
\((M,\Omega\,;\varphi)\) is an LCS manifold.

We now apply the same gluing principle to Hamiltonian dynamics. On each
symplectic chart \((U_\alpha,\Omega_\alpha)\), let
\begin{equation}\label{loc-Ham-Eq}
    \iota_{X_\alpha}\Omega_\alpha=dH_\alpha
\end{equation}
be the local Hamiltonian equation. Assume that the local vector fields glue to
a global vector field \(\mathcal E_H\), and that the local Hamiltonians glue
conformally, that is,
\begin{equation}\label{LCS-localHamEq}
    \mathcal E_H\vert_\alpha=X_\alpha,
    \qquad
    H\vert_\alpha
    =
    e^{\sigma_\alpha}H_\alpha .
\end{equation}
Then the local equations \eqref{loc-Ham-Eq} glue precisely to the global LCS Hamiltonian equation
\eqref{LCS-Ham-global}. Hence the LCS Hamiltonian equation is not an additional
prescription; it is the global form of conformally compatible local symplectic
Hamiltonian systems. 

\noindent \textbf{Cotangent Bundle.} 
Cotangent bundles provide the canonical geometric framework for symplectic
dynamics. Let \(Q\) be a smooth manifold, and consider its cotangent bundle
\(T^*Q\). The cotangent bundle carries the canonical, or Liouville, one-form
\(\Theta_Q\), whose negative exterior derivative defines the canonical
symplectic form. In induced canonical coordinates \((q^i,p_i)\), these forms
are written as
\begin{equation}\label{can-forms-coord}
    \Omega_Q=-d\Theta_Q,
    \qquad
    \Theta_Q\vert_\alpha=p_i\,dq^i,
    \qquad
    \Omega_Q\vert_\alpha=dq^i\wedge dp_i .
\end{equation}
For a Hamiltonian function \(H\), Hamilton's equation \eqref{HamEq} gives
\begin{equation}\label{Ham-symp-loc}
    \dot q^i
    =
    \frac{\partial H}{\partial p_i},
    \qquad
    \dot p_i
    =
    -\,\frac{\partial H}{\partial q^i}.
\end{equation}

Cotangent bundles also provide natural examples of LCS manifolds
\cite{ChMu17,OtimanStanciu2017}. Let \(\vartheta\) be a closed one-form on the
base manifold \(Q\). By pulling it back through the cotangent bundle projection
\(\pi_Q:T^*Q\to Q\), we obtain a closed semi-basic one-form
\begin{equation}\label{varphi}
    \varphi=(\pi_Q)^*\vartheta .
\end{equation}
We then define the LCS two-form on \(T^*Q\) by
\begin{equation}\label{Omega-phi}
    \Omega
    =
    -d_\varphi\Theta_Q
    =
    -d\Theta_Q+\varphi\wedge\Theta_Q
    =
    \Omega_Q+\varphi\wedge\Theta_Q .
\end{equation}
Since \(d\varphi=0\), this two-form satisfies \(d\Omega=\varphi\wedge\Omega\).
Hence \((T^*Q,\Omega;\varphi)\) is an LCS manifold. Notice that the Lee
one-form \(\varphi\) is defined by pullback from the base \(Q\), and therefore
depends only on the base variables.

\noindent \textbf{Coordinates in LCS Setting.}
On a cotangent chart \(U_\alpha\), there are two useful coordinate systems. The
first one is the canonical coordinate system \((q^i,p_i)\), in which the
canonical forms are given by \eqref{can-forms-coord}. The second one is the
local conformal Darboux coordinate system \((q^i,r_i)\), in which the local
symplectic forms are
\begin{equation}
    \Theta_\alpha=r_i\,dq^i,
    \qquad
    \Omega_\alpha=-d\Theta_\alpha=dq^i\wedge dr_i .
\end{equation}
Locally, the Lee form is exact, and we write
\begin{equation}\label{aks}
    \varphi\vert_\alpha
    =
    d\sigma_\alpha
    =
    \varphi_i(q)\,dq^i,
    \qquad
    \varphi_i(q)
    =
    \frac{\partial\sigma_\alpha}{\partial q^i}.
\end{equation}
The two coordinate systems are related by the conformal rescaling
\(r_i=e^{-\sigma_\alpha}p_i\). Accordingly,
\(\Theta_\alpha=e^{-\sigma_\alpha}\Theta_Q\vert_\alpha\) and
\(\Omega\vert_\alpha=e^{\sigma_\alpha}\Omega_\alpha\), in agreement with the
local-to-global construction above.

In canonical coordinates \((q^i,p_i)\), local realization of the global LCS two-form \eqref{Omega-phi} is 
\begin{equation}\label{glo-glue-symp-2form}
    \Omega\vert_\alpha
    =
    dq^i\wedge dp_i
    -
    \frac12 C_{ij}(q,p)\,dq^i\wedge dq^j ,
\end{equation}
where we have the following skew-symmetric quantities
\begin{equation}\label{Aij}
    C_{ij}(q,p)
    =
    \varphi_j(q)p_i-\varphi_i(q)p_j .
\end{equation}
The first term in \eqref{glo-glue-symp-2form} is the canonical symplectic form, while the second term is the
local contribution of the Lee form through \(\varphi\wedge\Theta_Q\).

The LCS Hamiltonian vector field \(\mathcal E_H\), defined by
\eqref{LCS-Ham-global}, has the following coordinate expression:
\begin{equation}\label{LCS-Ham-local}
    \frac{dq^i}{dt}
    =
    \frac{\partial H}{\partial p_i},
    \qquad
    \frac{dp_i}{dt}
    =
    -\frac{\partial H}{\partial q^i}
    +
    C_{ij}\frac{\partial H}{\partial p_j}
    +
    H\varphi_i .
\end{equation}
On the hypersurface \(H=0\), the conformal term \(H\varphi_i\) vanishes. Hence
the LCS Hamiltonian equations reduce on shell to the LCS evolution equations
associated with \eqref{evo-dyn-LCS}, namely
\begin{equation}\label{LCS-evo-local}
    \frac{dq^i}{dt}
    =
    \frac{\partial H}{\partial p_i},
    \qquad
    \frac{dp_i}{dt}
    =
    -\frac{\partial H}{\partial q^i}
    +
    C_{ij}\frac{\partial H}{\partial p_j}.
\end{equation}

\subsection{Continuum Dynamics on the LCS Setting}\label{Sec-Sym-Man}

In this subsection, we briefly recall the Lie--Poisson formulation of kinetic
dynamics on a locally conformally symplectic manifold, following
\cite{esen2023conformal}. We then analyze this dynamics through its continuity formulation.

Let $(M,\Omega\,;\varphi)$ 
be a \(2n\)-dimensional LCS manifold, and let
\(\mathfrak{X}_{\mathrm{LCS-ham}}(M)\) denote the space of LCS Hamiltonian
vector fields introduced in \eqref{LCS-Ham-vfield}. 
On the function space \(\mathcal{F}(M)\), the LCS Jacobi bracket
\eqref{LCS-Bracket} defines a Lie algebra structure. In light of the identity \eqref{sennur}, the assignment
\begin{equation}\label{Phi--iso}
\Phi:
\mathcal{F}(M)
\longrightarrow
\mathfrak{X}_{\mathrm{LCS-ham}}(M),
\qquad
H
\longmapsto
\mathcal{E}_H
\end{equation}
is a surjective Lie algebra homomorphism when the target is
equipped with the negative Jacobi--Lie bracket $[X,Y]_{-JL}=-[X,Y]$. Its kernel is the Lichnerowicz--de~Rham closed functions \eqref{twisted-closed-I}:
\begin{equation}
\ker\Phi
=
\bigl\{
H\in\mathcal{F}(M)
\; ; \;
d_\varphi H=0
\bigr\}.
\end{equation}
Hence,
\begin{equation}
\mathcal{F}(M)/\ker\Phi
\simeq
\mathfrak{X}_{\mathrm{LCS-ham}}(M)
\end{equation}
as Lie algebras.  The adjoint action on the function
representation is
\begin{equation}\label{Lie-brr-F}
\operatorname{ad}_H F
=
\{H,F\}^{\mathrm{LCS}}.
\end{equation}

We now pass to the dual description. The dual space of
\(\mathcal{F}(M)\) is identified with the space of densities
\(\operatorname{Den}(M)\) through the \(L_2\)-pairing
\begin{equation}\label{L-2-} 
\langle\bullet,\bullet\rangle_{L_2}:
\operatorname{Den}(M)\times\mathcal{F}(M)
 \longrightarrow
\mathbb{R}, \qquad
(f\,\mathrm{Vol}_\Omega,H)
 \mapsto
\int_M Hf\,\mathrm{Vol}_\Omega . 
\end{equation}
We fix the LCS volume form $\mathrm{Vol}_\Omega
=
\Omega^n$. 
Accordingly, a density \(f\,\mathrm{Vol}_\Omega\) is represented only by its density function
\(f\). 
Referring to the identity \eqref{LC-assoc}, we compute 
\begin{equation}\label{seko--}
\begin{split}
\big\langle \operatorname{ad}^*_H f ,F \big\rangle & = \big\langle f , \operatorname{ad}_H F \big\rangle = \big\langle f , \{H,F\}^{\mathrm{LCS}} \big\rangle 
\\ & = \big\langle 
\{f,H\}^{\mathrm{LCS}}
+
(n+2)\langle\varphi,\mathcal{E}_H\rangle f, F \big\rangle.  
\end{split}
\end{equation}
Hence the coadjoint action of \(\mathcal{F}(M)\) on its dual is obtained as 
\begin{equation}\label{seko}
\operatorname{ad}^*_H f = \{f,H\}^{\mathrm{LCS}}
+
(n+2)\langle\varphi,\mathcal{E}_H\rangle f. 
\end{equation}

The dual space \(\mathcal{F}^*(M)\) carries the corresponding Lie--Poisson
bracket. For two functionals
\(\mathfrak{F},\mathfrak{H}:\mathcal{F}^*(M)\longrightarrow\mathbb{R}\),
we define
\begin{equation}\label{LCS-Lie-Poisson-bracket}
\left\{
\mathfrak{F},
\mathfrak{H}
\right\}_{\mathcal{F}^*(M)}(f)
=
\int_M
f
\left\{
\frac{\delta\mathfrak{F}}{\delta f},
\frac{\delta\mathfrak{H}}{\delta f}
\right\}^{\mathrm{LCS}}
\mathrm{Vol}_\Omega .
\end{equation}
Here,
\(\delta\mathfrak{F}/\delta f\) and
\(\delta\mathfrak{H}/\delta f\) denote the Fr\'echet derivatives of the
functionals, which are identified with elements of the Lie algebra
\(\mathcal{F}(M)\). 
Given a Hamiltonian functional \(\mathfrak{H}\), the corresponding
Lie--Poisson equation is
\begin{equation}\label{LP-f}
\frac{\partial f}{\partial t}
=
-
\operatorname{ad}^*_{\delta\mathfrak{H}/\delta f}f =- \{f,H\}^{\mathrm{LCS}}
-
(n+2)\langle\varphi,\mathcal{E}_H\rangle f.
\end{equation}
We now choose the linear Hamiltonian functional
\begin{equation}\label{Ham-funcional-2}
\mathfrak{H}(f)
=
\int_M Hf\,\mathrm{Vol}_\Omega,
\end{equation}
for which $\delta\mathfrak{H}/\delta f =H$. 
Substituting the coadjoint action \eqref{seko} into
\eqref{LP-f}, we obtain what we call locally conformal Vlasov equation:
\begin{equation}\label{LC-Vlasov}
\frac{\partial f}{\partial t}
+
\{f,H\}^{\mathrm{LCS}}
=
-
(n+2)
\langle\varphi,\mathcal{E}_H\rangle f .
\end{equation}

\noindent \textbf{Continuity Form of the LCS Vlasov Equation.}
Using the relation \eqref{LCS-bracket-vector-relation} between the LCS
bracket and the LCS Hamiltonian vector field, the locally conformal Vlasov
equation can be written in terms of the LCS Hamiltonian vector field
\(\mathcal E_H\) as
\begin{equation}\label{LCS-Vlasov-transport-form}
\frac{\partial f}{\partial t}
+
\mathcal E_H(f)
=
-(n+1)
\langle\varphi,\mathcal E_H\rangle f .
\end{equation}
Recalling the divergence formula $\operatorname{div}_{\Omega}(\mathcal E_H)
=
n\langle\varphi,\mathcal E_H\rangle$ in \eqref{LCS-Ham-divergence}  
where the divergence is computed with respect to the volume form
\(\operatorname{Vol}_{\Omega}=\Omega^n\), equation
\eqref{LCS-Vlasov-transport-form} becomes
\begin{equation}\label{LCS-Vlasov-continuity-source-scalar}
\frac{\partial f}{\partial t}
+
\mathcal E_H(f)
+
\operatorname{div}_{\Omega}(\mathcal E_H) f
=
-\langle\varphi,\mathcal E_H\rangle f .
\end{equation}
Equivalently, in terms of densities with respect to the ordinary LCS volume
form \(\operatorname{Vol}_{\Omega}=\Omega^n\), one obtains the continuity
equation with a conformal source term
\begin{equation}\label{LCS-Vlasov-continuity-source}
\left(
\partial _ t
+
\mathcal L_{\mathcal E_H}
\right)
(f\,\operatorname{Vol}_{\Omega})
=
-\langle\varphi,\mathcal E_H\rangle
f\,\operatorname{Vol}_{\Omega}.
\end{equation}
Thus, the ordinary density \(f\,\operatorname{Vol}_{\Omega}\) is not conserved
by the LCS Hamiltonian flow. The defect from conservation is precisely
measured by the Lee contribution
\(\langle\varphi,\mathcal E_H\rangle\). In this sense, the Lee form records
the conformal expansion or contraction of the phase-space volume along the
Hamiltonian motion.

For the following global weighted formulation, we assume that the
Lee form is exact and fix a positive conformal unit $\mathbb I$,
as explained in Appendix~\ref{sec-LdR}. There are two equivalent ways to remove this source term.
For this, we introduce the conformally weighted unit \(\mathbb I\) given in \eqref{twisted-closed-I}, 
characterized by
\begin{equation}
d_{\varphi}\mathbb I=0,
\qquad
d\mathbb I=\varphi\,\mathbb I .
\end{equation}
Then its directional derivative along \(\mathcal E_H\) satisfies
\begin{equation}
\mathcal E_H(\mathbb I)
=
\langle\varphi,\mathcal E_H\rangle \mathbb I .
\end{equation}
Accordingly, we define the conformally weighted volume form
\begin{equation}\label{LCS-weighted-volume}
\operatorname{Vol}_{\Omega,\varphi}
=
\mathbb I\,\operatorname{Vol}_{\Omega}
=
\mathbb I\,\Omega^n .
\end{equation}
Using the preceding identity together with the divergence formula, one obtains
\begin{equation}
\mathcal L_{\mathcal E_H}\operatorname{Vol}_{\Omega,\varphi}
=
(n+1)\langle\varphi,\mathcal E_H\rangle
\operatorname{Vol}_{\Omega,\varphi}.
\end{equation}
Hence \eqref{LCS-Vlasov-transport-form} is equivalently written as the
source-free continuity equation
\begin{equation}\label{LCS-Vlasov-continuity-weighted}
\left(
\partial _ t
+
\mathcal L_{\mathcal E_H}
\right)
\left(
f\,\operatorname{Vol}_{\Omega,\varphi}
\right)
=
0.
\end{equation}
Therefore, the quantity transported by the LCS Hamiltonian flow is naturally
measured with respect to the conformally weighted volume
\(\mathbb I\Omega^n\), rather than the ordinary volume \(\Omega^n\).

The second, equivalent, way is to keep the ordinary volume form
\(\operatorname{Vol}_{\Omega}\), but to absorb the conformal weight into the
density function. We define
\begin{equation}\label{physical-density-def}
\mathfrak f
=
\mathbb I f .
\end{equation}
Then, using
\(\mathcal E_H(\mathbb I)
=
\langle\varphi,\mathcal E_H\rangle\mathbb I\), we compute
\begin{equation}
\begin{split}
\frac{\partial \mathfrak f}{\partial t}
+
\mathcal E_H(\mathfrak f)
+
\operatorname{div}_{\Omega}(\mathcal E_H)\mathfrak f
&=
\mathbb I
\left(
\frac{\partial f}{\partial t}
+
\mathcal E_H(f)
\right)
+
f\,\mathcal E_H(\mathbb I)
+
\operatorname{div}_{\Omega}(\mathcal E_H)\mathbb I f
\\
&=
-(n+1)\langle\varphi,\mathcal E_H\rangle \mathbb I f
+
\langle\varphi,\mathcal E_H\rangle \mathbb I f
+
n\langle\varphi,\mathcal E_H\rangle \mathbb I f
=0 .
\end{split}
\end{equation}
Consequently, \(\mathfrak f\) satisfies the ordinary source-free continuity
equation
\begin{equation}\label{LCS-Vlasov-continuity-physical}
\left(
\partial _ t
+
\mathcal L_{\mathcal E_H}
\right)
(\mathfrak f\,\operatorname{Vol}_{\Omega})
=
0.
\end{equation}

Thus the two formulations are equivalent:
\begin{equation}
f\,\operatorname{Vol}_{\Omega,\varphi}
=
f\,\mathbb I\,\operatorname{Vol}_{\Omega}
=
\mathfrak f\,\operatorname{Vol}_{\Omega}.
\end{equation}
The function \(f\) is the density coefficient adapted to the conformal LCS
volume, whereas \(\mathfrak f=\mathbb I f\) is the density coefficient with
respect to the ordinary volume \(\Omega^n\). In physical terms, the LCS
Hamiltonian flow does not preserve the ordinary phase-space volume; the Lee
form produces a conformal rescaling of volume elements. The apparent source
term in \eqref{LCS-Vlasov-continuity-source} is therefore not a creation or
annihilation of particles. Rather, it expresses the mismatch between the
ordinary volume \(\Omega^n\) and the natural conformally weighted volume
\(\mathbb I\Omega^n\). The conserved particle number over a material domain
\(D_t\) transported by \(\mathcal E_H\) is
\begin{equation}
\frac{d}{dt}
\int_{D_t}
f\,\operatorname{Vol}_{\Omega,\varphi}
=
0,
\end{equation}
or, equivalently,
\begin{equation}
\frac{d}{dt}
\int_{D_t}
\mathfrak f\,\operatorname{Vol}_{\Omega}
=
0.
\end{equation}

\noindent \textbf{Local Realizations.} Let \((q^i,p_i)\) be local coordinates on the cotangent bundle realization of
the LCS manifold, with the Lee form \(\varphi\) as in \eqref{aks} and the LCS
form \(\Omega\) as in \eqref{glo-glue-symp-2form}. In these coordinates, the
locally conformal Vlasov equation becomes
\begin{equation}\label{local-LCS-Vlasov}
\frac{\partial f}{\partial t}
+
\frac{\partial f}{\partial q^i}
\frac{\partial H}{\partial p_i}
-
\frac{\partial f}{\partial p_i}
\frac{\partial H}{\partial q^i}
=
-
H\varphi_i
\frac{\partial f}{\partial p_i}
-
C_{ij}
\frac{\partial f}{\partial p_i}
\frac{\partial H}{\partial p_j}
-
(n+1)
\varphi_i
\frac{\partial H}{\partial p_i}
f .
\end{equation}

Equivalently, the physical density $\mathfrak f=\mathbb I f$, defined in \eqref{physical-density-def}, has the local form 
\begin{equation}\label{local-LCS-Vlasov-mf}
\frac{\partial \mathfrak f}{\partial t}
+\frac{\partial \mathfrak f}{\partial q^i}\frac{\partial H}{\partial p_i}
-\frac{\partial \mathfrak f}{\partial p_i}\frac{\partial H}{\partial q^i}
=-H\varphi_i \frac{\partial \mathfrak f}{\partial p_i}
-C_{ij}\frac{\partial \mathfrak f}{\partial p_i}\frac{\partial H}{\partial p_j}
-n\varphi_i\frac{\partial H}{\partial p_i}\mathfrak f.
\end{equation}

The zero-energy hypersurface $H=0$ is invariant under the LCS Hamiltonian
flow, since
\begin{equation}
\mathcal L_{\mathcal E_H}H
=
H\,\varphi(\mathcal E_H).
\end{equation}
Restricting \eqref{local-LCS-Vlasov} to this hypersurface gives
\begin{equation}\label{local-LCS-Vlasov-on-shell}
\frac{\partial f}{\partial t}
+
\frac{\partial f}{\partial q^i}
\frac{\partial H}{\partial p_i}
-
\frac{\partial f}{\partial p_i}
\frac{\partial H}{\partial q^i}
=
-
C_{ij}
\frac{\partial f}{\partial p_i}
\frac{\partial H}{\partial p_j}
-
(n+1)
\varphi_i
\frac{\partial H}{\partial p_i}
f .
\end{equation}

If the Lee form vanishes identically, then the LCS structure reduces to the
symplectic one. Accordingly, the LCS bracket becomes the canonical Poisson
bracket, \(C_{ij}=0\), and the conformal source term disappears. In this case,
the continuity equation \eqref{LCS-Vlasov-continuity-weighted} reduces to the
ordinary symplectic continuity equation
\begin{equation}\label{Symplectic-continuity}
\left(
\partial_t
+
\mathcal L_{X_H}
\right)
\left(
f\,\operatorname{Vol}_{\Omega}
\right)
=
0,
\end{equation}
where \(\operatorname{Vol}_{\Omega}\) is now the symplectic volume. Equivalently,
the LCS--Vlasov equation \eqref{LC-Vlasov} reduces to the symplectic Vlasov
equation
\begin{equation}\label{Symplectic-Vlasov}
\frac{\partial f}{\partial t}
+
\{f,H\}^{(S)}
=
0.
\end{equation}
In canonical coordinates, this is simply
\begin{equation}\label{vlasovpos}
\frac{\partial f}{\partial t}
+
\frac{\partial f}{\partial q^i}
\frac{\partial H}{\partial p_i}
-
\frac{\partial f}{\partial p_i}
\frac{\partial H}{\partial q^i}
=
0.
\end{equation}
Thus, the ordinary symplectic Vlasov equation is recovered as the
zero-Lee-form limit.

\noindent \textbf{Entropy evolution in the LCS Vlasov equation.} Consider now the total entropy as volume integral of entropy density $\eta(\mathfrak f)$ with $\eta:\mathbb{R}\rightarrow\mathbb{R}$ being a smooth concave function satisfying $\eta(0)=0$. For instance, we could think of $\eta$ as $-\mathfrak{f} \ln \mathfrak{f}$ in the case of Boltzmann entropy \cite{Boltzmann}. How does the entropy evolve in the LCS setting?

The change of entropy in time reads
\begin{align}
    \frac{d S}{d t} =& \frac{d}{d t} \int_{D_t} \eta(\mathfrak{f}) \, \Vol_\Omega = \int_{D_t} \left[ \frac{\partial}{\partial t} \left( \eta(\mathfrak{f}) \, \Vol_\Omega \right) + \mathcal{L}_{\mathcal{E}_H} \left( \eta(\mathfrak{f}) \, \Vol_\Omega \right) \right]\nonumber\\
    =& \int_{D_t} \left[ \eta'(\mathfrak{f}) \left( \frac{\partial \mathfrak{f}}{\partial t} + \mathcal{E}_H(\mathfrak{f}) \right) + \eta(\mathfrak{f}) \operatorname{div}_\Omega(\mathcal{E}_H) \right] \Vol_\Omega,
\end{align}
using the Cartan formula.

From the evolution equation of $\mathfrak{f}$, Equation \eqref{LCS-Vlasov-continuity-physical}, it follows that 
\begin{align}\label{eq.entropy}
    \frac{d S}{d t} =& \int_{D_t} \left[ \eta'(\mathfrak{f}) \left( -\mathfrak{f} \operatorname{div}_\Omega(\mathcal{E}_H) \right) + \eta(\mathfrak{f}) \operatorname{div}_\Omega(\mathcal{E}_H) \right] \Vol_\Omega \nonumber \\
    =& \int_{D_t} \left( \eta(\mathfrak{f}) - \mathfrak{f} \eta'(\mathfrak{f}) \right) \operatorname{div}_\Omega(\mathcal{E}_H) \, \Vol_\Omega
    = n \int_{D_t} \left( \eta(\mathfrak{f}) - \mathfrak{f} \eta'(\mathfrak{f}) \right) \langle \varphi, \mathcal{E}_H \rangle \, \Vol_\Omega.
\end{align}
Since
\(\eta(\mathfrak f)-\mathfrak f\eta'(\mathfrak f)\geq0\)
for a concave entropy density satisfying \(\eta(0)=0\), the entropy is
non-decreasing if
\(\langle\varphi,\mathcal E_H\rangle\geq0\) throughout \(D_t\), and
non-increasing if
\(\langle\varphi,\mathcal E_H\rangle\leq0\) throughout \(D_t\).
If this quantity changes sign, the sign of the entropy production is
determined by its weighted phase-space average.

\subsection{Liouville Equation for a Thermostatted Ensemble}\label{sec-NH}

In this application, by a thermostat we mean the classical Nos\'e--Hoover
thermostat. We first show that the Nos\'e--Hoover equations fit naturally into
the LCS formalism as an LCS Hamiltonian dynamics restricted to the invariant
hypersurface \(H=0\). Then, following Subsection~\ref{Sec-Sym-Man}, we pass to
the kinetic description and present the Liouville equation for an ensemble of
thermostatted particles.

\noindent \textbf{Nos\'e--Hoover Thermostat as LCS Hamiltonian Dynamics.}
The classical Nos\'e--Hoover thermostat is a deterministic constant-temperature
model obtained from Nos\'e's extended-phase-space approach and Hoover's
time-scaling-free reformulation \cite{Hoover1985,Nose1984}. In Hoover's
formulation, the equations of motion may be written as
\begin{equation}\label{eq:NH-literature}
    \dot q^i
    =
    \frac{p_i}{m_i},
    \qquad
    \dot p_i
    =
    F_i(q)-y\,p_i,
    \qquad
    \dot y
    =
    \frac{1}{W}
    \left(
    \sum_i\frac{p_i^2}{m_i}-gkT
    \right),
\end{equation}
where \(F_i(q)=-\partial \phi/\partial q^i\), \(m_i\) denotes the particle
masses, \(W>0\) is the thermostat mass, \(g\) is the number of momentum degrees
of freedom, \(k\) is the Boltzmann constant, \(T\) is the prescribed
temperature, and \(y\) is the thermostat variable. We use the normalized form
obtained by choosing \(m_i=1\), \(W=1\), and denoting the constant combination
\(gkT\) by \(C\). Then the Nos\'e--Hoover system becomes
\begin{equation}\label{eq:NH-classical}
    \dot q^i
    =
    \delta^{ij}p_j,
    \qquad
    \dot p_i
    =
    -\frac{\partial\phi}{\partial q^i}
    -
    y\,p_i,
    \qquad
    \dot y
    =
    \delta^{ij}p_ip_j-C .
\end{equation}
Equivalently, by introducing an auxiliary variable \(x\) through \(\dot x=y\),
the system can be written as Newton's equations
\begin{equation}\label{eq:NH-second-order}
    \ddot q^i
    =
    -\delta^{ij}\frac{\partial\phi}{\partial q^j}
    -
    \dot x\,\dot q^i,
    \qquad
    \ddot x
    =
    \delta_{ij}\dot q^i\dot q^j-C .
\end{equation}
Accordingly, the force on the mechanical variables \(q^i\) splits into the
conservative part
\(F^i_{\mathrm{cons}}
=
-\delta^{ij}\partial\phi/\partial q^j\)
and the dissipative part
\(F^i_{\mathrm{diss}}=-\dot x\,\dot q^i\). 
The auxiliary variable $x$ is the time integral of the thermostat
variable $y$.
Thus, the Nos\'e--Hoover dynamics
may be regarded as a Newtonian system with a potential force together with a
velocity-dependent dissipative force governed by the auxiliary variable \(x\).

To realize this system geometrically, consider the extended base manifold
\(\bar Q=Q\times\mathbb R\), with local coordinates \((q^i,x)\), and its
cotangent bundle
\(T^*\bar Q=T^*Q\times\mathbb R^2\), with Darboux coordinates
\((q^i,p_i,x,y)\) on a cotangent chart. We choose the Lee one-form to be
\begin{equation}
    \varphi=-dx .
\end{equation}
Referring to \eqref{Omega-phi}, the corresponding LCS two-form is
\begin{equation}\label{Omega-bar-NH}
    \Omega_{\bar Q}^{\varphi}
    =
    -d_\varphi\Theta_{\bar Q}
    =
    -d\Theta_{\bar Q}
    +
    \varphi\wedge\Theta_{\bar Q}
    =
    dq^i\wedge dp_i
    +
    dx\wedge dy
    -
    p_i\,dx\wedge dq^i .
\end{equation}
Thus \((T^*\bar Q,\Omega_{\bar Q}^{\varphi};\varphi)\) is an LCS manifold. We
consider the Hamiltonian function
\begin{equation}\label{eq:NH-extended-H}
    H(q,p,x,y)
    =
    \frac12\,\delta^{ij}p_ip_j
    +
    \frac12\,y^2
    +
    \phi(q)
    +
    Cx .
\end{equation}
Then the LCS Hamiltonian equation \eqref{LCS-Ham-global}, or equivalently its
coordinate realization \eqref{LCS-Ham-local}, gives
\begin{equation}\label{eq:NH-twisted-system}
    \dot q^i
    =
    \delta^{ij}p_j,
    \qquad
    \dot p_i
    =
    -\frac{\partial\phi}{\partial q^i}
    -
    y\,p_i,
    \qquad
    \dot x
    =
    y,
    \qquad
    \dot y
    =
    \delta^{ij}p_ip_j-C-H .
\end{equation}
Thus the thermostat equation acquires the additional term \(-H\). This term
has a clear geometric meaning. By \eqref{cons-LCS-Ham}, namely
\(\mathcal L_{\mathcal E_H}H=H\varphi(\mathcal E_H)\), the Hamiltonian
function is not preserved in the LCS Hamiltonian formalism. Since
\(\varphi=-dx\) and \(\mathcal E_H(x)=y\), we obtain
\begin{equation}
    \mathcal L_{\mathcal E_H}H
    =
    -Hy .
\end{equation}
It follows that the hypersurface \(H=0\) is invariant under the flow of
\(\mathcal E_H\). If we restrict \eqref{eq:NH-twisted-system} to this invariant
level set, then the extra term \(-H\) disappears and we recover the normalized
Nos\'e--Hoover system \eqref{eq:NH-classical}. Hence the projected subsystem on
\((q^i,p_i,y)\) is exactly the classical Nos\'e--Hoover thermostat.

This also shows that the Nos\'e--Hoover system can be understood directly as
an LCS evolution dynamics. Indeed, the LCS evolution equation
\eqref{evo-dyn-LCS}, in its coordinate form \eqref{LCS-evo-local}, contains no
additional Hamiltonian term in the thermostat equation. Therefore, the
normalized Nos\'e--Hoover equations arise directly from the LCS evolution
dynamics on the full extended phase space, while the LCS Hamiltonian dynamics
recovers the same equations after restriction to the invariant hypersurface
\(H=0\). 

\noindent \textbf{Kinetic Dynamics of the Nos\'e--Hoover System.}
We now pass to the kinetic description of an ensemble of thermostatted
particles. Let \(f=f(q,p,x,y)\) be the density coefficient on
\(T^*\bar Q\). We assume that \(\dim Q=n\), so that
\(\dim \bar Q=n+1\). Since \(\varphi=-dx\) and
\(\mathcal E_H(x)=y\), one has
\(\langle\varphi,\mathcal E_H\rangle=-y\). Therefore, the LCS Vlasov equation
\eqref{LCS-Vlasov-transport-form} becomes, for the Nos\'e--Hoover Hamiltonian
system \eqref{eq:NH-twisted-system},
\begin{equation}\label{NH-LCS-Vlasov-f}
\begin{split}
    \frac{\partial f}{\partial t}
    +
    \delta^{ij}p_j\frac{\partial f}{\partial q^i}
    +
    \left(
    -\frac{\partial\phi}{\partial q^i}
    -
    y\,p_i
    \right)
    \frac{\partial f}{\partial p_i}
    +
    y\frac{\partial f}{\partial x}
    +
    \left(
    \delta^{ij}p_ip_j-C-H
    \right)
    \frac{\partial f}{\partial y}
    =
    (n+2)y f .
\end{split}
\end{equation}
Evaluating the ambient LCS Vlasov equation along the invariant hypersurface
\(H=0\) removes the explicit term involving \(H\) and gives
\begin{equation}\label{NH-LCS-Vlasov-f-onshell}
\begin{split}
    \frac{\partial f}{\partial t}
    +
    \delta^{ij}p_j\frac{\partial f}{\partial q^i}
    +
    \left(
    -\frac{\partial\phi}{\partial q^i}
    -
    y\,p_i
    \right)
    \frac{\partial f}{\partial p_i}
    +
    y\frac{\partial f}{\partial x}
    +
    \left(
    \delta^{ij}p_ip_j-C
    \right)
    \frac{\partial f}{\partial y}
    =
    (n+2)y f .
\end{split}
\end{equation}
Equivalently, let \(\mathfrak f=\mathbb I f\) be the density coefficient with
respect to the ordinary LCS volume, as in \eqref{physical-density-def}. Then
\(\mathfrak f\) satisfies the source-free continuity equation
\eqref{LCS-Vlasov-continuity-physical}. In the present coordinates this reads
\begin{equation}\label{NH-LCS-Vlasov-mf}
\begin{split}
    \frac{\partial \mathfrak f}{\partial t}
    +
    \delta^{ij}p_j\frac{\partial \mathfrak f}{\partial q^i}
    +
    \left(
    -\frac{\partial\phi}{\partial q^i}
    -
    y\,p_i
    \right)
    \frac{\partial \mathfrak f}{\partial p_i}
    +
    y\frac{\partial \mathfrak f}{\partial x}
    +
    \left(
    \delta^{ij}p_ip_j-C-H
    \right)
    \frac{\partial \mathfrak f}{\partial y}
    -
    (n+1)y\,\mathfrak f
    =
    0 .
\end{split}
\end{equation}
The coefficient \(n+1\) in \eqref{NH-LCS-Vlasov-mf} is the phase-space
compression rate of the full LCS Hamiltonian vector field. Indeed, the
friction term contributes \(-ny\), while the additional term \(-H\) in the
\(y\)-equation contributes
\[
-\frac{\partial H}{\partial y}
=
-y.
\]
Hence,
\[
\operatorname{div}_{\Omega}(\mathcal E_H)
=
-(n+1)y.
\]

On the invariant hypersurface \(H=0\), the particle equations reduce to the
normalized Nos\'e--Hoover equations. However, the Liouville equation of the
Nos\'e--Hoover vector field must be computed from the restricted particle
dynamics itself, rather than by merely substituting \(H=0\) into the ambient
continuity equation. The normalized Nos\'e--Hoover vector field is
\begin{equation}
\label{NH-vector-field}
\begin{split}
X_{\mathrm{NH}}
={}&
\delta^{ij}p_j\frac{\partial}{\partial q^i}
+
\left(
-\frac{\partial\phi}{\partial q^i}
-
yp_i
\right)
\frac{\partial}{\partial p_i}
+
y\frac{\partial}{\partial x}
+
\left(
\delta^{ij}p_ip_j-C
\right)
\frac{\partial}{\partial y}.
\end{split}
\end{equation}
Its divergence with respect to the ordinary phase-space volume is
\begin{equation}
\label{NH-divergence}
\operatorname{div}_{\Omega}(X_{\mathrm{NH}})
=
-ny.
\end{equation}
Therefore, the density \(\mathfrak f_{\mathrm{NH}}\) transported by the
normalized Nos\'e--Hoover flow satisfies the Liouville equation
\begin{equation}
\label{NH-Liouville-equation}
\begin{split}
    \frac{\partial \mathfrak f_{\mathrm{NH}}}{\partial t}
    &+
    \delta^{ij}p_j
    \frac{\partial\mathfrak f_{\mathrm{NH}}}{\partial q^i}
    +
    \left(
    -\frac{\partial\phi}{\partial q^i}
    -
    yp_i
    \right)
    \frac{\partial\mathfrak f_{\mathrm{NH}}}{\partial p_i}
\\
    &+
    y\frac{\partial\mathfrak f_{\mathrm{NH}}}{\partial x}
    +
    \left(
    \delta^{ij}p_ip_j-C
    \right)
    \frac{\partial\mathfrak f_{\mathrm{NH}}}{\partial y}
    -
    ny\,\mathfrak f_{\mathrm{NH}}
    =
    0.
\end{split}
\end{equation}
Thus, the factor \(n+1\) belongs to the full LCS Hamiltonian dynamics,
whereas the phase-space compression factor of the normalized Nos\'e--Hoover
flow is \(n\).

 \noindent \textbf{Evolution of entropy.}
Let \(\mathfrak f_{\mathrm{NH}}\) be the density satisfying
\eqref{NH-Liouville-equation}, and define
\begin{equation}
S(t)
=
\int_{D_t}
\eta(\mathfrak f_{\mathrm{NH}})
\operatorname{Vol}_{\Omega},
\end{equation}
where \(D_t\) is transported by \(X_{\mathrm{NH}}\). Using
\eqref{NH-divergence}, we obtain
\begin{equation}\label{NH-entropy-evolution}
\frac{dS}{dt}
=
-n
\int_{D_t}
y
\left(
\eta(\mathfrak f_{\mathrm{NH}})
-
\mathfrak f_{\mathrm{NH}}
\eta'(\mathfrak f_{\mathrm{NH}})
\right)
\operatorname{Vol}_{\Omega}.
\end{equation}
For a concave entropy density satisfying \(\eta(0)=0\), the factor
\[
\eta(\mathfrak f_{\mathrm{NH}})
-
\mathfrak f_{\mathrm{NH}}
\eta'(\mathfrak f_{\mathrm{NH}})
\]
is non-negative. Consequently, the sign of the entropy production is
determined by a weighted phase-space average of the thermostat variable
\(y\), rather than by its pointwise value alone. In particular, the entropy
is non-increasing if \(y\geq0\) on the support of the distribution and
non-decreasing if \(y\leq0\) there.

For the Boltzmann entropy density
\(\eta(\mathfrak f_{\mathrm{NH}})
=-\mathfrak f_{\mathrm{NH}}\log\mathfrak f_{\mathrm{NH}}\),
equation \eqref{NH-entropy-evolution} reduces to
\begin{equation}
\frac{dS}{dt}
=
-n
\int_{D_t}
y\,\mathfrak f_{\mathrm{NH}}
\operatorname{Vol}_{\Omega}.
\end{equation}

\subsection{LCS-driven Vlasov Dynamics with a Momentum-dependent Gyroscopic Term}\label{sec-LCS-MV}

In this application, we first recall the Hamiltonian realization of the
Lorentz force law through a magnetic symplectic structure, which provides the
standard geometric basis for the Maxwell--Vlasov dynamics. We then replace the
magnetic symplectic structure by a locally conformally symplectic structure,
thereby obtaining a Lorentz-type law driven by the LCS geometry. Physically, the skew-symmetric tensor induced by the Lee form occupies,
in the equations of motion, the position held by the magnetic field in
the Lorentz force law. The analogy stops there: the resulting field
depends on the momentum and is thus a phase-space object rather than a
magnetic field on the configuration space, and the underlying two-form
is not closed. We therefore call this contribution gyroscopic, since it
is orthogonal to the momentum and performs no work. Passing from the particle picture to
an ensemble of particles governed by this LCS--driven Lorentz-type law, we
derive the associated kinetic continuum equation and refer to the resulting
model as an LCS-driven Vlasov model with a momentum-dependent gyroscopic term.

\noindent \textbf{Lorentz Law from a Magnetic Symplectic Form.}
Let \(Q\) be a three-dimensional configuration manifold. Consider the
cotangent bundle \(T^*Q\) with induced canonical coordinates
\((\boldsymbol q,\boldsymbol p)\) on a cotangent chart
\(U_\alpha\subset T^*Q\). Let
\(A=\boldsymbol A\cdot d\boldsymbol q\) be a local magnetic vector
potential. We define the magnetic two-form by its exterior derivative:
\begin{equation}
    B=dA, \qquad
    \boldsymbol B=\nabla\times\boldsymbol A,
\end{equation}
where \(B=(1/2)B_{ij}dq^i\wedge dq^j\) and
\(B_{ij}=\epsilon_{ijk}B^k\). The associated magnetic symplectic form is
\begin{equation}\label{magma}
    \Omega_B\vert_\alpha
    =
    dq^i\wedge dp_i
    -
    \frac12 B_{ij}(q)\,dq^i\wedge dq^j .
\end{equation}
Since \(B=dA\), one has \(dB=0\). Hence \(d\Omega_B=0\), and
\(\Omega_B\) is symplectic. We take the mechanical Hamiltonian
\begin{equation}\label{mech-ham}
    H(q,p)
    =
    \frac12 |\boldsymbol p|^2+V(q),
\end{equation}
where \(V\) is the electric potential. The Hamiltonian vector field \(X_H\),
determined by \(\iota_{X_H}\Omega_B=dH\), gives
\begin{equation}
    \dot{\boldsymbol q}
    =
    \boldsymbol p,
    \qquad
    \dot{\boldsymbol p}
    =
    -\nabla_{\boldsymbol q}V
    +
    \boldsymbol p\times\boldsymbol B .
\end{equation}
Since \(\dot{\boldsymbol q}=\boldsymbol p\), the Newton equation is
\begin{equation}\label{aassa}
    \ddot{\boldsymbol q}
    =
    -\nabla_{\boldsymbol q}V
    +
    \dot{\boldsymbol q}\times\boldsymbol B .
\end{equation}
Equivalently, this is the Lorentz force law
\begin{equation}
    \ddot{\boldsymbol q}
    =
    \boldsymbol E
    +
    \dot{\boldsymbol q}\times\boldsymbol B,
    \qquad
    \boldsymbol E=-\nabla_{\boldsymbol q}V,
    \qquad
    \boldsymbol B=\nabla\times\boldsymbol A .
\end{equation}
Thus, if \(\boldsymbol F=\ddot{\boldsymbol q}\) denotes the total force, then
\(\boldsymbol F_{\mathrm{el}}=\boldsymbol E\) and
\(\boldsymbol F_{\mathrm{mag}}=\dot{\boldsymbol q}\times\boldsymbol B\).
The magnetic contribution is therefore carried by the symplectic form, whereas
the electric contribution is carried by the Hamiltonian potential.

This geometric construction also explains the homogeneous Maxwell relations
in the static setting. Since
\(\boldsymbol B=\nabla\times\boldsymbol A\) and
\(\boldsymbol E=-\nabla_{\boldsymbol q}V\), one has
\begin{equation}
    \nabla\cdot\boldsymbol B=0,
    \qquad
    \nabla\times\boldsymbol E=\boldsymbol 0,
    \qquad
    \nabla\cdot\boldsymbol E=\rho,
    \qquad
    \nabla\times\boldsymbol B=\boldsymbol J .
\end{equation}
The first two identities correspond to the closedness of the magnetic
two-form and to the exactness of the electric one-form in the time-independent
case. The remaining equations, up to conventional physical constants and, in
the non-static case, the displacement-current term, are not determined by the
particle Hamiltonian system alone. They require additional field equations, or
equivalently a prescription of the charge and current sources. Hence, in the
present example, the symplectic geometry accounts for the Lorentz force law
and for the closedness of the magnetic two-form, while the full Maxwell system
belongs to the coupled field theory.

\noindent \textbf{Symplectic Continuity Equation and the Vlasov Equation.}
We now pass from the single-particle Lorentz dynamics to the corresponding
kinetic description. Let \(f=f(\boldsymbol q,\boldsymbol p)\) be a distribution
function on the magnetic phase space \((T^*Q,\Omega_B)\). Its evolution is
determined by transport along the Hamiltonian vector field \(X_H\), or
equivalently by the symplectic continuity equation
\begin{equation}
    \frac{\partial f}{\partial t}
    +
    \operatorname{div}_{\Omega_B}(fX_H)
    =
    0,
\end{equation}
where the divergence is computed with respect to the Liouville volume induced
by \(\Omega_B\). Since \(X_H\) is Hamiltonian, it preserves this volume, and
hence
\begin{equation}
    \frac{\partial f}{\partial t}
    +
    \{f,H\}^{(S)}
    =
    0.
\end{equation}
Using the Lorentz Hamiltonian vector field above, this gives the Vlasov
equation
\begin{equation}\label{magnetic-Vlasov-eq}
    \frac{\partial f}{\partial t}
    +
    \boldsymbol p\cdot\nabla_{\boldsymbol q} f
    +
    \big(
        \boldsymbol E
        +
        \boldsymbol p\times\boldsymbol B
    \big)\cdot\nabla_{\boldsymbol p} f
    =
    0 .
\end{equation}
Thus, the Vlasov equation is the continuity equation on the magnetic
symplectic phase space written in terms of the Lorentz force field. To obtain
the Maxwell--Vlasov system, this kinetic equation is coupled to Maxwell's
equations by taking the charge and current densities as moments of \(f\):
\begin{equation}\label{evee}
    \rho(\boldsymbol q)
    =
    \int_{\mathbb R^3} f(\boldsymbol q,\boldsymbol p)\,d^3p,
    \qquad
    \boldsymbol J(\boldsymbol q)
    =
    \int_{\mathbb R^3}
    \boldsymbol p\,
    f(\boldsymbol q,\boldsymbol p)\,d^3p .
\end{equation}
With normalized physical constants, the Maxwell--Vlasov dynamics is therefore
the coupled system consisting of \eqref{magnetic-Vlasov-eq} together with
Maxwell's equations for \(\boldsymbol E\) and \(\boldsymbol B\), with sources
\(\rho\) and \(\boldsymbol J\). In this way, the magnetic symplectic form
encodes the Lorentz part of the particle motion, while the electromagnetic
fields are determined self-consistently through their coupling to the
distribution function.

\noindent \textbf{LCS-driven Gyroscopic Dynamics.}
We now repeat the passage from particle motion to continuum dynamics, but this
time we replace the magnetic symplectic two-form \(\Omega_B\) in
\eqref{magma} by the LCS two-form
\begin{equation}
\begin{split}
   & \Omega = \Omega_Q + \varphi \wedge \Theta_Q
    =
    dq^i\wedge dp_i
    -
    \frac12 C_{ij}(q,p)\,dq^i\wedge dq^j,
    \\ 
   &  \qquad   C_{ij}(q,p)
    =
    \varphi_j(q)p_i-\varphi_i(q)p_j
\end{split}
\end{equation}
appearing in \eqref{glo-glue-symp-2form}, where $\Omega_Q = dq^i \wedge dp_i$ is the standard symplectic two-form on $T^*Q$ and $\Theta_Q$ is the canonical one-form. Note that this is compatible with the LCS definition \eqref{Omega-varphi}, as $d \varphi =0 $.

Consider the cotangent bundle \(T^*Q\) equipped with the LCS structure
\((\Omega,\varphi)\). In three
dimensions, both the Lee form \(\varphi\) and the skew-symmetric tensor
\(C_{ij}\) may be represented by the vectors \(\boldsymbol\varphi\) and
\(\boldsymbol C\), respectively, where
\begin{equation}
    C_{ij}
    =
    \epsilon_{ijk}C^k,
    \qquad
    \boldsymbol C
    =
    \boldsymbol p\times\boldsymbol\varphi .
\end{equation}
We consider the mechanical Hamiltonian exactly as in \eqref{mech-ham}. The
LCS Hamiltonian vector field \(\mathcal E_H\) is defined by the twisted
Hamiltonian equation
\(\iota_{\mathcal E_H}\Omega=d_\varphi H\). Using \eqref{LCS-Ham-local}, we
obtain, in vector notation,
\begin{equation}\label{doguu}
    \dot{\boldsymbol q}
    =
    \boldsymbol p,
    \qquad
    \dot{\boldsymbol p}
    =
    -\nabla V
    +
    \boldsymbol p\times
    (\boldsymbol p\times\boldsymbol\varphi)
    +
    H\boldsymbol\varphi .
\end{equation}
Since \(\dot{\boldsymbol q}=\boldsymbol p\), the Newton form of the LCS
Hamiltonian dynamics is
\begin{equation}\label{C-driven-Ham-Newton}
    \ddot{\boldsymbol q}
    =
    \boldsymbol E
    +
    \dot{\boldsymbol q}\times\boldsymbol B_C^{\mathrm{eff}}
    +
    H\boldsymbol\varphi,
    \qquad
    \boldsymbol E=-\nabla V,
    \qquad
    \boldsymbol B_C^{\mathrm{eff}}
    =
    \dot{\boldsymbol q}\times\boldsymbol\varphi .
\end{equation}
Thus the ordinary magnetic field is replaced by the effective LCS-driven field
\(\boldsymbol B_C^{\mathrm{eff}}\), which depends on the velocity and on the
Lee form. Equivalently, if
\(\boldsymbol F^{\mathrm{Ham}}=\ddot{\boldsymbol q}\) denotes the total force,
then the electric part is
\(\boldsymbol F_{\mathrm{el}}=\boldsymbol E\), the LCS-driven gyroscopic part is
\begin{equation}
    \boldsymbol F_C^{\mathrm{Ham}}
    =
    \dot{\boldsymbol q}\times\boldsymbol B_C^{\mathrm{eff}}
    =
    \dot{\boldsymbol q}\times
    (\dot{\boldsymbol q}\times\boldsymbol\varphi)
    =
    \langle\boldsymbol\varphi,\dot{\boldsymbol q}\rangle\dot{\boldsymbol q}
    -
    |\dot{\boldsymbol q}|^2\boldsymbol\varphi ,
\end{equation}
and the conformal part is
\(\boldsymbol F_{\mathrm{conf}}=H\boldsymbol\varphi\). Thus
\(\boldsymbol F^{\mathrm{Ham}}
=
\boldsymbol F_{\mathrm{el}}
+
\boldsymbol F_C^{\mathrm{Ham}}
+
\boldsymbol F_{\mathrm{conf}}\).
Moreover,
\begin{equation}
    \left\langle
    \dot{\boldsymbol q},
    \boldsymbol F_C^{\mathrm{Ham}}
    \right\rangle
    =
    0.
\end{equation}
Hence this contribution is gyroscopic: it is orthogonal to the velocity,
performs no work on the particle and does not change its kinetic energy.
It occupies in \eqref{C-driven-Ham-Newton} the position held by the
magnetic force in the Lorentz law; the two should nevertheless not be
identified, since the effective field depends on the momentum.
 
On the zero-energy hypersurface \(H=0\), the conformal force
\(H\boldsymbol\varphi\) vanishes. Therefore, on shell, the LCS Hamiltonian
dynamics reduces to the corresponding LCS evolution dynamics, and
\eqref{C-driven-Ham-Newton} becomes the Lorentz-type equation and the LCS-driven gyroscopic contribution.

It is worth noting that the effective momentum-dependent phase-space field
\[
\mathbf B_C^{\mathrm{eff}}
=
\dot{\mathbf q}\times\boldsymbol\varphi
\]
and the resulting equations of motion depend only on the Lee form \(\varphi\), rather than on the particular choice of its local potentials. Indeed, replacing a local Lee potential by \(\sigma_\alpha\mapsto\sigma_\alpha+c_\alpha\), where \(c_\alpha\) is locally constant, leaves \(\varphi|_\alpha=d\sigma_\alpha\) unchanged. Hence the local expressions of the LCS-driven force law are compatible on chart overlaps and determine a globally well-defined dynamics. This independence from additive constants in the local Lee potentials should not be confused with electromagnetic gauge invariance.

\noindent \textbf{LCS--driven Continuum Dynamics.}
In the classical Maxwell--Vlasov theory, the Vlasov equation is obtained by
transporting a one-particle distribution function along the Hamiltonian
particle dynamics induced by the Lorentz force law. In the present setting, we
perform the same construction after replacing the magnetic symplectic
structure by the locally conformally symplectic structure. 

Let \(f=f(q,p)\) be the density coefficient on \(T^*Q\). We consider the LCS
Hamiltonian vector field \(\mathcal E_H\) associated with the mechanical
Hamiltonian \(H=H(q,p)\) in \eqref{mech-ham}. As stated in
Subsection~\ref{Sec-Sym-Man}, the LCS Vlasov equation is equivalently the
continuity equation \eqref{LCS-Vlasov-continuity-source} with the conformal
source term. Using the local expression \eqref{doguu} of \(\mathcal E_H\), this
gives
\begin{equation}\label{LCS-driven-Vlasov-f}
\begin{split}
    \frac{\partial f}{\partial t}
    +
    \boldsymbol p\cdot\nabla_{\boldsymbol q}f
    +
    \left(
    \boldsymbol E
    +
    \langle\boldsymbol\varphi,\boldsymbol p\rangle\boldsymbol p
    -
    |\boldsymbol p|^2\boldsymbol\varphi
    +
    H\boldsymbol\varphi
    \right)
    \cdot
    \nabla_{\boldsymbol p}f
    =
    -
    4
    \langle\boldsymbol\varphi,\boldsymbol p\rangle f .
\end{split}
\end{equation}
Here \(\boldsymbol E=-\nabla_{\boldsymbol q}V\), and the term
\(\langle\boldsymbol\varphi,\boldsymbol p\rangle\boldsymbol p
-
|\boldsymbol p|^2\boldsymbol\varphi\) is the LCS-driven gyroscopic force,
namely
\(\boldsymbol p\times(\boldsymbol p\times\boldsymbol\varphi)\).

Let \(\mathfrak f=\mathbb I f\) be the density coefficient defined in
\eqref{physical-density-def} with respect to the ordinary LCS volume. Then
\(\mathfrak f\) satisfies the source-free continuity equation
\eqref{LCS-Vlasov-continuity-physical}. In local coordinates, this becomes
\begin{equation}\label{LCS-driven-Vlasov-mf}
\begin{split}
    \frac{\partial \mathfrak f}{\partial t}
    +
    \boldsymbol p\cdot\nabla_{\boldsymbol q}\mathfrak f
    +
    \left(
    \boldsymbol E
    +
    \langle\boldsymbol\varphi,\boldsymbol p\rangle\boldsymbol p
    -
    |\boldsymbol p|^2\boldsymbol\varphi
    +
    H\boldsymbol\varphi
    \right)
    \cdot
    \nabla_{\boldsymbol p}\mathfrak f
    +
    3
    \langle\boldsymbol\varphi,\boldsymbol p\rangle
    \mathfrak f
    =
    0 .
\end{split}
\end{equation}
We regard \eqref{LCS-driven-Vlasov-mf} as the full LCS-driven Vlasov model
with a momentum-dependent gyroscopic term. It describes the transport of the
weighted density \(\mathfrak f\) on the full phase space \(T^*Q\). The
zero-energy restriction discussed below is not taken as the definition of the
kinetic model.

\begin{remark}[Restriction to the invariant zero-energy hypersurface]
The zero-energy hypersurface \(H=0\) is invariant under the LCS Hamiltonian
flow, since
\[
\mathcal L_{\mathcal E_H}H
=
H\langle\varphi,\mathcal E_H\rangle.
\]
Consequently, the full phase-space equation
\eqref{LCS-driven-Vlasov-mf} may be evaluated along this invariant
hypersurface. On \(H=0\), the explicit conformal force
\(H\boldsymbol\varphi\) vanishes, and the restricted equation becomes
\begin{equation}
\label{LCS-driven-Vlasov-mf-zero-energy-restriction}
\begin{split}
    \frac{\partial \mathfrak f}{\partial t}
    +
    \boldsymbol p\cdot\nabla_{\boldsymbol q}\mathfrak f
    +
    \left(
    \boldsymbol E
    +
    \langle\boldsymbol\varphi,\boldsymbol p\rangle\boldsymbol p
    -
    |\boldsymbol p|^2\boldsymbol\varphi
    \right)
    \cdot
    \nabla_{\boldsymbol p}\mathfrak f
    +
    3
    \langle\boldsymbol\varphi,\boldsymbol p\rangle
    \mathfrak f
    =
    0 .
\end{split}
\end{equation}
Equation
\eqref{LCS-driven-Vlasov-mf-zero-energy-restriction}
should be understood as the restriction of the ambient phase-space equation
to an invariant hypersurface, rather than as an independently derived
continuity equation for a density intrinsically supported on \(H=0\).
An intrinsic kinetic theory on the energy hypersurface would require an
appropriate induced hypersurface measure and a corresponding surface
continuity equation.
\end{remark}

The charge density \(\rho\) and the current density \(\boldsymbol J\) are taken
as the moments of the conserved density \(\mathfrak f\), namely
\begin{equation}\label{LCS-rho-J}
    \rho(\boldsymbol q)
    =
    \int_{\mathbb R^3}
    \mathfrak f(\boldsymbol q,\boldsymbol p)\,d^3p,
    \qquad
    \boldsymbol J(\boldsymbol q)
    =
    \int_{\mathbb R^3}
    \boldsymbol p\,
    \mathfrak f(\boldsymbol q,\boldsymbol p)\,d^3p .
\end{equation}
The effective LCS-driven magnetic-type field is
\[
\boldsymbol B_C
=
\boldsymbol p\times\boldsymbol\varphi(\boldsymbol q).
\]
It is therefore not an ordinary magnetic field on the configuration space,
but a phase-space field. Since
\(\boldsymbol E=-\nabla_{\boldsymbol q}V\), one has
\begin{equation}
\label{LCS-electric-curl-identity}
    \nabla_{\boldsymbol q}\times\boldsymbol E
    =
    \boldsymbol 0.
\end{equation}
If a self-consistent electrostatic Poisson coupling is imposed, one
additionally requires
\begin{equation}
\label{LCS-electrostatic-Poisson-relation}
    \nabla_{\boldsymbol q}\cdot\boldsymbol E
    =
    \rho.
\end{equation}

For comparison with Maxwell's equations, we record the differential
identities satisfied by \(\boldsymbol B_C\). Since the Lee form is closed,
and with \(\boldsymbol p\) held fixed, one obtains
\begin{equation}
\label{LCS-phase-space-q-identities}
    \nabla_{\boldsymbol q}\cdot\boldsymbol B_C
    =
    0,
    \qquad
    \nabla_{\boldsymbol q}\times\boldsymbol B_C
    =
    \boldsymbol p\,
    \nabla_{\boldsymbol q}\cdot\boldsymbol\varphi
    -
    (\boldsymbol p\cdot\nabla_{\boldsymbol q})
    \boldsymbol\varphi .
\end{equation}
Its momentum-space derivatives satisfy
\begin{equation}
\label{LCS-phase-space-p-identities}
    \nabla_{\boldsymbol p}\cdot\boldsymbol B_C
    =
    0,
    \qquad
    \nabla_{\boldsymbol p}\times\boldsymbol B_C
    =
    -2\boldsymbol\varphi .
\end{equation}

Thus, the kinetic evolution is governed by the full LCS Vlasov equation
\eqref{LCS-driven-Vlasov-mf}. The associated phase-space field
\(\boldsymbol B_C\) obeys the differential identities
\eqref{LCS-phase-space-q-identities} with respect to the configuration
variables and \eqref{LCS-phase-space-p-identities} with respect to the
momentum variables. These identities are structural consequences of the
definition
\[
\boldsymbol B_C
=
\boldsymbol p\times\boldsymbol\varphi
\]
and of the closedness of the Lee form. They should therefore not be
interpreted as independent Maxwell field equations or as evolution equations
for \(\boldsymbol B_C\).

The charge density and the current moment are defined by
\eqref{LCS-rho-J}. However, in the present formulation,
\(\boldsymbol B_C\) is determined by the prescribed Lee form rather than
being sourced by the current \(\boldsymbol J\). Accordingly, the resulting
system should be interpreted as an LCS-driven Vlasov model on a geometric
background, rather than as a self-consistent Maxwell--Vlasov system. In
contrast with the ordinary magnetic field of classical Maxwell--Vlasov
theory, the gyroscopic response is generated by the phase-space field
\(\boldsymbol B_C\), whose momentum dependence records the essential
geometric difference from the classical case.

 Physically, the Lee form plays the role of a geometric background which turns
the particle momentum into an effective magnetic-type response. Hence the
magnetic force is not produced by an external field on the configuration
space, but by the conformal geometry of phase space itself. The additional
terms in the Vlasov-type equation record both this momentum-dependent magnetic
response and the conformal expansion or contraction of the transported
phase-space density.
  
\subsection{A Thermostatted LCS-driven Gyroscopic Model}\label{sec:NH-LCS-EM} 
The preceding thermostat and gyroscopic models can be combined on the extended configuration manifold \(\bar Q=Q\times\mathbb R\), with coordinates $(\qq,x)$ and cotangent coordinates $(\qq,x,\pp,y)$. 
Let $\varphi$ be a closed one-form on $Q$, and introduce the closed Lee form
\begin{equation}
    \Phi=\varphi-dx
\end{equation}
on $\bar{Q}$. The corresponding LCS two-form is then
\begin{equation}
    \bar\Omega 
    = -d_\Phi\Theta_{\bar Q}
    = dq^i\wedge dp_i+dx\wedge dy+\Phi\wedge(p_i\,dq^i+y\,dx),
\end{equation}
where $\Theta_{\bar Q}=p_i\,dq^i+y\,dx$. Combining the Hamiltonians of the two preceding models, we take the extended Hamiltonian
\begin{equation}\label{eq:NH-LCS-EM-Hamiltonian}
    H = \frac{1}{2}|\boldsymbol p|^2 + \frac12y^2 + V(\boldsymbol q) + Cx.
\end{equation}

Consequently, the LCS Hamilton's equations \eqref{LCS-Ham-local} are
\begin{equation}\label{eq:NH-LCS-EM-particle}
\begin{split}
    \dot{\boldsymbol q}
    & =\boldsymbol p,
    \qquad
    \dot x=y, \\
    \dot{\boldsymbol p}
    & = \boldsymbol E -y\boldsymbol p +\boldsymbol p\times(\boldsymbol p\times\varphi) +(H-y^2)\boldsymbol\varphi, \\
    \dot y
    & = |\boldsymbol p|^2-C  +y\langle\boldsymbol\varphi,\boldsymbol p\rangle-H, 
    \qquad \boldsymbol E=-\nabla V.
\end{split}
\end{equation}
The thermostat contributes the dissipative force $-y\boldsymbol p$, whereas $\pp\times(\pp\times \varphi)$ is the LCS-driven gyroscopic contribution, which is orthogonal to $\pp$.

The Hamiltonian is not conserved in general. Equation \eqref{cons-LCS-Ham} tells that
\begin{equation}\label{eq:NH-LCS-EM-H-evolution}
    \frac{dH}{dt}
    = H\,\Phi(\mathcal E_H) = H\, H_{p_i} \Phi_i
    = H
    \left(
    \langle\boldsymbol\varphi,\boldsymbol p\rangle-y
    \right).
\end{equation}
Thus, the energy changes through the conformal factor $\langle\boldsymbol\varphi,\boldsymbol p\rangle-y$, while the gyroscopic term itself makes no direct contribution. 
In particular, the hypersurface \(H=0\) is invariant. When \(\varphi=0\), Equation \eqref{eq:NH-LCS-EM-particle} reduces to the extended Nosé--Hoover system \eqref{eq:NH-twisted-system}. On $H=0$, the combined model retains both the gyroscopic force
and the additional coupling terms $-y^2\boldsymbol\varphi$ and
$y\langle\boldsymbol\varphi,\boldsymbol p\rangle$.
The normalized Nos\'e--Hoover equations are recovered when
$\varphi=0$.
 
The kinetic description of the Nosé--Hoover thermostat with gyroscopic forces can be obtained by plugging Equations \eqref{eq:NH-LCS-EM-particle} into the kinetic equation \eqref{local-LCS-Vlasov-mf}. Here $n=\dim Q$; in the three-dimensional vector notation used
above, $n=3$, so the coefficient $n+1$ equals $4$. 
Since $\langle\Phi,\mathcal E_H\rangle =\langle\boldsymbol\varphi,\boldsymbol p\rangle-y$, we obtain that
\begin{equation}\label{eq:NH-LCS-EM-Vlasov}
\begin{split}
    \frac{\partial\mathfrak f}{\partial t}
    &+ \boldsymbol p\cdot\nabla_{\boldsymbol q}\mathfrak f
    +\left(
    \boldsymbol E -y\boldsymbol p +\boldsymbol p\times(\boldsymbol p\times\boldsymbol\varphi) +(H-y^2)\boldsymbol\varphi
    \right)\cdot\nabla_{\boldsymbol p}\mathfrak f\\
    &
    +y\frac{\partial\mathfrak f}{\partial x}
    +\left(
    |\boldsymbol p|^2-C +y\langle\boldsymbol\varphi,\boldsymbol p\rangle-H
    \right)\frac{\partial\mathfrak f}{\partial y}
    +(n+1) \left(\langle\boldsymbol\varphi,\boldsymbol p\rangle-y\right)\mathfrak f 
    =0.
\end{split}
\end{equation}
This equation combines the Nos\'e--Hoover phase-space compression with the conformal changes induced by the gyroscopic LCS background.

Applying the general entropy balance \eqref{eq.entropy} to the
combined model gives
\begin{equation}\label{eq:NH-LCS-EM-entropy-evolution}
\frac{dS}{dt}
=(n+1)\int_{D_t}
\left(\eta(\mathfrak f)-\mathfrak f\eta'(\mathfrak f)\right)
\left(\langle\boldsymbol\varphi,\boldsymbol p\rangle-y\right)
\operatorname{Vol}_{\bar\Omega}.
\end{equation}
The first factor in the integrand is non-negative. Thus, (concave) entropy is non-decreasing when $\langle\boldsymbol\varphi,\boldsymbol p\rangle-y\geq0$, and non-increasing when this quantity is
non-positive. For the Boltzmann entropy density
$\eta(\mathfrak f)=-\mathfrak f\log\mathfrak f$, this becomes
\begin{equation}
\frac{dS}{dt}
=(n+1)\int_{D_t}
\left(\langle\boldsymbol\varphi,\boldsymbol p\rangle-y\right)\mathfrak f\,\operatorname{Vol}_{\bar\Omega}.
\end{equation}

\section{Conformal Analysis of Cosymplectic Continuum Dynamics}\label{sec:3}

\subsection{Cosymplectic Dynamics}\label{Sec:CosDyn}

\noindent \textbf{(Almost) Cosymplectic Dynamics.}
Let \(M\) be a \((2n+1)\)-dimensional manifold. It is called an almost
cosymplectic manifold if it is equipped with a one-form \(\tau\) and a two-form
\(\Omega\) satisfying the non-degeneracy condition
\begin{equation}\label{non-deg-cos}
    \mathrm{Vol}_{\tau,\Omega}
    =
    \tau\wedge\Omega^n
    \neq 0 .
\end{equation}
We denote such a manifold by \((M,\tau,\Omega)\). The non-degeneracy condition
\eqref{non-deg-cos} induces the \(\tau\)-flat mapping
\begin{equation}\label{flat-map-LCC}
    \flat^{\tau}:\mathfrak{X}(M)\longrightarrow \Gamma^1(M),
    \qquad
    X\longmapsto \iota_X\Omega+(\iota_X\tau)\tau .
\end{equation}
Its inverse is the \(\tau\)-sharp mapping
\begin{equation}\label{sharp-tau}
    \sharp^{\tau}: \Gamma^1(M) \longrightarrow \mathfrak{X}(M),
    \qquad
    \sharp^{\tau}=(\flat^{\tau})^{-1}.
\end{equation}
In this geometry, there exists a distinguished vector field, called the
\(\tau\)-Reeb vector field and denoted by \(\mathcal R\), defined by
\begin{equation}\label{Reeb-cosymp}
    \iota_{\mathcal R}\tau=1,
    \qquad
    \iota_{\mathcal R}\Omega=0,
    \qquad
    \text{or equivalently}
    \qquad
    \mathcal R=\sharp^{\tau}\tau .
\end{equation}

A direct computation gives the following identity for the \(\tau\)-sharp
mapping:
\begin{equation}\label{sharp-dualization}
    \langle \alpha,\sharp^\tau\beta\rangle
    +
    \langle \beta,\sharp^\tau\alpha\rangle
    =
    2\langle \tau,\sharp^\tau\alpha\rangle
    \langle \tau,\sharp^\tau\beta\rangle .
\end{equation}
In particular, choosing \(\beta=\tau\), one obtains the duality correspondence
between the Reeb vector field and the one-form \(\tau\):
\begin{equation}\label{Reeb-sharp-correspondence}
    \langle \alpha,\mathcal R\rangle
    =
    \langle \tau,\sharp^\tau\alpha\rangle .
\end{equation}

We now introduce a bivector field \(\Lambda^\tau\) in terms of the
\(\tau\)-sharp mapping by
\begin{equation}\label{Lambda-cos}
    \Lambda^\tau(\alpha,\beta)
    =
    \Omega(\sharp^\tau \alpha, \sharp^\tau \beta).
\end{equation}
By right contraction, this bivector induces the musical mapping
\begin{equation}\label{sharp-cos-L}
    \sharp_{\Lambda^\tau}\beta
    =
    \iota_\beta \Lambda^\tau
    =
    \Lambda^\tau(\bullet,\beta).
\end{equation}
Unlike the \(\tau\)-sharp mapping \(\sharp^\tau\), the musical mapping
\(\sharp_{\Lambda^\tau}\) is not an isomorphism. Indeed, its one-dimensional
kernel is spanned by \(\tau\). The relation between the two musical mappings is
\begin{equation}\label{sharp-Lambda-sharp-tau}
    \sharp_{\Lambda^\tau}\beta
    =
    \sharp^\tau \beta
    -
    \langle \beta,\mathcal R \rangle \mathcal R .
\end{equation}
Thus \(\sharp_{\Lambda^\tau}\) is obtained from \(\sharp^\tau\) by projecting
out the Reeb direction.

Referring to the bivector field \(\Lambda^\tau\), we define a bracket on
smooth functions by
\begin{equation}\label{cos-bracket}
    \{F,H\}^{\tau}
    =
    \Lambda^\tau(dF,dH).
\end{equation}
This bracket is skew-symmetric and satisfies the Leibniz identity. In general,
however, it does not satisfy the Jacobi identity, and therefore defines only an
almost Poisson bracket. Moreover, it is degenerate, since the associated
musical mapping \(\sharp_{\Lambda^\tau}\) has kernel spanned by \(\tau\).

\noindent \textbf{Cosymplectic Manifold.}
An almost cosymplectic manifold \((M,\tau,\Omega)\) is called cosymplectic if
both defining forms are closed, that is,
\begin{equation}\label{closure-cosym}
    d\tau=0,
    \qquad
    d\Omega=0 .
\end{equation}
In this case, \((M,\tau,\Omega)\) is the natural odd-dimensional counterpart of
a symplectic manifold; see, for example, \cite{Cape}. The closedness of
\(\tau\) and \(\Omega\) turns the almost cosymplectic structure into a Poisson
framework. Indeed, when \eqref{closure-cosym} holds, the bivector field
\(\Lambda^\tau\) satisfies
\begin{equation}
    [\Lambda^\tau,\Lambda^\tau]=0 .
\end{equation}
Accordingly, the bracket \(\{\bullet,\bullet\}^{\tau}\) satisfies the Jacobi
identity and becomes a genuine Poisson bracket.

\noindent \textbf{Extended Cotangent Bundle.}
We consider the global coordinate \(t\) on \(\mathbb R\), interpreted as the
time variable. Pulling back the canonical symplectic two-form \(\Omega_Q\) on
\(T^*Q\) to the product manifold \(\mathbb R\times T^*Q\), we obtain the pair
\begin{equation}\label{cancosypair}
    \tau=dt,
    \qquad
    \Omega_Q=-d\Theta_Q,
\end{equation}
which defines a cosymplectic structure on \(\mathbb R\times T^*Q\). Here, to
avoid notational inflation, we use the same symbols for the canonical forms on
\(T^*Q\) and for their pullbacks to \(\mathbb R\times T^*Q\). In local
coordinates \((t,q^i,p_i)\), we have
\begin{equation}\label{cancosypair-coord}
    \tau|_\alpha=dt,
    \qquad
    \Theta_Q|_\alpha=p_i\,dq^i,
    \qquad
    \Omega_Q|_\alpha=dq^i\wedge dp_i,
    \qquad
    \mathcal R=\frac{\partial}{\partial t}.
\end{equation}

\noindent \textbf{Hamiltonian Dynamics.}
For a cosymplectic manifold \((M,\tau,\Omega)\), the Hamiltonian vector field
associated with a Hamiltonian function \(H\) is defined by
\begin{equation}\label{cos-Ham-musical}
    X_H
    =
    \sharp_{\Lambda^\tau}(dH)
    =
    \sharp^\tau(dH)
    -
    \mathcal R(H)\mathcal R .
\end{equation}
Equivalently, it is characterized by the contraction identities
\begin{equation}\label{cos-Ham-contractions}
    \iota_{X_H}\tau=0,
    \qquad
    \iota_{X_H}\Omega
    =
    dH-\mathcal R(H)\tau .
\end{equation}
Thus \(X_H\) is tangent to the codimension-one distribution \(\ker\tau\). In
Darboux coordinates \((t,q^i,p_i)\), one obtains the Hamilton's equation 
\begin{equation}\label{cos-three-vf-Darboux} 
    \dot t=0,
    \qquad  \dot q^i
    =
    \frac{\partial H}{\partial p_i},
    \qquad
    \dot p_i
    =
    -
    \frac{\partial H}{\partial q^i}.
\end{equation}

Since the cosymplectic bracket is skew-symmetric, the Hamiltonian function is
preserved along its own Hamiltonian flow:
\begin{equation}\label{cos-Ham-H-preservation}
    \mathcal L_{X_H}H
    =
    X_H(H)
    =
    \{H,H\}^{\tau}
    =
    0 .
\end{equation}
Moreover, using \(d\tau=0\), \(d\Omega=0\), and
\eqref{cos-Ham-contractions}, one obtains
\begin{equation}\label{cos-Ham-form-evolution}
    \mathcal L_{X_H}\tau=0,
    \qquad
    \mathcal L_{X_H}\Omega
    =
    -d\bigl(\mathcal R(H)\bigr)\wedge\tau .
\end{equation}
Although \(\Omega\) is not necessarily preserved, the natural cosymplectic
volume form is preserved:
\begin{equation}\label{cos-Ham-volume-preservation}
    \mathcal L_{X_H}\mathrm{Vol}_{\tau,\Omega}
    =
    0 .
\end{equation}

The family of cosymplectic Hamiltonian vector fields is closed under the
Jacobi--Lie bracket. More precisely,
\begin{equation}\label{cos-Ham-Lie-closure}
    [X_F,X_H]
    =
    -X_{\{F,H\}^{\tau}} .
\end{equation}
Consequently, the cosymplectic Hamiltonian vector fields form a Lie algebra,
and the characteristic distribution generated by them is involutive and hence
integrable.

\noindent \textbf{Evolution Vector Field.}
The evolution vector field associated with \(H\) is defined by
\begin{equation}\label{cos-evolution-musical}
    E_H
    =
    X_H+\mathcal R
    =
    \sharp_{\Lambda^\tau}(dH)+\mathcal R
    =
    \sharp^\tau(dH)
    +
    \bigl(1-\mathcal R(H)\bigr)\mathcal R .
\end{equation}
Equivalently, it is characterized by
\begin{equation}\label{cos-evolution-contractions}
    \iota_{E_H}\tau=1,
    \qquad
    \iota_{E_H}\Omega
    =
    dH-\mathcal R(H)\tau .
\end{equation}
The normalization \(\iota_{E_H}\tau=1\) distinguishes \(E_H\) from the
Hamiltonian vector field and makes it suitable for time-dependent dynamics. In
Darboux coordinates \((t,q^i,p_i)\), one obtains evolution equation
\begin{equation}\label{cos-evolution-vf-Darboux} 
    \dot t=1,
    \qquad
    \dot q^i
    =
    \frac{\partial H}{\partial p_i},
    \qquad
    \dot p_i
    =
    -
    \frac{\partial H}{\partial q^i}.
\end{equation}
The identity \(\dot t=1\) identifies the Darboux coordinate \(t\) with the
independent evolution parameter.

For the evolution dynamics, the Hamiltonian function evolves along \(E_H\) as
\begin{equation}\label{cos-evolution-H-evolution}
    \mathcal L_{E_H}H
    =
    E_H(H)
    =
    \mathcal R(H).
\end{equation}
Therefore, \(E_H\) preserves \(H\) if and only if \(\mathcal R(H)=0\). Thus the
Hamiltonian is conserved along the evolution flow precisely when it is constant
in the Reeb direction. The evolution vector field preserves the natural
cosymplectic volume form. Indeed,
\begin{equation}\label{cos-evolution-form-evolution}
    \mathcal L_{E_H}\tau=0,
    \qquad
    \mathcal L_{E_H}\Omega
    =
    -d\bigl(\mathcal R(H)\bigr)\wedge\tau ,
\end{equation}
and hence
\begin{equation}\label{cos-evolution-volume-preservation}
    \mathcal L_{E_H}\mathrm{Vol}_{\tau,\Omega}
    =
    0 .
\end{equation}
The evolution vector fields do not form a Lie algebra in general. 

\subsection{Locally Conformally Cosymplectic Manifolds} \label{sec:LCCosM}

In Subsection~\ref{sec:lcs}, we presented locally conformally symplectic geometry
from a global point of view and then justified the global expressions through
their local conformal realizations. In the present locally conformally
cosymplectic setting, we proceed in the opposite direction: we start from local
cosymplectic structures and glue them by conformal compatibility conditions.
The reason for this approach is that, in this paper, we aim to refine the
definition of a locally conformally cosymplectic manifold used in
\cite{AtEsLe23,ChLeMa91} by replacing it with a formulation that we find more
natural.

\noindent \textbf{From Local-to-Global.}
Let \(M\) be a \((2n+1)\)-dimensional manifold equipped with an atlas
\(\{U_\alpha\}\). We assume that every chart \(U_\alpha\) carries a local
cosymplectic structure \((\tau_\alpha,\Omega_\alpha)\), that is,
\begin{equation}\label{local-cosymplectic-conditions}
d\tau_\alpha=0,
\qquad
d\Omega_\alpha=0,
\qquad
\tau_\alpha\wedge\Omega_\alpha^n\neq 0 .
\end{equation}
As in the locally conformally symplectic construction, we do not assume that
the local cosymplectic structures glue directly to a global cosymplectic pair.
Instead, we assume that, on every non-empty overlap \(U_\alpha\cap U_\beta\),
they are conformally related by
\begin{equation}\label{lcc-loc-rel}
e^{\sigma_\alpha}\tau_\alpha
=
e^{\sigma_\beta}\tau_\beta,
\qquad
e^{\sigma_\alpha}\Omega_\alpha
=
e^{\sigma_\beta}\Omega_\beta,
\end{equation}
for a family of smooth local functions
\(\sigma_\alpha\in C^\infty(U_\alpha)\).
At this point, our construction departs from the convention commonly used in
the literature on locally conformally cosymplectic manifolds; see, for
example, \cite{AtEsLe23,ChLeMa91}. In the traditional formulation, the local
one-form and two-form are assigned different conformal weights\footnote{In the
traditional LCCos formulation, one has
\begin{equation*}\label{traditional-lcc-gluing}
\tau|_\alpha
=
e^{\sigma_\alpha}\tau_\alpha,
\qquad
\Omega|_\alpha
=
e^{2\sigma_\alpha}\Omega_\alpha .
\end{equation*}
Accordingly, the resulting global forms satisfy
$d\tau
=
\varphi\wedge\tau$ and 
$d\Omega
=
2\varphi\wedge\Omega$.}.
In the present paper, we instead assign the same conformal weight to
\(\tau_\alpha\) and \(\Omega_\alpha\), as in \eqref{lcc-loc-rel}. The relevance of this choice will become clearer in the discussion of the
corresponding Hamiltonian dynamics, where we comment further on this weight
choice.

Here, however, we proceed with the conformal compatibility condition
\eqref{lcc-loc-rel}. This condition implies that the exterior derivatives of
the conformal factors agree on overlaps:
\begin{equation}
d\sigma_\alpha
=
d\sigma_\beta
\qquad
\text{on } U_\alpha\cap U_\beta .
\end{equation}
Thus the local exact one-forms \(d\sigma_\alpha\) patch together to define a
global closed one-form \(\varphi\), called the Lee one-form, by
\begin{equation}\label{lcc-Lee-local}
\varphi|_\alpha
=
d\sigma_\alpha .
\end{equation}

Although the local cosymplectic data
\(\{(\tau_\alpha,\Omega_\alpha)\}\) do not, in general, glue directly to a
global cosymplectic pair, the conformal compatibility relations
\eqref{lcc-loc-rel} show that the rescaled local forms
\begin{equation}\label{local-forms-cosymplectic}
\tau|_\alpha
=
e^{\sigma_\alpha}\tau_\alpha,
\qquad
\Omega|_\alpha
=
e^{\sigma_\alpha}\Omega_\alpha
\end{equation}
patch together to define a global one-form \(\tau\) and a global two-form
\(\Omega\) on \(M\). The non-degeneracy of the resulting global pair follows
locally from
\begin{equation}
\left.
\tau\wedge\Omega^n
\right|_\alpha
=
e^{(n+1)\sigma_\alpha}
\tau_\alpha\wedge\Omega_\alpha^n
\neq 0.
\end{equation}
Moreover, using \eqref{local-cosymplectic-conditions},
\eqref{lcc-Lee-local}, and \eqref{local-forms-cosymplectic}, we compute
\begin{equation}
\begin{split}
d\tau|_\alpha
&=
d\bigl(e^{\sigma_\alpha}\tau_\alpha\bigr)
=
d\sigma_\alpha\wedge e^{\sigma_\alpha}\tau_\alpha
=
\varphi|_\alpha\wedge\tau|_\alpha,
\\
d\Omega|_\alpha
&=
d\bigl(e^{\sigma_\alpha}\Omega_\alpha\bigr)
=
d\sigma_\alpha\wedge e^{\sigma_\alpha}\Omega_\alpha
=
\varphi|_\alpha\wedge\Omega|_\alpha .
\end{split}
\end{equation}

In light of the above calculation, we shall call a manifold equipped with a
pair \((\tau,\Omega)\) and a closed one-form \(\varphi\) satisfying
\begin{equation}\label{new-lcc-structure-equations}
d\tau
=
\varphi\wedge\tau,
\qquad
d\Omega
=
\varphi\wedge\Omega,
\qquad
\text{or equivalently,}
\qquad
d_\varphi\tau=0,
\qquad
d_\varphi\Omega=0
\end{equation}
a locally conformally cosymplectic (LCCos) manifold. Thus, this definition
arises naturally from the conformal gluing of local cosymplectic structures.

Let \(\mathcal R_\alpha\) be the Reeb vector field of the local cosymplectic
structure \((U_\alpha,\tau_\alpha,\Omega_\alpha)\), that is,
\begin{equation}
\iota_{\mathcal R_\alpha}\tau_\alpha=1,
\qquad
\iota_{\mathcal R_\alpha}\Omega_\alpha=0 .
\end{equation}
If these local Reeb vector fields are glued by the conformal rule
\begin{equation}\label{globalReeb-LCCos}
\mathcal R|_\alpha
=
e^{-\sigma_\alpha}\mathcal R_\alpha ,
\end{equation}
then the resulting global vector field satisfies the Reeb conditions
\eqref{Reeb-cosymp} for the underlying almost cosymplectic pair
\((\tau,\Omega)\) of the LCCos manifold. Indeed, using
\(\tau|_\alpha=e^{\sigma_\alpha}\tau_\alpha\),
\(\Omega|_\alpha=e^{\sigma_\alpha}\Omega_\alpha\), and
\(\mathcal R|_\alpha=e^{-\sigma_\alpha}\mathcal R_\alpha\), one obtains
\begin{equation}\label{Reebok}
\iota_{\mathcal R}\tau
=
\iota_{e^{-\sigma_\alpha}\mathcal R_\alpha}
\bigl(e^{\sigma_\alpha}\tau_\alpha\bigr)
=
\iota_{\mathcal R_\alpha}\tau_\alpha
=
1,
\qquad
\iota_{\mathcal R}\Omega
=
\iota_{e^{-\sigma_\alpha}\mathcal R_\alpha}
\bigl(e^{\sigma_\alpha}\Omega_\alpha\bigr)
=
\iota_{\mathcal R_\alpha}\Omega_\alpha
=
0 .
\end{equation}
 
\noindent \textbf{Jacobi Structure.}
Referring to the local cosymplectic structure
\((U_\alpha,\tau_\alpha,\Omega_\alpha)\), we have the local Poisson bivector
field \(\Lambda_\alpha^\tau\) given by
\begin{equation}
    \Lambda_\alpha^\tau(dH_\alpha,dF_\alpha)=\{H_\alpha,F_\alpha\}_\alpha^\tau. 
\end{equation}
On overlapping charts \(U_\alpha\cap U_\beta\), these local Poisson bivector
fields satisfy the following conformal relation:
\begin{equation}
 e^{-\sigma_\alpha}\Lambda_\alpha^{\tau} =  e^{-\sigma_\beta}\Lambda_\beta^{\tau}.
\end{equation}
Notice that, this time, the conformal factor appears with a negative sign.
Accordingly, we obtain a global bivector field \(\Lambda^\tau\) whose local
realization is
\begin{equation}\label{LCCos-local-bivector}
\Lambda^\tau|_\alpha
=
e^{-\sigma_\alpha}\Lambda_\alpha^{\tau}.
\end{equation}
We now show that \(\Lambda^\tau\) is not Poisson in general, that is,
\([\Lambda^\tau,\Lambda^\tau]_{SN}\) is not identically zero. To see this, we
compute
\begin{equation}\label{hoplou}
\begin{split}
   [\Lambda^\tau_\alpha,\Lambda^\tau_\alpha]_{SN} &=[e^{\sigma_\alpha}\Lambda^\tau|_\alpha,e^{\sigma_\alpha}\Lambda^\tau|_\alpha]_{SN}    
   \\  &
   = e^{\sigma_\alpha}[e^{\sigma_\alpha}\Lambda^\tau|_\alpha,\Lambda^\tau|_\alpha]_{SN} +e^{2\sigma_\alpha} \iota_{d\sigma_\alpha} (\Lambda^\tau|_\alpha)\wedge \Lambda^\tau|_\alpha  
   \\  &
   = e^{\sigma_\alpha}[\Lambda^\tau|_\alpha,e^{\sigma_\alpha}\Lambda^\tau|_\alpha]_{SN} +e^{2\sigma_\alpha} \iota_{d\sigma_\alpha} (\Lambda^\tau|_\alpha)\wedge \Lambda^\tau|_\alpha 
   \\  & 
   = e^{\sigma_\alpha}\big(e^{\sigma_\alpha} [\Lambda^\tau|_\alpha,\Lambda^\tau|_\alpha]_{SN} + e^{\sigma_\alpha}  \iota_{d\sigma_\alpha} \Lambda^\tau|_\alpha \wedge \Lambda^\tau|_\alpha \big) + e^{2\sigma_\alpha} \iota_{d\sigma_\alpha} (\Lambda^\tau|_\alpha)\wedge \Lambda^\tau|_\alpha 
   \\  &
   = e^{2\sigma_\alpha} \big( [\Lambda^\tau|_\alpha,\Lambda^\tau|_\alpha]_{SN}+ 2 \iota_{d\sigma_\alpha} (\Lambda^\tau|_\alpha)\wedge \Lambda^\tau|_\alpha \big),
\end{split}
\end{equation}
where we have used the identity in \eqref{X-FY-2} and the commutation of
bivector fields in the Schouten--Nijenhuis bracket. Since
\(\Lambda_\alpha^\tau\) is a Poisson bivector field, the left-hand side is
identically zero, namely
\([\Lambda_\alpha^\tau,\Lambda_\alpha^\tau]_{SN}=0\). We define the vector
field associated with the Lee one-form by
\begin{equation}\label{LCCos-Lee-vector}
\mathcal{V}_{\varphi}
 =
\iota_{\varphi}\Lambda^\tau
=
\sharp_{\Lambda^\tau}(\varphi).
\end{equation}
Then, gluing these local identities, the calculation in \eqref{hoplou} gives
\begin{equation}
    [\Lambda^\tau,\Lambda^\tau]_{SN}=- 2 \iota_{\varphi} (\Lambda^\tau)\wedge \Lambda^\tau 
\end{equation}
Furthermore, we compute
\begin{equation}\label{LCCos-Schouten-second-calc}
\begin{split}
[\Lambda^\tau,\mathcal{V}_{\varphi}]_{SN} 
& =
[\Lambda^\tau,\iota_{\varphi}\Lambda^\tau]_{SN}
 = \frac{1}{2} \iota_{\varphi} [\Lambda^\tau,\Lambda^\tau]_{SN} \\&  =- \iota_{\varphi} \big(\iota_{ \varphi } (\Lambda^\tau)\wedge \Lambda^\tau \big) 
 \\& = - \iota_{\varphi} \iota_{\varphi}  (\Lambda^\tau)\wedge \Lambda^\tau - \big(\iota_{\varphi}  (\Lambda^\tau)\big)\wedge \big(\iota_{\varphi}  (\Lambda^\tau)\big) 
=0,
\end{split}
\end{equation}
where we have used the identities \eqref{iota-SN-2-2-X+} and
\eqref{iota-SN-2-2-X}, together with the facts that
\((\iota_{\varphi})^2=0\) and that the wedge product of a vector field with
itself vanishes. As a result,
\begin{equation}\label{LCCos-Jacobi-pair}
\bigl(\Lambda^\tau,-\mathcal{V}_{\varphi}\bigr)
\end{equation}
is a Jacobi pair.
Accordingly, the locally conformally cosymplectic bracket is computed to be
\begin{equation}\label{Pois-LCC}
\begin{split}
\{F,H\}^{\mathrm{LCCos}}
&=
\Lambda^\tau(d_{\varphi}F,d_{\varphi}H)
\\
&=
\Lambda^\tau(dF,dH)
+
F\,\mathcal{V}_{\varphi}(H)
-
H\,\mathcal{V}_{\varphi}(F).
\end{split}
\end{equation}
This has the form of the Jacobi bracket in \eqref{Jacobi-bracket} associated
with the Jacobi pair
\(\bigl(\Lambda^\tau,-\mathcal{V}_{\varphi}\bigr)\).

\noindent \textbf{Extended Cotangent Bundle as an LCCos Manifold.}
A natural source of locally conformally cosymplectic structures is the extended cotangent bundle \(\mathbb{R}\times T^*Q\) where $t$ is the global coordinate on $\mathbb{R}$. To construct such a structure, let \(\varphi\) be a closed one-form on the extended configuration space \(Q\times \mathbb{R}\), and pull it back to \(\mathbb{R}\times T^*Q\). For simplicity, we denote the pull-back again by \(\varphi\). We then define
\begin{equation}\label{LCC-can}
\tau =d_{\varphi}t=dt-t\varphi,
\qquad
\Omega=-d_{\varphi}\Theta_Q=\Omega_Q+\varphi\wedge\Theta_Q,
\end{equation}
where \(\Theta_Q\) is the canonical Liouville one-form on \(T^*Q\), and \(\Omega_Q=-d\Theta_Q\) is the canonical symplectic two-form. 
Here and below, we use the same notation for the canonical forms on \(T^*Q\)
and for their pullbacks to \(\mathbb{R}\times T^*Q\).  A direct calculation shows that the pair $(\tau,\Omega)$ defined in
\eqref{LCC-can} satisfies the LCCos structure equations
\eqref{new-lcc-structure-equations}. However, the non-degeneracy
condition must be imposed separately. Writing
$\varphi=\phi\,dt+\varphi_i\,dq^i$, we obtain
\begin{equation}
\tau\wedge\Omega^n
=
(1-t\phi)\,dt\wedge\Omega_Q^n.
\end{equation}
Consequently, the pair defines an LCCos structure on the open subset
\begin{equation}
U
=
\bigl\{
(t,q,p)\in\mathbb R\times T^*Q
:
1-t\phi(t,q)\neq0
\bigr\}.
\end{equation}
All subsequent coordinate expressions involving $(1-t\phi)^{-1}$
are understood on this open subset. In particular, when the Lee form
has no $dt$-component, so that $\phi=0$, the structure is
non-degenerate on the whole extended cotangent bundle.

Let \(U_{\alpha}\subset \mathbb{R}\times T^*Q\) be an open chart, and consider two Darboux coordinates
\begin{equation}\label{babanis}
(t,q^i,p_i),
\qquad
(s,q^i,r_i).
\end{equation}
The first one is Darboux coordinates for the canonical cosymplectic structure \((dt,\Omega_Q)\) on \(\mathbb{R}\times T^*Q\), as given in \eqref{cancosypair-coord}. The second Darboux coordinates \eqref{babanis} is the one for the local cosymplectic structure \((U_\alpha,\tau_\alpha,\Omega_\alpha) \). Hence we get 
\begin{equation}\label{babanis-2}
\tau_\alpha=ds,\qquad \Omega_\alpha= dq^i\wedge dr_i . 
\end{equation}
Since the Lee form is pulled back from the extended configuration
space, its coefficients depend only on $(t,q^i)$, and we write 
\begin{equation}\label{haaland}
\varphi\vert_\alpha=d\sigma_{\alpha}(t,q)= \phi\,dt + \varphi_i(t,q)\,dq^i=\frac{\partial \sigma_{\alpha}}{\partial t} dt + \frac{\partial \sigma_{\alpha}}{\partial q^i} dq^i.
\end{equation}
Then, the LCCos pair $(\tau,\Omega)$ in \eqref{LCC-can} on \(\mathbb{R}\times T^*Q\) written in terms of the Darboux coordinates $(t,q^i,p_i)$ as 
\begin{equation}\label{LCC-RTQ-forms}
\begin{split}
    \tau  &=(1-t\phi)dt-t\varphi_i\,dq^i,
\\ 
\Omega  &=dq^i\wedge dp_i-\frac{1}{2}C_{ij}\,dq^i\wedge dq^j + p_i\phi \, dt \wedge dq^i, \qquad C_{ij}(q,p)
    =
    \varphi_jp_i-\varphi_ip_j
\end{split}
\end{equation}
where the coefficients \(C_{ij}\) are of the same form as those defined in \eqref{Aij}.

\subsection{Locally Conformal Cosymplectic Dynamics}\label{sec:LCCosH}

In this subsection, we revisit the Hamiltonian and evolution dynamics
introduced in Subsection~\ref{Sec:CosDyn}, and formulate their locally
conformal cosymplectic counterparts. Thus, the cosymplectic Hamiltonian and
evolution vector fields are replaced by their LCCos analogues, obtained by
using the Lichnerowicz--de~Rham differential and the conformally modified
geometric data.

\noindent \textbf{Locally Conformal Cosymplectic Hamiltonian Dynamics.}
Local cosymplectic Hamiltonian dynamics can be glued under the conformal
compatibility introduced above. On each local cosymplectic chart
\((U_\alpha,\tau_\alpha,\Omega_\alpha)\), a local Hamiltonian function
\(H_\alpha\) determines the local Hamiltonian vector field by
\begin{equation}\label{localHam-LCCos}
    \iota_{X_{H_\alpha}}\tau_\alpha=0,
    \qquad
    \iota_{X_{H_\alpha}}\Omega_\alpha=dH_\alpha-\mathcal{R}_\alpha(H_\alpha)\tau_\alpha.
\end{equation}
We assume that these local Hamiltonian vector fields glue to a global
dynamical vector field \(\mathcal E_H\), namely,
\begin{equation}\label{globalHam-LCCos}
    \mathcal{E}_H|_\alpha=X_{H_\alpha}.
\end{equation}

The reason for the conformal weight choice in \eqref{lcc-loc-rel} becomes
clear at this point. In the formulations commonly used in the
literature\footnote{In the common convention in the literature, the
corresponding gluing rule is instead
\(\mathcal E_H|_\alpha=e^{-\sigma_\alpha}X_{H_\alpha}\). We do not adopt this
choice here, since our aim is to make the local Hamiltonian vector fields glue
directly to the global dynamical vector field.},
where the one-form and the two-form are assigned different conformal weights
\cite{AtEsLe23,ChLeMa91}, the local cosymplectic Hamiltonian vector fields do
not glue directly to a global vector field of the form
\eqref{globalHam-LCCos}. In other words, although the dynamics is
cosymplectic on each local chart, this property cannot be transferred
unchanged to the global level because of the conformal obstruction in the
gluing data. By contrast, the present choice assigns the same conformal weight
to \(\tau_\alpha\) and \(\Omega_\alpha\), and this makes the gluing of the
local Hamiltonian dynamics compatible with the global LCCos structure. This is
the same geometric principle that appears in the locally conformally
symplectic case through \eqref{LCS-localHamEq}: the locally conformal viewpoint
treats global dynamics as the result of gluing local dynamical systems whose
geometric data are conformally related. In this sense, our definition is not
merely a formal change in the structure equations; it is chosen so that the
local tensorial geometry and the corresponding local Hamiltonian equations
glue coherently. As we shall see below, the same perspective becomes even more
relevant for the evolution vector field, where the gluing procedure leads
naturally to a more intrinsic formulation of the dynamics. 

After this remark, we impose the corresponding compatibility condition on the
local Hamiltonian functions. Namely, we assume that
\begin{equation}\label{local-Hamiltonians-LCC}
    e^{\sigma_\alpha}H_\alpha=e^{\sigma_\beta}H_\beta
    \qquad
    \text{on } U_\alpha\cap U_\beta .
\end{equation}
Then the local Hamiltonians glue to a global function \(H\in C^\infty(M)\)
defined by
\begin{equation}\label{global-H-LCC}
    H|_\alpha=e^{\sigma_\alpha}H_\alpha .
\end{equation}
The conformal factor in \eqref{global-H-LCC} is essential: it matches the
conformal weight of the local two-forms \(\Omega_\alpha\) and produces the
twisted differential \(d_\varphi H\) in the global equation. Substituting the
gluing relations for \(\tau_\alpha\), \(\Omega_\alpha\), \(\mathcal R_\alpha\),
\(X_{H_\alpha}\), and \(H_\alpha\) into \eqref{localHam-LCCos}, 
one obtains the LCCos Hamiltonian dynamics
\begin{equation}\label{LCCos-Ham-global}
    \iota_{\mathcal{E}_H}\tau=0,
    \qquad
    \iota_{\mathcal{E}_H}\Omega
    =
    d_{\varphi}H-\langle d_{\varphi}H,\mathcal{R}\rangle\tau .
\end{equation}
Thus, \eqref{LCCos-Ham-global} is precisely the global form of the local
cosymplectic Hamiltonian equations \eqref{localHam-LCCos} under the conformal
gluing described above. When the Lee one-form vanishes, the conformal
contribution disappears and \eqref{LCCos-Ham-global} reduces to the usual
cosymplectic Hamiltonian formalism.

Referring to \eqref{bra-X-J}, the LCCos Hamiltonian vector field
\(\mathcal E_H\) and the LCCos bracket
\(\{\bullet,\bullet\}^{\mathrm{LCCos}}\) in \eqref{Pois-LCC} are related by
\begin{equation}
\mathcal E_H(F)
=
\{F,H\}^{\mathrm{LCCos}}
-
F\,\mathcal V_\varphi(H).
\end{equation}
Since
\(\mathcal V_\varphi(H)=-\langle\varphi,\mathcal E_H\rangle\), this becomes
\begin{equation}\label{LCCos-vector-bracket-relation}
    \mathcal E_H(F)
    =
    \{F,H\}^{\mathrm{LCCos}}
    +
    F\langle\varphi,\mathcal E_H\rangle .
\end{equation}
This relation is the LCCos counterpart of the bracket--vector field identity
in the LCS setting. It will be used below to rewrite the Lie--Poisson kinetic
equation in transport and continuity form.

The Hamiltonian function is not preserved in general. Indeed, since
\(\iota_{\mathcal E_H}\Omega(\mathcal E_H)=0\) and
\(\iota_{\mathcal E_H}\tau=0\), equation \eqref{LCCos-Ham-global} gives
\begin{equation}\label{LCCos-Ham-H-evolution}
    \mathcal L_{\mathcal E_H}H
    =
    H\langle\varphi,\mathcal E_H\rangle .
\end{equation}
Hence the zero-energy hypersurface \(H=0\) is invariant under the LCCos
Hamiltonian flow. The defining LCCos forms are also not preserved in general.
Using \(d\tau=\varphi\wedge\tau\), \(d\Omega=\varphi\wedge\Omega\), and
\eqref{LCCos-Ham-global}, one obtains
\begin{equation}\label{LCCos-Ham-form-evolution}
    \mathcal L_{\mathcal E_H}\tau
    =
    \langle\varphi,\mathcal E_H\rangle\tau,
    \qquad
    \mathcal L_{\mathcal E_H}\Omega
    =
    \langle\varphi,\mathcal E_H\rangle\Omega
    -
    d\langle d_\varphi H,\mathcal R\rangle\wedge\tau .
\end{equation}
Thus, the LCCos Hamiltonian flow preserves neither \(\tau\) nor \(\Omega\)
strictly; rather, it preserves their conformal class up to the additional
Reeb-direction correction encoded by
\(\langle d_\varphi H,\mathcal R\rangle\). Consequently, the natural almost
cosymplectic volume form
\(\mathrm{Vol}_{\tau,\Omega}=\tau\wedge\Omega^n\) evolves by
\begin{equation}\label{LCCos-Ham-volume-evolution}
    \mathcal L_{\mathcal E_H}\mathrm{Vol}_{\tau,\Omega}
    =
    (n+1)
    \langle\varphi,\mathcal E_H\rangle
    \mathrm{Vol}_{\tau,\Omega}.
\end{equation}
Equivalently, the divergence of \(\mathcal E_H\) with respect to
\(\mathrm{Vol}_{\tau,\Omega}\) is
\begin{equation}\label{LCCos-Ham-divergence}
    \operatorname{div}_{\tau,\Omega}(\mathcal E_H)
    =
    (n+1)
    \langle\varphi,\mathcal E_H\rangle .
\end{equation}
Therefore, the LCCos Hamiltonian dynamics is not volume-preserving with
respect to the ordinary almost cosymplectic volume; the failure of volume
preservation is measured by the Lee form. 

 The Jacobi pair \((\Lambda^\tau,-\mathcal V_{\varphi})\) in
\eqref{LCCos-Jacobi-pair} controls the Lie-algebraic behaviour of the LCCos
Hamiltonian vector fields. There are two points to emphasize. First, the
LCCos Hamiltonian vector fields are closed under the negative Jacobi--Lie
bracket. Indeed, by referring to the LCCos bracket in \eqref{Pois-LCC}, one
has
\begin{equation}\label{LCCos-Ham-Lie-closure}
    [\mathcal E_F,\mathcal E_H]_{-JL}
    =
    \mathcal E_{\{F,H\}^{\mathrm{LCCos}}} .
\end{equation}
Thus the space of LCCos Hamiltonian vector fields is a Lie algebra, denoted by
\begin{equation}\label{LCCos-Ham-vfield}
    \mathfrak X_{\mathrm{LCCos-ham}}(M)
    =
    \bigl\{
    \mathcal E_H\in\mathfrak X(M):
    \iota_{\mathcal E_H}\tau=0,\ 
    \iota_{\mathcal E_H}\Omega
    =
    d_\varphi H
    -
    \langle d_\varphi H,\mathcal R\rangle\tau
    \bigr\}.
\end{equation}
Consequently, the characteristic distribution generated by the LCCos
Hamiltonian vector fields is involutive and hence integrable.

Second, this Lie algebra is obtained as the image of a natural Lie algebra
homomorphism from the function algebra. More precisely, since
\(\mathcal F(M)\) is a Lie algebra with respect to the LCCos bracket
\eqref{Pois-LCC}, the assignment
\[
    \mathcal F(M)
    \longrightarrow
    \mathfrak X_{\mathrm{LCCos-ham}}(M),
    \qquad
    H
    \longmapsto
    \mathcal E_H
\]
is a Lie algebra homomorphism. Its kernel consists of those Hamiltonian
functions whose twisted differential has only a Reeb component:
\begin{equation}\label{LCCos-Ham-kernel}
    \ker(H\mapsto\mathcal E_H)
    =
    \left\{
    H\in\mathcal F(M):
    d_\varphi H
    =
    \langle d_\varphi H,\mathcal R\rangle\tau
    \right\}.
\end{equation}
Equivalently, \(H\) lies in the kernel precisely when \(d_\varphi H\) vanishes
on the distribution \(\ker\tau\). In particular, the conformal unit
\(\mathbb I\) belongs to this kernel, since \(d_\varphi\mathbb I=0\), and
therefore \(\mathcal E_{\mathbb I}=0\). Hence
\begin{equation}\label{LCCos-Ham-quotient}
    \mathfrak X_{\mathrm{LCCos-ham}}(M)
    \simeq
    \mathcal F(M)\big/
    \ker(H\mapsto\mathcal E_H).
\end{equation}

\noindent \textbf{Coordinate Realization of LCCos Hamiltonian Dynamics.} 
We next compute the Hamiltonian dynamics determined by this LCCos structure in
the Darboux coordinates \((t,q^i,p_i)\). For this purpose, we use the Lee form
given in \eqref{haaland} and the coordinate expressions of the LCCos forms
\(\tau_\varphi\) and \(\Omega_\varphi\) in \eqref{LCC-RTQ-forms}. For a
Hamiltonian function \(H=H(t,q,p)\), the associated Hamiltonian vector field is
\begin{equation}\label{xiH-RTQ-local}
\mathcal{E}_H
=
\frac{t\varphi_i}{1-t\phi} \frac{\partial H}{\partial p_i}\frac{\partial}{\partial t}
+
\frac{\partial H}{\partial p_i}\frac{\partial}{\partial q^i}
+
\Bigl(
-\frac{\partial H}{\partial q^i}
+\frac{C_{ij}}{1-t\phi}\frac{\partial H}{\partial p_j}
+ \frac{\varphi_i}{1-t\phi}H
-\frac{t\varphi_i}{1-t\phi}\frac{\partial H}{\partial t}
\Bigr)\frac{\partial}{\partial p_i}.
\end{equation}
Therefore, the LCCos Hamilton equations are
\begin{equation}\label{LCC-eom-RTQ}
\dot t=\frac{t\varphi_i}{1-t\phi}\frac{\partial H}{\partial p_i},
\qquad
\dot q^i=\frac{\partial H}{\partial p_i},
\qquad
\dot p_i
=
-\frac{\partial H}{\partial q^i}
+\frac{C_{ij}}{1-t\phi}\frac{\partial H}{\partial p_j}
+ \frac{\varphi_i}{1-t\phi}H
-\frac{t\varphi_i}{1-t\phi}\frac{\partial H}{\partial t}.
\end{equation}
These equations show that the conformal modification affects not only the
momentum equation but also the evolution of the time coordinate itself. In
particular, the condition \(\tau_\varphi(\mathcal E_H)=0\) does not imply
\(\dot t=0\), because \(\tau_\varphi=dt-t\varphi\). Hence, unlike the ordinary
cosymplectic Hamiltonian vector field, the LCCos Hamiltonian vector field may
have a non-trivial \(t\)-component.

If \(\sigma=\sigma(q)\), then \(\phi=0\). The Hamiltonian vector field $\mathcal{E}_H$ in \eqref{xiH-RTQ-local} takes
the simpler form
\begin{equation}\label{xiH-RTQ-local2}
\mathcal{E}_H
=
t\varphi_i \frac{\partial H}{\partial p_i}\frac{\partial}{\partial t}
+
\frac{\partial H}{\partial p_i}\frac{\partial}{\partial q^i}
+
\Bigl(
-\frac{\partial H}{\partial q^i}
+C_{ij}\frac{\partial H}{\partial p_j}
+\varphi_i H
-t\varphi_i\frac{\partial H}{\partial t}
\Bigr)\frac{\partial}{\partial p_i}.
\end{equation}
In this case, the LCCos Hamilton equations become
\begin{equation}\label{LCC-eom-RTQ2}
\dot t=t\varphi_i\frac{\partial H}{\partial p_i},
\qquad
\dot q^i=\frac{\partial H}{\partial p_i},
\qquad
\dot p_i
=
-\frac{\partial H}{\partial q^i}
+C_{ij}\frac{\partial H}{\partial p_j}
+ \varphi_i H
-t\varphi_i\frac{\partial H}{\partial t}.
\end{equation}
Thus, even in the case where the Lee form has no \(dt\)-component, the LCCos
Hamiltonian dynamics still differs from the canonical cosymplectic dynamics.
The skew term
\(C_{ij}\partial H/\partial p_j\) is the analogue of the gyroscopic correction already encountered in the LCS setting. The term
\(\varphi_i H\) comes from the twisted differential
\(d_\varphi H=dH-H\varphi\), whereas the term
\(-t\varphi_i\partial H/\partial t\) is a genuinely cosymplectic correction
coming from the modified one-form \(\tau_\varphi=dt-t\varphi\). Therefore, the
LCCos Hamiltonian equations combine the LCS-type conformal correction in the
momentum variables with an additional coupling to the time direction.

Let us recall the energy evolution identity
\(\mathcal L_{\mathcal E_H}H=H\langle\varphi,\mathcal E_H\rangle\) in
\eqref{LCCos-Ham-H-evolution}. This identity implies that the zero-energy
hypersurface \(H=0\) is invariant under the LCCos Hamiltonian flow. Hence, the
on-shell dynamics on \(H=0\) can be studied as a distinguished restriction of
the full LCCos Hamiltonian dynamics. On this invariant hypersurface, the term
\(\varphi_iH/(1-t\phi)\) vanishes. Therefore, in the general case, the
on-shell LCCos Hamiltonian equations reduce to
\begin{equation} 
\dot t
=
\frac{t\varphi_i}{1-t\phi}\frac{\partial H}{\partial p_i} , \qquad 
 \dot q^i
=
\frac{\partial H}{\partial p_i} ,  \qquad 
 \dot p_i
=
-\frac{\partial H}{\partial q^i}
+
\frac{C_{ij}}{1-t\phi}\frac{\partial H}{\partial p_j}
-
\frac{t\varphi_i}{1-t\phi}\frac{\partial H}{\partial t}.
\end{equation}
In the special case \(\phi=0\), this becomes
\begin{equation}
\dot t
=
t\varphi_i\frac{\partial H}{\partial p_i},
\qquad
\dot q^i
=
\frac{\partial H}{\partial p_i},
\qquad
\dot p_i
=
-\frac{\partial H}{\partial q^i}
+
C_{ij}\frac{\partial H}{\partial p_j}
-
t\varphi_i\frac{\partial H}{\partial t}.
\end{equation}
Thus, even on shell, the dynamics does not reduce completely to the LCS-type
momentum correction unless \(H\) is independent of \(t\); the last term records
the remaining coupling produced by the modified one-form
\(\tau_\varphi=dt-t\varphi\).

\noindent \textbf{Locally Conformal Cosymplectic Evolution Dynamics.} In accordance with the convention stated in
Appendix~\ref{sec-LdR}, we work in the exact-Lee sector throughout the
global scalar formulation of the LCCos evolution dynamics.

We now apply the same conformal gluing principle to the cosymplectic evolution
vector field. On a local cosymplectic chart
\((U_\alpha,\tau_\alpha,\Omega_\alpha)\), the Hamiltonian evolution vector
field \(E_{H_\alpha}\) is determined by
\begin{equation}\label{local-evolution-cos}
    \iota_{E_{H_\alpha}}\tau_\alpha=1_\alpha,
    \qquad
    \iota_{E_{H_\alpha}}\Omega_\alpha
    =
    dH_\alpha-\mathcal{R}_\alpha(H_\alpha)\tau_\alpha ,
\end{equation}
where \(1_\alpha\) denotes the local unit. In the locally conformal setting,
this unit is not a globally fixed constant. It is the local representative of
the normalization of the evolution field, and this normalization must transform
with the same conformal weight as \(\tau_\alpha\). Hence, on overlaps
\(U_\alpha\cap U_\beta\), we impose
\begin{equation}\label{local-unit-gluing}
    e^{\sigma_\alpha}1_\alpha
    =
    e^{\sigma_\beta}1_\beta .
\end{equation}
The local units therefore glue to a global conformal unit \(\mathbb I\)
defined by
\begin{equation}\label{global-weighted-unit}
    \mathbb I|_\alpha=e^{\sigma_\alpha}1_\alpha .
\end{equation}
Equivalently, \(\mathbb I\) is a nowhere-vanishing global scalar element of the kernel of the
Lichnerowicz--de~Rham differential already introduced in
\eqref{twisted-closed-I}; that is,
\begin{equation}\label{twisted-closed-I-}
d_\varphi \mathbb I=0,
\qquad
\text{or equivalently}
\qquad
d\mathbb I=\varphi \mathbb I .
\end{equation}

Assume that the local Hamiltonian evolution vector fields glue to a global
vector field \(E_H\), namely,
\begin{equation}\label{global-evolution-field}
    E_H|_\alpha=E_{H_\alpha}.
\end{equation}
Then, using the gluing rules \eqref{local-forms-cosymplectic} and
\eqref{global-weighted-unit}, the first equation in
\eqref{local-evolution-cos} gives
\begin{equation}
    \iota_{E_H}\tau\big|_\alpha
    =
    \iota_{E_{H_\alpha}}\bigl(e^{\sigma_\alpha}\tau_\alpha\bigr)
    =
    e^{\sigma_\alpha}\iota_{E_{H_\alpha}}\tau_\alpha
    =
    e^{\sigma_\alpha}1_\alpha
    =
    \mathbb I|_\alpha .
\end{equation}
Thus the local cosymplectic evolution equations glue to the global LCCos
evolution equations
\begin{equation}\label{global-LCCos-evolution}
    \iota_{E_H}\tau=\mathbb I,
    \qquad
    \iota_{E_H}\Omega
    =
    d_{\varphi}H
    -
    \langle d_{\varphi}H,\mathcal R\rangle\tau .
\end{equation}
Hence the local normalization
\(\iota_{E_{H_\alpha}}\tau_\alpha=1_\alpha\) globalizes not as
\(\iota_{E_H}\tau=1\), but as the weighted normalization
\(\iota_{E_H}\tau=\mathbb I\). The two-form equation, on the other hand,
globalizes with the same twisted Hamiltonian right-hand side as in the LCCos
Hamiltonian vector field case.

The definition of the LCCos evolution vector field proposed here differs from
the one commonly used\footnote{In
the usual formulation of LCCos dynamics, the second equation in
\eqref{global-LCCos-evolution} is kept unchanged, whereas the first one is
written as \(\iota_{E_H}\tau=1\).} in the literature \cite{AtEsLe23,ChLeMa91}. More precisely, in
\eqref{global-LCCos-evolution}, the right-hand side of the first equation is
the conformal unit \(\mathbb I\), which belongs to the kernel of the
Lichnerowicz--de~Rham differential. In the standard formulation, this role is
played instead by the ordinary constant function \(1\), which belongs to the
kernel of the de~Rham differential. This distinction is a direct manifestation
of the conformal weight convention adopted in \eqref{lcc-loc-rel}: the
normalization of the evolution vector field must transform compatibly with the
local one-forms \(\tau_\alpha\), and therefore globalizes as
\(\iota_{E_H}\tau=\mathbb I\), rather than as \(\iota_{E_H}\tau=1\).

Let us now record the basic conservation properties of the LCCos evolution
vector field \(E_H\). 
Evaluating the second equation in \eqref{global-LCCos-evolution} on \(E_H\), and using
\(\iota_{E_H}\Omega(E_H)=0\), we obtain
\begin{equation}\label{LCCos-evolution-H-evolution}
    \mathcal L_{E_H}H
    =
    H\langle\varphi,E_H\rangle
    +
    \mathbb I
    \langle d_\varphi H,\mathcal R\rangle .
\end{equation}
Thus, unlike the LCCos Hamiltonian vector field, the LCCos evolution vector
field does not preserve the zero-energy hypersurface in general. On \(H=0\),
one still has
\(\mathcal L_{E_H}H=\mathbb I\langle d_\varphi H,\mathcal R\rangle\), so the
restriction to \(H=0\) is invariant only under the additional condition
\(\langle d_\varphi H,\mathcal R\rangle=0\) along that hypersurface.

The defining LCCos forms are not preserved in general. Since
\(d\tau=\varphi\wedge\tau\), \(d\Omega=\varphi\wedge\Omega\), and
\(d\mathbb I=\varphi\mathbb I\), one obtains
\begin{equation}\label{LCCos-evolution-form-evolution}
    \mathcal L_{E_H}\tau
    =
    \langle\varphi,E_H\rangle\tau,
    \qquad
    \mathcal L_{E_H}\Omega
    =
    \langle\varphi,E_H\rangle\Omega
    -
    d\langle d_\varphi H,\mathcal R\rangle\wedge\tau .
\end{equation}
Consequently, the almost cosymplectic volume form
\(\mathrm{Vol}_{\tau,\Omega}=\tau\wedge\Omega^n\) evolves according to
\begin{equation}\label{LCCos-evolution-volume-evolution}
    \mathcal L_{E_H}\mathrm{Vol}_{\tau,\Omega}
    =
    (n+1)
    \langle\varphi,E_H\rangle
    \mathrm{Vol}_{\tau,\Omega}.
\end{equation}
Equivalently,
\begin{equation}\label{LCCos-evolution-divergence}
    \operatorname{div}_{\tau,\Omega}(E_H)
    =
    (n+1)
    \langle\varphi,E_H\rangle .
\end{equation}
Finally, the LCCos evolution vector fields do not form a Lie algebra in
general. Indeed, using \(\iota_{E_H}\tau=\mathbb I\) and
\(\mathcal L_{E_F}\tau=\langle\varphi,E_F\rangle\tau\), we find
\begin{equation}\label{LCCos-evolution-bracket-tau}
    \iota_{[E_F,E_H]}\tau
    =
    \mathcal L_{E_F}(\iota_{E_H}\tau)
    -
    \iota_{E_H}(\mathcal L_{E_F}\tau)
    =
    E_F(\mathbb I)
    -
    \langle\varphi,E_F\rangle\mathbb I
    =
    0 .
\end{equation}
Hence the bracket \([E_F,E_H]\) is tangent to \(\ker\tau\), whereas an
evolution vector field must satisfy
\(\iota_{E_K}\tau=\mathbb I\). Therefore, the bracket of two LCCos evolution
vector fields cannot, in general, be another LCCos evolution vector field.
This is the locally conformal counterpart of the usual fact that evolution
vector fields do not form a Lie algebra. 
 
\noindent \textbf{Comparison of LCCos Hamiltonian and Evolution Dynamics.}
We now record the relation between the LCCos Hamiltonian vector field and the
LCCos evolution vector field:
\begin{equation}\label{LCCos-evolution-Ham-Reeb}
    E_H
    =
    \mathcal E_H+\mathbb I\,\mathcal R .
\end{equation}
Indeed, a direct computation proves that
\begin{equation}
\begin{split}
    \iota_{\mathcal E_H+\mathbb I\mathcal R}\tau
    & =
    \iota_{\mathcal E_H}\tau
    +
    \mathbb I\,\iota_{\mathcal R}\tau
    =
    \mathbb I, \\
    \iota_{\mathcal E_H+\mathbb I\mathcal R}\Omega
    & =
    \iota_{\mathcal E_H}\Omega
    +
    \mathbb I\,\iota_{\mathcal R}\Omega
    =
    d_\varphi H
    -
    \langle d_\varphi H,\mathcal R\rangle\tau .
\end{split}
\end{equation}
Thus, the LCCos evolution vector field is obtained from the LCCos Hamiltonian
vector field by adding the conformal Reeb contribution \(\mathbb I\mathcal R\).

It is important, however, to distinguish this global construction from the
normalized comparison obtained by adding only the Reeb vector field. If one
sets
\[
    \widetilde E_H
    =
    \mathcal E_H+\mathcal R ,
\]
then the resulting vector field satisfies
\[
    \iota_{\widetilde E_H}\tau=1,
    \qquad
    \iota_{\widetilde E_H}\Omega
    =
    d_\varphi H
    -
    \langle d_\varphi H,\mathcal R\rangle\tau .
\]
Both $E_H$ and $\widetilde E_H$ are globally defined.
Our choice $E_H=\mathcal E_H+\mathbb I\mathcal R$ is distinguished
by its compatibility with the conformal gluing of the local
normalizations in \eqref{local-unit-gluing}.

\noindent \textbf{Coordinate Realization of LCCos Evolution Dynamics.}
We next compute the evolution dynamics determined by the same LCCos structure
in the Darboux coordinates \((t,q^i,p_i)\). As before, we use the Lee form
given in \eqref{haaland} and the coordinate expressions of the LCCos forms
\(\tau_\varphi\) and \(\Omega_\varphi\) in \eqref{LCC-RTQ-forms}. In the local
chart under consideration, we write
\[
    \varphi|_\alpha=d\sigma_\alpha,
    \qquad
    \mathbb I|_\alpha=e^{\sigma_\alpha(t,q)} .
\]
Thus the conformal unit is kept in the global restricted form. For a
Hamiltonian function \(H=H(t,q,p)\), the associated LCCos evolution vector
field is
\begin{equation}\label{EH-RTQ-local}
\begin{split}
    E_H \vert_\alpha
& =
\frac{1}{1-t\phi}
\left(
e^{\sigma_\alpha(t,q)}
+
t\varphi_i\frac{\partial H}{\partial p_i}
\right)
\frac{\partial}{\partial t}
+
\frac{\partial H}{\partial p_i}
\frac{\partial}{\partial q^i}
\\
&\qquad
+
\left(
-\frac{\partial H}{\partial q^i}
+
\frac{C_{ij}}{1-t\phi}
\frac{\partial H}{\partial p_j}
+
\frac{\varphi_i}{1-t\phi}H
-
\frac{t\varphi_i}{1-t\phi}
\frac{\partial H}{\partial t}
+
\frac{p_i\phi}{1-t\phi}
e^{\sigma_\alpha(t,q)}
\right)
\frac{\partial}{\partial p_i}.
\end{split}
\end{equation}
Therefore, the corresponding LCCos evolution equations are
\begin{equation}\label{LCCos-evo-eq}
\begin{split}
    \dot t
    &=
    \frac{1}{1-t\phi}
    \left(
    e^{\sigma_\alpha(t,q)}
    +
    t\varphi_i\frac{\partial H}{\partial p_i}
    \right),
\\
    \dot q^i
    &=
    \frac{\partial H}{\partial p_i},
\\
    \dot p_i
    &=
    -\frac{\partial H}{\partial q^i}
    +
    \frac{C_{ij}}{1-t\phi}
    \frac{\partial H}{\partial p_j}
    +
    \frac{\varphi_i}{1-t\phi}H
    -
    \frac{t\varphi_i}{1-t\phi}
    \frac{\partial H}{\partial t}
    +
    \frac{p_i\phi}{1-t\phi}
    e^{\sigma_\alpha(t,q)} .
\end{split}
\end{equation}
These equations should be compared with the LCCos Hamiltonian equations in
\eqref{LCC-eom-RTQ}. The most visible difference is the normalization of the
time component: the evolution vector field satisfies
\(\iota_{E_H}\tau_\varphi=\mathbb I\), rather than
\(\iota_{\mathcal E_H}\tau_\varphi=0\). Consequently, the equation for
\(\dot t\) contains the conformal unit contribution
\(e^{\sigma_\alpha(t,q)}\). Moreover, the momentum equation acquires the extra
conformal Reeb term
\(p_i\phi e^{\sigma_\alpha(t,q)}/(1-t\phi)\), which comes from the
\(dt\wedge dq^i\) part of \(\Omega_\varphi\) together with the nonzero
\(\partial/\partial t\)-component of the evolution vector field.

If instead one considers only the normalized comparison field
\(\widetilde E_H=\mathcal E_H+\mathcal R\), then the same coordinate
calculation would be obtained from \eqref{EH-RTQ-local} by replacing the
conformal unit contribution \(e^{\sigma_\alpha(t,q)}\) by \(1\) in the Reeb
terms. This field is also globally defined, but it uses the ordinary
unit normalization rather than the conformally weighted one. The evolution field adopted here is the one in
\eqref{EH-RTQ-local}, where the Reeb contribution is multiplied by the
restricted conformal unit \(\mathbb I|_\alpha=e^{\sigma_\alpha(t,q)}\).

If \(\sigma=\sigma(q)\), then \(\phi=0\), and the evolution vector field
\eqref{EH-RTQ-local} takes the simpler form
\begin{equation}\label{EH-RTQ-local2}
\begin{split}
    E_H \vert_\alpha
& =
\left(
e^{\sigma_\alpha(q)}
+
t\varphi_i\frac{\partial H}{\partial p_i}
\right)
\frac{\partial}{\partial t}
+
\frac{\partial H}{\partial p_i}
\frac{\partial}{\partial q^i}
\\
&\qquad
+
\left(
-\frac{\partial H}{\partial q^i}
+
C_{ij}\frac{\partial H}{\partial p_j}
+
\varphi_i H
-
t\varphi_i
\frac{\partial H}{\partial t}
\right)
\frac{\partial}{\partial p_i}.
\end{split}
\end{equation}
In this case, the LCCos evolution equations become
\begin{equation}\label{LCCos-evo-eq2}
\begin{split}
    \dot t
    &=
    e^{\sigma_\alpha(q)}
    +
    t\varphi_i\frac{\partial H}{\partial p_i},
\\
    \dot q^i
    &=
    \frac{\partial H}{\partial p_i},
\\
    \dot p_i
    &=
    -\frac{\partial H}{\partial q^i}
    +
    C_{ij}\frac{\partial H}{\partial p_j}
    +
    \varphi_i H
    -
    t\varphi_i
    \frac{\partial H}{\partial t}.
\end{split}
\end{equation}
Thus, in contrast with the Hamiltonian vector field, the evolution vector
field keeps the time-direction contribution already present in cosymplectic
dynamics, but this contribution is conformally weighted by
\(e^{\sigma_\alpha}\) in the global restricted representation. The conformal
terms then modify both the time equation and the momentum equation. The term
\(C_{ij}\partial H/\partial p_j\) is the LCS-type skew correction, the term
\(\varphi_iH\) comes from the twisted differential \(d_\varphi H\), and the
term \(-t\varphi_i\partial H/\partial t\) records the coupling caused by the
modified one-form \(\tau_\varphi=dt-t\varphi\). For the normalized comparison field \(\widetilde E_H=\mathcal E_H+\mathcal R\),
the same expressions are obtained by replacing the conformal time contribution
\(e^{\sigma_\alpha}\) by \(1\) in the Reeb part.

Finally, unlike the LCCos Hamiltonian vector field, the LCCos evolution vector
field does not in general preserve the zero-energy hypersurface. Indeed, as
recorded in \eqref{LCCos-evolution-H-evolution}, the energy evolves according
to
\(\mathcal L_{E_H}H
=
H\langle\varphi,E_H\rangle
+
\mathbb I\langle d_\varphi H,\mathcal R\rangle\).
Therefore, even on \(H=0\), one generally has
\(\mathcal L_{E_H}H
=
\mathbb I\langle d_\varphi H,\mathcal R\rangle\), which need not vanish.
Hence there is no natural on-shell restriction for the LCCos evolution
dynamics unless the additional condition
\(\langle d_\varphi H,\mathcal R\rangle=0\) is imposed along \(H=0\).

\subsection{Continuum Dynamics on the LCCos Setting}\label{Sec-LCCos-Man}
 
We begin with a \((2n+1)\)-dimensional LCCos manifold
\((M,\tau,\Omega\,;\varphi)\). We recall the associated Jacobi structure
\((\Lambda^\tau,-\mathcal V_{\varphi})\) given in
\eqref{LCCos-Jacobi-pair}. The existence of this Jacobi structure endows the
space \(\mathcal F(M)\) of smooth functions with a Lie algebra structure via
the LCCos bracket \eqref{Pois-LCC}. Hence, the adjoint action of
\(\mathcal F(M)\) on itself is written as
\begin{equation}\label{Lie-brr-F-LCCos}
\operatorname{ad}_H F
=
\{H,F\}^{\mathrm{LCCos}}
=
\Lambda^\tau(dH,dF)
+
H\,\mathcal V_{\varphi}(F)
-
F\,\mathcal V_{\varphi}(H).
\end{equation}
The dual space of \(\mathcal F(M)\) is identified with the space of densities
\(\operatorname{Den}(M)\) through the \(L_2\)-pairing
\begin{equation}\label{LCCos-L-2}
\langle\bullet,\bullet\rangle_{L_2}:
\operatorname{Den}(M)\times\mathcal F(M)
\longrightarrow
\mathbb R,
\qquad
(f\,\operatorname{Vol}_{\tau,\Omega},H)
\longmapsto
\int_M Hf\,\operatorname{Vol}_{\tau,\Omega}.
\end{equation}
Here
\(\operatorname{Vol}_{\tau,\Omega}=\tau\wedge\Omega^n\) is the natural almost
cosymplectic volume form in \eqref{non-deg-cos}. Accordingly, a density
\(f\,\operatorname{Vol}_{\tau,\Omega}\) is represented by its density
coefficient \(f\).

We now derive the integration identity needed for the coadjoint action. Recall
the relation \eqref{LCCos-vector-bracket-relation} between the LCCos bracket
\eqref{Pois-LCC} and the LCCos Hamiltonian vector field
\eqref{LCCos-Ham-global}. We shall use it in the equivalent form
\[
\{H,F\}^{\mathrm{LCCos}}
=
-\mathcal E_H(F)
+
F\langle\varphi,\mathcal E_H\rangle .
\]
Therefore, assuming that the boundary terms vanish, we compute
\begin{equation}
\begin{split}
\int_M
f\{H,F\}^{\mathrm{LCCos}}
\operatorname{Vol}_{\tau,\Omega}
&=
-\int_M
f\,\mathcal E_H(F)
\operatorname{Vol}_{\tau,\Omega}
+
\int_M
fF\langle\varphi,\mathcal E_H\rangle
\operatorname{Vol}_{\tau,\Omega}
\\
&=
\int_M
F
\left(
\mathcal E_H(f)
+
f\,\operatorname{div}_{\tau,\Omega}(\mathcal E_H)
+
f\langle\varphi,\mathcal E_H\rangle
\right)
\operatorname{Vol}_{\tau,\Omega}.
\end{split}
\end{equation}
On the other hand, applying again the bracket--vector-field relation gives
\[
\mathcal E_H(f)
=
\{f,H\}^{\mathrm{LCCos}}
+
f\langle\varphi,\mathcal E_H\rangle .
\]
Moreover, by \eqref{LCCos-Ham-divergence},
\[
\operatorname{div}_{\tau,\Omega}(\mathcal E_H)
=
(n+1)\langle\varphi,\mathcal E_H\rangle .
\]
Substituting these two identities into the previous expression, we obtain
\[
\mathcal E_H(f)
+
f\,\operatorname{div}_{\tau,\Omega}(\mathcal E_H)
+
f\langle\varphi,\mathcal E_H\rangle
=
\{f,H\}^{\mathrm{LCCos}}
+
(n+3)\langle\varphi,\mathcal E_H\rangle f .
\]
Thus we arrive at the following integration identity for the LCCos bracket:
\begin{equation}\label{LCCos-assoc+}
\int_M
f\{H,F\}^{\mathrm{LCCos}}
\operatorname{Vol}_{\tau,\Omega}
=
\int_M
\left(
\{f,H\}^{\mathrm{LCCos}}
+
(n+3)\langle\varphi,\mathcal E_H\rangle f
\right)
F\,\operatorname{Vol}_{\tau,\Omega}.
\end{equation}

\noindent \textbf{Lie-Poisson Dynamics.}
Now we are ready to compute the coadjoint action of the Lie algebra
\(\mathcal F(M)\) on its dual space. With respect to the \(L_2\)-pairing, we
have
\begin{equation}
\begin{split}
\big\langle \operatorname{ad}^*_H f,F\big\rangle_{L^2}
&=
\int_M
(\operatorname{ad}^*_H f)F
\operatorname{Vol}_{\tau,\Omega}
\\
&=
\int_M
f\,\operatorname{ad}_H F
\operatorname{Vol}_{\tau,\Omega}
\\
&=
\int_M
f\{H,F\}^{\mathrm{LCCos}}
\operatorname{Vol}_{\tau,\Omega}
\\
&=
\int_M
\left(
\{f,H\}^{\mathrm{LCCos}}
+
(n+3)\langle\varphi,\mathcal E_H\rangle f
\right)
F\,\operatorname{Vol}_{\tau,\Omega},
\end{split}
\end{equation}
where we have used \eqref{LCCos-assoc+} in the last line. Since this holds for
all \(F\in\mathcal F(M)\), the coadjoint action on \(\mathcal F^*(M)\) is
given by
\begin{equation}\label{LCCos-coadjoint}
\operatorname{ad}^*_H f
=
\{f,H\}^{\mathrm{LCCos}}
+
(n+3)\langle\varphi,\mathcal E_H\rangle f .
\end{equation}

The dual space \(\mathcal F^*(M)\) carries the corresponding Lie--Poisson
bracket. For two functionals
\(\mathfrak F,\mathfrak H:\mathcal F^*(M)\longrightarrow\mathbb R\), we define
\begin{equation}\label{LCCos-Lie-Poisson-bracket}
\left\{
\mathfrak F,
\mathfrak H
\right\}_{\mathcal F^*(M)}(f)
=
\int_M
f
\left\{
\frac{\delta\mathfrak F}{\delta f},
\frac{\delta\mathfrak H}{\delta f}
\right\}^{\mathrm{LCCos}}
\operatorname{Vol}_{\tau,\Omega}.
\end{equation}
Given a Hamiltonian functional \(\mathfrak H\), the corresponding
Lie--Poisson equation is
\begin{equation}\label{LP-f-LCCos}
\dot f
=
-
\operatorname{ad}^*_{\delta\mathfrak H/\delta f}f .
\end{equation}
We choose the linear Hamiltonian functional
\begin{equation}\label{Ham-functional-LCCos}
\mathfrak H(f)
=
\int_M
Hf\,\operatorname{Vol}_{\tau,\Omega},
\end{equation}
for which \(\delta\mathfrak H/\delta f=H\). Substituting
\eqref{LCCos-coadjoint} into \eqref{LP-f-LCCos}, we obtain the LCCos Vlasov
equation
\begin{equation}\label{LCCos-Vlasov}
\dot f
+
\{f,H\}^{\mathrm{LCCos}}
=
-
(n+3)\langle\varphi,\mathcal E_H\rangle f .
\end{equation}

\noindent \textbf{Continuity Dynamics.}
Using the relation between the LCCos bracket and the LCCos Hamiltonian vector
field, equation \eqref{LCCos-Vlasov} can be written in transport form as
\begin{equation}\label{LCCos-Vlasov-transport-form}
\dot f
+
\mathcal E_H(f)
=
-
(n+2)
\langle\varphi,\mathcal E_H\rangle f .
\end{equation}
Recalling
\[
\operatorname{div}_{\tau,\Omega}(\mathcal E_H)
=
(n+1)\langle\varphi,\mathcal E_H\rangle,
\]
where the divergence is computed with respect to
\(\operatorname{Vol}_{\tau,\Omega}=\tau\wedge\Omega^n\), we obtain
\begin{equation}\label{LCCos-Vlasov-continuity-source-scalar}
\dot f
+
\mathcal E_H(f)
+
\operatorname{div}_{\tau,\Omega}(\mathcal E_H)f
=
-
\langle\varphi,\mathcal E_H\rangle f .
\end{equation}
Equivalently,
\begin{equation}\label{LCCos-Vlasov-continuity-source}
\left(
\partial_\lambda
+
\mathcal L_{\mathcal E_H}
\right)
\left(
f\,\operatorname{Vol}_{\tau,\Omega}
\right)
=
-
\langle\varphi,\mathcal E_H\rangle
f\,\operatorname{Vol}_{\tau,\Omega}.
\end{equation}
Thus, the ordinary density \(f\,\operatorname{Vol}_{\tau,\Omega}\) is not
conserved by the LCCos Hamiltonian flow. The defect from conservation is again
measured by the Lee contribution
\(\langle\varphi,\mathcal E_H\rangle\).

As in the LCS case, the following global weighted formulation
assumes an exact Lee form and a fixed positive conformal unit.
First, we modify the volume form by using the conformally weighted unit
\(\mathbb I\), characterized by
\begin{equation}\label{LCCos-I-closed}
d_\varphi\mathbb I=0,
\qquad
d\mathbb I=\varphi\,\mathbb I .
\end{equation}
Then
\begin{equation}\label{LCCos-I-derivative}
\mathcal E_H(\mathbb I)
=
\langle\varphi,\mathcal E_H\rangle\mathbb I .
\end{equation}
Accordingly, we define the conformally weighted LCCos volume form by
\begin{equation}\label{LCCos-weighted-volume}
\operatorname{Vol}_{\tau,\Omega,\varphi}
=
\mathbb I\,\operatorname{Vol}_{\tau,\Omega}
=
\mathbb I\,\tau\wedge\Omega^n .
\end{equation}
Using \eqref{LCCos-I-derivative} and the divergence formula, one obtains
\begin{equation}
\mathcal L_{\mathcal E_H}
\operatorname{Vol}_{\tau,\Omega,\varphi}
=
(n+2)\langle\varphi,\mathcal E_H\rangle
\operatorname{Vol}_{\tau,\Omega,\varphi}.
\end{equation}
Hence \eqref{LCCos-Vlasov-transport-form} is equivalently written as the
source-free continuity equation
\begin{equation}\label{LCCos-Vlasov-continuity-weighted}
\left(
\partial_\lambda
+
\mathcal L_{\mathcal E_H}
\right)
\left(
f\,\operatorname{Vol}_{\tau,\Omega,\varphi}
\right)
=
0.
\end{equation}
 
The second equivalent way is to keep the ordinary almost cosymplectic volume
form \(\operatorname{Vol}_{\tau,\Omega}\), but to absorb the conformal weight
into the density function. We define
\begin{equation}\label{LCCos-physical-density-def}
\mathfrak f
=
\mathbb I f .
\end{equation}
Then, using \eqref{LCCos-Vlasov-transport-form},
\eqref{LCCos-I-derivative}, and the divergence formula, one obtains
\begin{equation}
\begin{split}
\frac{\partial \mathfrak f}{\partial \lambda}
+
\mathcal E_H(\mathfrak f)
+
\operatorname{div}_{\tau,\Omega}(\mathcal E_H)\mathfrak f
&=
\mathbb I
\left(
\frac{\partial f}{\partial \lambda}
+
\mathcal E_H(f)
\right)
+
f\,\mathcal E_H(\mathbb I)
+
\operatorname{div}_{\tau,\Omega}(\mathcal E_H)\mathbb I f
\\
&=
-(n+2)\langle\varphi,\mathcal E_H\rangle\mathbb I f
+
\langle\varphi,\mathcal E_H\rangle\mathbb I f
+
(n+1)\langle\varphi,\mathcal E_H\rangle\mathbb I f
=
0 .
\end{split}
\end{equation}
Consequently, \(\mathfrak f\) satisfies the ordinary source-free continuity
equation
\begin{equation}\label{LCCos-Vlasov-continuity-physical}
\left(
\partial_\lambda
+
\mathcal L_{\mathcal E_H}
\right)
\left(
\mathfrak f\,\operatorname{Vol}_{\tau,\Omega}
\right)
=
0.
\end{equation}
Equivalently, using \eqref{LCCos-vector-bracket-relation}, this equation can
be written explicitly as
\begin{equation}\label{LCCos-Vlasov-physical-bracket}
\dot{\mathfrak f}
+
\{\mathfrak f,H\}^{\mathrm{LCCos}}
=
-
(n+2)\langle\varphi,\mathcal E_H\rangle \mathfrak f .
\end{equation}

Thus the two formulations are equivalent:
\begin{equation}
f\,\operatorname{Vol}_{\tau,\Omega,\varphi}
=
f\,\mathbb I\,\operatorname{Vol}_{\tau,\Omega}
=
\mathfrak f\,\operatorname{Vol}_{\tau,\Omega}.
\end{equation}
The function \(f\) is the density coefficient adapted to the conformal LCCos
volume, whereas \(\mathfrak f=\mathbb I f\) is the density coefficient with
respect to the ordinary volume \(\tau\wedge\Omega^n\). In physical terms, the
Lee form produces a conformal rescaling of the phase-space volume. The source
term in \eqref{LCCos-Vlasov-continuity-source} expresses the mismatch between
the ordinary volume \(\tau\wedge\Omega^n\) and the natural conformally weighted
volume \(\mathbb I\,\tau\wedge\Omega^n\).

For a material domain \(D_\lambda\) transported by \(\mathcal E_H\), the
conserved particle number can be written as
\begin{equation}
\frac{d}{d\lambda}
\int_{D_\lambda}
f\,\operatorname{Vol}_{\tau,\Omega,\varphi}
=
0,
\end{equation}
or, equivalently,
\begin{equation}
\frac{d}{d\lambda}
\int_{D_\lambda}
\mathfrak f\,\operatorname{Vol}_{\tau,\Omega}
=
0.
\end{equation}

\noindent \textbf{Coordinate Formulations.}
Let \((t,q^i,p_i)\) be local coordinates on the extended cotangent bundle
realization of the LCCos manifold. We write the Lee form as
\(\varphi=\phi\,dt+\varphi_i\,dq^i\) and use the local LCCos forms given in
\eqref{LCC-RTQ-forms}. Then the LCCos Vlasov equation
\eqref{LCCos-Vlasov} admits the following coordinate realization:
\begin{equation}\label{local-LCCos-Vlasov-expanded}
\begin{split}
\frac{\partial f}{\partial \lambda}
+
\frac{\partial H}{\partial p_i}
\frac{\partial f}{\partial q^i}
-
\frac{\partial H}{\partial q^i}
\frac{\partial f}{\partial p_i}
={}&
-
\frac{t\varphi_i}{1-t\phi}
\frac{\partial H}{\partial p_i}
\frac{\partial f}{\partial t}
\\
&-
\left(
\frac{C_{ij}}{1-t\phi}
\frac{\partial H}{\partial p_j}
+
\frac{\varphi_i}{1-t\phi}H
-
\frac{t\varphi_i}{1-t\phi}
\frac{\partial H}{\partial t}
\right)
\frac{\partial f}{\partial p_i}
\\
&-
\frac{n+2}{1-t\phi}
\varphi_i
\frac{\partial H}{\partial p_i}
f .
\end{split}
\end{equation}
The zero-energy hypersurface \(H=0\) is invariant under the LCCos Hamiltonian
flow, since
\(\mathcal L_{\mathcal E_H}H
=
H\langle\varphi,\mathcal E_H\rangle\).
Restricting \eqref{local-LCCos-Vlasov-expanded} to this invariant hypersurface
removes the explicit \(H\)-term and gives
\begin{equation}\label{local-LCCos-Vlasov-on-shell}
\begin{split}
\frac{\partial f}{\partial \lambda}
+
\frac{\partial H}{\partial p_i}
\frac{\partial f}{\partial q^i}
-
\frac{\partial H}{\partial q^i}
\frac{\partial f}{\partial p_i}
={}&
-
\frac{t\varphi_i}{1-t\phi}
\frac{\partial H}{\partial p_i}
\frac{\partial f}{\partial t}
\\
&-
\left(
\frac{C_{ij}}{1-t\phi}
\frac{\partial H}{\partial p_j}
-
\frac{t\varphi_i}{1-t\phi}
\frac{\partial H}{\partial t}
\right)
\frac{\partial f}{\partial p_i}
\\
&-
\frac{n+2}{1-t\phi}
\varphi_i
\frac{\partial H}{\partial p_i}
f .
\end{split}
\end{equation}
In the special case \(\phi=0\), equation
\eqref{local-LCCos-Vlasov-expanded} becomes
\begin{equation}\label{local-LCCos-Vlasov-phi-zero}
\begin{split}
\frac{\partial f}{\partial \lambda}
+
\frac{\partial H}{\partial p_i}
\frac{\partial f}{\partial q^i}
-
\frac{\partial H}{\partial q^i}
\frac{\partial f}{\partial p_i}
={}&
-
t\varphi_i
\frac{\partial H}{\partial p_i}
\frac{\partial f}{\partial t}
\\
&-
\left(
C_{ij}
\frac{\partial H}{\partial p_j}
+
\varphi_i H
-
t\varphi_i
\frac{\partial H}{\partial t}
\right)
\frac{\partial f}{\partial p_i}
\\
&-
(n+2)
\varphi_i
\frac{\partial H}{\partial p_i}
f .
\end{split}
\end{equation}

For the density coefficient \(\mathfrak f=\mathbb I f\) with respect to the
ordinary almost cosymplectic volume form \(\operatorname{Vol}_{\tau,\Omega}\),
the corresponding coordinate equation is
\begin{equation}\label{local-LCCos-Vlasov-physical-expanded}
\begin{split}
\frac{\partial \mathfrak f}{\partial \lambda}
+
\frac{\partial H}{\partial p_i}
\frac{\partial \mathfrak f}{\partial q^i}
-
\frac{\partial H}{\partial q^i}
\frac{\partial \mathfrak f}{\partial p_i}
={}&
-
\frac{t\varphi_i}{1-t\phi}
\frac{\partial H}{\partial p_i}
\frac{\partial \mathfrak f}{\partial t}
\\
&-
\left(
\frac{C_{ij}}{1-t\phi}
\frac{\partial H}{\partial p_j}
+
\frac{\varphi_i}{1-t\phi}H
-
\frac{t\varphi_i}{1-t\phi}
\frac{\partial H}{\partial t}
\right)
\frac{\partial \mathfrak f}{\partial p_i}
\\
&-
\frac{n+1}{1-t\phi}
\varphi_i
\frac{\partial H}{\partial p_i}
\mathfrak f .
\end{split}
\end{equation}
Here the only difference between the \(f\)-equation and the
\(\mathfrak f\)-equation is the last conformal source coefficient: it is
\(n+2\) for \(f\), whereas it is \(n+1\) for \(\mathfrak f\). If one restricts
to \(H=0\), the explicit \(H\)-term in the momentum correction disappears also
in \eqref{local-LCCos-Vlasov-physical-expanded}. If \(\phi=0\), all factors
\((1-t\phi)^{-1}\) disappear and the remaining conformal corrections are
controlled only by the spatial Lee components \(\varphi_i\). Finally, if the
Lee form vanishes identically, then \(C_{ij}=0\), \(\mathbb I=1\),
\(\mathfrak f=f\), and both equations reduce to the ordinary cosymplectic
Vlasov equation. 
 
\subsection{Time-Dependent Mass Dynamics in LCCos Geometry}\label{sec-time-mass}

Consider the time-dependent Hamiltonian function
\begin{equation}\label{TD-geodesic-H}
H(t,q,p)
=
\frac{1}{2}g^{ij}(t,q)p_i p_j,
\end{equation}
where \(g^{ij}(t,q)\) is the inverse of a time-dependent Riemannian metric
\(g_{ij}(t,q)\) on \(Q\). Equivalently, \(g^{ij}(t,q)\) may be interpreted as
a time-dependent inverse mass tensor. We compare the ordinary cosymplectic
Hamiltonian dynamics with its locally conformal cosymplectic deformation, and
then write the corresponding Vlasov equations.

\noindent \textbf{Cosymplectic Dynamics.} In the ordinary cosymplectic case on \(\mathbb R\times T^*Q\), with
\(\tau=dt\) and \(\Omega_Q=dq^i\wedge dp_i\), the Hamiltonian vector field is
tangent to \(\ker dt\). Hence \(\dot t=0\), and the Hamiltonian equations are
\begin{equation}\label{TD-geodesic-cos-Ham}
    \dot t=0,
    \qquad
    \dot q^i
    =
    g^{ij}(t,q)p_j,
    \qquad
    \dot p_i
    =
    -\frac{1}{2}
    \frac{\partial g^{jk}}{\partial q^i}(t,q)p_jp_k .
\end{equation}
Equivalently, since \(p_i=g_{ij}(t,q)\dot q^j\) and \(t\) is constant along
the Hamiltonian flow, the projected dynamics on \(Q\) is the geodesic equation
for the metric \(g(t,\cdot)\) at fixed \(t\):
\begin{equation}\label{TD-geodesic-cos-Newton}
    \ddot q^\ell
    +
    \Gamma^\ell_{ij}(t,q)\dot q^i\dot q^j
    =
    0 .
\end{equation}
Here \(\Gamma^\ell_{ij}(t,q)\) denotes the Christoffel symbols of the metric
\(g(t,\cdot)\), namely
\begin{equation}\label{TD-Christoffel}
\Gamma^\ell_{ij}
=
\frac{1}{2}g^{\ell k}
\left(
\frac{\partial g_{jk}}{\partial q^i}
+
\frac{\partial g_{ik}}{\partial q^j}
-
\frac{\partial g_{ij}}{\partial q^k}
\right).
\end{equation}

\noindent \textbf{LCCos Hamiltonian Dynamics.} 
We now pass to the LCCos Hamiltonian dynamics. Let the Lee form be written
locally as
\(\varphi=\phi(t,q)\,dt+\varphi_i(t,q)\,dq^i\), and let
\(C_{ij}=\varphi_jp_i-\varphi_i p_j\). Substituting
\eqref{TD-geodesic-H} into the local LCCos Hamilton equations gives
\begin{equation}\label{TD-geodesic-LCCos-Ham}
\begin{split}
    \dot t
    &=
    \frac{t\varphi_i}{1-t\phi}
    g^{ij}p_j,
\\
    \dot q^i
    &=
    g^{ij}p_j,
\\
    \dot p_i
    &=
    -
    \frac{1}{2}
    \frac{\partial g^{jk}}{\partial q^i}p_jp_k
    +
    \frac{C_{ij}}{1-t\phi}g^{jk}p_k
    +
    \frac{\varphi_i}{2(1-t\phi)}g^{jk}p_jp_k
    -
    \frac{t\varphi_i}{2(1-t\phi)}
    \frac{\partial g^{jk}}{\partial t}p_jp_k .
\end{split}
\end{equation}
Thus, unlike in the ordinary cosymplectic Hamiltonian case, the variable \(t\)
is no longer frozen in general. The LCCos deformation couples the geodesic
motion on \(Q\) to the distinguished time direction. Let \(v^i=\dot q^i\), so
that \(p_i=g_{ij}v^j\). Then the LCCos Hamilton equations can be rewritten in
Newton form as
\begin{equation}\label{TD-geodesic-LCCos-Newton}
    \ddot q^\ell
    +
    \Gamma^\ell_{ij}(t,q)\dot q^i\dot q^j
    =
    -
    g^{\ell i}
    \frac{\partial g_{ij}}{\partial t}
    \dot t\,\dot q^j
    +
    \frac{1}{1-t\phi}
    \left(
    \varphi_i\dot q^i\,\dot q^\ell
    -
    \frac{1}{2}g^{\ell i}\varphi_i
    g_{ab}\dot q^a\dot q^b
    +
    \frac{t}{2}g^{\ell i}\varphi_i
    \frac{\partial g_{ab}}{\partial t}\dot q^a\dot q^b
    \right),
\end{equation}
where
\begin{equation}\label{TD-geodesic-LCCos-tdot}
    \dot t
    =
    \frac{t}{1-t\phi}
    \varphi_i\dot q^i .
\end{equation}
If \(\phi=0\), that is, if the Lee form has no \(dt\)-component, then
\eqref{TD-geodesic-LCCos-Newton} reduces to
\begin{equation}\label{TD-geodesic-LCCos-Newton-phi-zero}
\begin{split}
    \ddot q^\ell
    +
    \Gamma^\ell_{ij}(t,q)\dot q^i\dot q^j
    &=
    -
    g^{\ell i}
    \frac{\partial g_{ij}}{\partial t}
    \dot t\,\dot q^j
    +
    \varphi_i\dot q^i\,\dot q^\ell
\\
    &\quad
    -
    \frac{1}{2}g^{\ell i}\varphi_i
    g_{ab}\dot q^a\dot q^b
    +
    \frac{t}{2}g^{\ell i}\varphi_i
    \frac{\partial g_{ab}}{\partial t}\dot q^a\dot q^b,
    \qquad
    \dot t=t\varphi_i\dot q^i .
\end{split}
\end{equation}

\noindent \textbf{Continuum Dynamics for a 3D Model.} 
As a concrete specialization, let \(Q\) be three-dimensional and consider the
time-dependent mass metric
\begin{equation}\label{TD-mass-metric-3D}
    g_{ij}(t)=m(t)\delta_{ij},
    \qquad
    g^{ij}(t)=\frac{1}{m(t)}\delta^{ij}.
\end{equation}
Then
\begin{equation}\label{TD-mass-H-3D}
    H(t,\boldsymbol q,\boldsymbol p)
    =
    \frac{1}{2m(t)}
    |\boldsymbol p|^2 .
\end{equation}
Assume moreover that \(\phi=0\) and that the spatial part of the Lee form is
constant,
\[
    \varphi=\varphi_i\,dq^i,
    \qquad
    \boldsymbol\varphi=(\varphi_1,\varphi_2,\varphi_3).
\]
Since the metric has no \(q\)-dependence, the Christoffel symbols vanish. The
LCCos Hamilton equations become
\begin{equation}\label{TD-mass-LCCos-Ham-3D}
    \dot t
    =
    t\frac{\langle\boldsymbol\varphi,\boldsymbol p\rangle}{m(t)},
    \qquad
    \dot{\boldsymbol q}
    =
    \frac{\boldsymbol p}{m(t)},
    \qquad
    \dot{\boldsymbol p}
    =
    \frac{\langle\boldsymbol\varphi,\boldsymbol p\rangle}{m(t)}
    \boldsymbol p
    +
    \frac{|\boldsymbol p|^2}{2m(t)}
    \left(
    t\frac{\dot m(t)}{m(t)}-1
    \right)
    \boldsymbol\varphi .
\end{equation}
Equivalently, the projected Newton equation is
\begin{equation}\label{TD-mass-LCCos-Newton-3D-reduced}
    \ddot{\boldsymbol q}
    =
    \left(
    1
    -
    t\frac{\dot m(t)}{m(t)}
    \right)
    \left(
    \langle\boldsymbol\varphi,\dot{\boldsymbol q}\rangle
    \dot{\boldsymbol q}
    -
    \frac{1}{2}
    |\dot{\boldsymbol q}|^2
    \boldsymbol\varphi
    \right).
\end{equation}
Thus, in three dimensions, the LCCos Hamiltonian deformation produces a
directional force determined by the Lee vector \(\boldsymbol\varphi\). The
factor \(1-t\dot m(t)/m(t)\) records the additional effect of the
time-dependent mass. If \(m(t)\) is constant, this factor becomes \(1\).
For the LCCos Hamiltonian vector field \eqref{TD-mass-LCCos-Ham-3D}, since
\(Q\) is three-dimensional, one has \(n=3\). Therefore the LCCos Vlasov
equation for \(f\) becomes
\begin{equation}\label{3D-LCCos-mass-Vlasov}
\begin{split}
\frac{\partial f}{\partial \lambda}
+
\frac{\boldsymbol p}{m(t)}
\cdot
\nabla_{\boldsymbol q}f
={}&
-
t\frac{\langle\boldsymbol\varphi,\boldsymbol p\rangle}{m(t)}
\frac{\partial f}{\partial t}
\\
&-
\left[
\frac{\langle\boldsymbol\varphi,\boldsymbol p\rangle}{m(t)}
\boldsymbol p
+
\frac{|\boldsymbol p|^2}{2m(t)}
\left(
t\frac{\dot m(t)}{m(t)}-1
\right)
\boldsymbol\varphi
\right]
\cdot
\nabla_{\boldsymbol p}f
\\
&-
5\frac{\langle\boldsymbol\varphi,\boldsymbol p\rangle}{m(t)}f .
\end{split}
\end{equation}
For the density coefficient \(\mathfrak f=\mathbb I f\) with respect to the
ordinary almost cosymplectic volume, the same transport terms occur, while the
last conformal source coefficient changes from \(5\) to \(4\):
\begin{equation}\label{3D-LCCos-mass-Vlasov-physical}
\begin{split}
\frac{\partial \mathfrak f}{\partial \lambda}
+
\frac{\boldsymbol p}{m(t)}
\cdot
\nabla_{\boldsymbol q}\mathfrak f
={}&
-
t\frac{\langle\boldsymbol\varphi,\boldsymbol p\rangle}{m(t)}
\frac{\partial \mathfrak f}{\partial t}
\\
&-
\left[
\frac{\langle\boldsymbol\varphi,\boldsymbol p\rangle}{m(t)}
\boldsymbol p
+
\frac{|\boldsymbol p|^2}{2m(t)}
\left(
t\frac{\dot m(t)}{m(t)}-1
\right)
\boldsymbol\varphi
\right]
\cdot
\nabla_{\boldsymbol p}\mathfrak f
\\
&-
4\frac{\langle\boldsymbol\varphi,\boldsymbol p\rangle}{m(t)}
\mathfrak f .
\end{split}
\end{equation}
Hence, in the three-dimensional time-dependent mass example, the LCCos
deformation produces three additional effects in the Vlasov equation: transport
in the distinguished time direction, momentum-space drift controlled by the Lee
vector \(\boldsymbol\varphi\), and a conformal source term. In the zero-Lee
limit, these terms disappear:
\begin{equation}\label{3D-cosymplectic-mass-Vlasov}
\frac{\partial f}{\partial \lambda}
+
\frac{\boldsymbol p}{m(t)}
\cdot
\nabla_{\boldsymbol q}f
=
0 .
\end{equation}
Thus the distribution is transported along the free spatial flow on each fixed
time-slice.

\noindent \textbf{Evolution Dynamics.}
For the ordinary cosymplectic evolution vector field \(E_H\), the normalization
is \(\iota_{E_H}\tau=1\). Therefore the time variable itself becomes the
evolution parameter. The cosymplectic evolution equations are
\begin{equation}\label{TD-geodesic-cos-Evo}
    \dot t=1,
    \qquad
    \dot q^i
    =
    g^{ij}(t,q)p_j,
    \qquad
    \dot p_i
    =
    -\frac{1}{2}
    \frac{\partial g^{jk}}{\partial q^i}(t,q)p_jp_k .
\end{equation}
In Newton form, this gives
\begin{equation}\label{TD-geodesic-cos-Evo-Newton}
    \ddot q^\ell
    +
    \Gamma^\ell_{ij}(t,q)\dot q^i\dot q^j
    =
    -
    g^{\ell i}
    \frac{\partial g_{ij}}{\partial t}
    \dot q^j .
\end{equation}
Thus, in contrast with the Hamiltonian vector field, the evolution vector field
sees the explicit time-dependence of the metric.

\noindent \textbf{LCCos Evolution Dynamics.}
Finally, we consider the LCCos evolution dynamics for the same Hamiltonian. In
a local chart with \(\mathbb I|_\alpha=e^{\sigma_\alpha(t,q)}\), substituting
\eqref{TD-geodesic-H} into the LCCos evolution equations gives
\begin{equation}\label{TD-geodesic-LCCos-Evo}
\begin{split}
    \dot t
    &=
    \frac{1}{1-t\phi}
    \left(
    e^{\sigma_\alpha(t,q)}
    +
    t\varphi_i g^{ij}p_j
    \right),
\\
    \dot q^i
    &=
    g^{ij}p_j,
\\
    \dot p_i
    &=
    -
    \frac{1}{2}
    \frac{\partial g^{jk}}{\partial q^i}p_jp_k
    +
    \frac{C_{ij}}{1-t\phi}g^{jk}p_k
    +
    \frac{\varphi_i}{2(1-t\phi)}g^{jk}p_jp_k
    -
    \frac{t\varphi_i}{2(1-t\phi)}
    \frac{\partial g^{jk}}{\partial t}p_jp_k
    +
    \frac{p_i\phi e^{\sigma_\alpha(t,q)}}{1-t\phi}.
\end{split}
\end{equation}
Compared with \eqref{TD-geodesic-LCCos-Ham}, the evolution dynamics contains
the additional unit contribution in the \(\dot t\)-equation and the additional term
$p_i\phi e^{\sigma_\alpha(t,q)}/(1-t\phi)$
in the momentum equation. These are precisely the
terms caused by the normalization
\(\iota_{E_H}\tau_\varphi=\mathbb I\).

In Newton form, with \(v^i=\dot q^i\), this gives
\begin{equation}\label{TD-geodesic-LCCos-Evo-Newton}
\begin{split}
    \ddot q^\ell
    +
    \Gamma^\ell_{ij}(t,q)\dot q^i\dot q^j
    &=
    -
    g^{\ell i}
    \frac{\partial g_{ij}}{\partial t}
    \dot t\,\dot q^j
\\
    &\quad
    +
    \frac{1}{1-t\phi}
    \left(
    \varphi_i\dot q^i\,\dot q^\ell
    -
    \frac{1}{2}g^{\ell i}\varphi_i
    g_{ab}\dot q^a\dot q^b
    +
    \frac{t}{2}g^{\ell i}\varphi_i
    \frac{\partial g_{ab}}{\partial t}\dot q^a\dot q^b
    +
    \phi\,\dot q^\ell e^{\sigma_\alpha(t,q)}
    \right),
\end{split}
\end{equation}
where now
\begin{equation}\label{TD-geodesic-LCCos-Evo-tdot}
    \dot t
    =
    \frac{1}{1-t\phi}
    \left(
    e^{\sigma_\alpha(t,q)}
    +
    t\varphi_i\dot q^i
    \right).
\end{equation}
Thus the LCCos evolution vector field keeps an intrinsic time-propagation term.
This is the locally conformal analogue of the ordinary cosymplectic evolution
dynamics, where \(\dot t=1\). Consequently, the evolution dynamics describes a
genuine time-dependent geodesic motion, deformed by the Lee form and by the
conformal normalization of the evolution vector field.

\section{Summary and Discussion}

In this work, we have developed locally conformal formulations of continuum
dynamics in the symplectic and cosymplectic settings. The common geometric
ingredient is the Lichnerowicz--de~Rham differential, recalled in
Subsection~\ref{sec-LdR}, which replaces ordinary closedness by twisted
closedness and allows local geometric data to glue conformally into global
structures.

On the LCS side, we recalled the Hamiltonian dynamics, the associated Jacobi
bracket, and the corresponding bracket--vector field relation in
Subsection~\ref{sec:lcs}. Passing to densities in
Subsection~\ref{Sec-Sym-Man} led to the locally conformal Vlasov equation
\eqref{LC-Vlasov}. The main point emphasized there is that the apparent source
term in the continuity equation should not be interpreted as a genuine
creation or loss of particles. Rather, it measures the mismatch between the
ordinary LCS volume and the conformally weighted volume determined by the Lee
form. Equivalently, after introducing the weighted density
\(\mathfrak f=\mathbb I f\), the same dynamics becomes the source-free
continuity equation \eqref{LCS-Vlasov-continuity-physical} with respect to the
ordinary volume. This framework was then applied to three models: the
Nos\'e--Hoover thermostatted ensemble in Subsection~\ref{sec-NH}, where the
classical normalized thermostat is recovered on the invariant hypersurface
\(H=0\), and the LCS-driven Vlasov model with a momentum-dependent gyroscopic term in
Subsection~\ref{sec-LCS-MV}, where the Lee form produces a momentum-dependent gyroscopic force,
and their combination in Subsection~\ref{sec:NH-LCS-EM}.

On the LCCos side, we introduced in Subsection~\ref{sec:LCCosM} a formulation
adapted to the conformal gluing of local cosymplectic Hamiltonian systems. The
main advantage of this choice is that local Hamiltonian vector fields glue
directly into a global LCCos Hamiltonian vector field. The resulting Jacobi
structure and bracket--vector field relation were then used in
Subsections~\ref{sec:LCCosH} and~\ref{Sec-LCCos-Man} to derive the LCCos
Hamiltonian, evolution, and kinetic equations, including the LCCos Vlasov
equation \eqref{LCCos-Vlasov}. As in the LCS case, the conformal source term
admits an equivalent weighted-density interpretation. The coordinate
realizations show that the Lee form produces transport in the distinguished
time direction, momentum-space drift, and conformal source terms. The
time-dependent mass Hamiltonian in Subsection~\ref{sec-time-mass} illustrates
these effects at both the particle and continuum levels and clarifies the
difference between the Hamiltonian and evolution pictures.

Overall, the paper shows that locally conformal geometry provides a natural
framework for continuum theories in which the failure of ordinary volume
preservation is not treated as an external source mechanism, but as a
geometric effect of conformal gluing. The LCS and LCCos constructions
therefore fit into a common picture: the Lee form governs the conformal
deformation of the underlying geometric structure, while the corresponding
weighted-density formulation restores the continuity character of the kinetic
equations.

Two natural directions remain open. The first is to develop the locally
conformal version of contact Hamiltonian and kinetic dynamics. The second is
to extend the present continuum construction to cocontact geometry, where
locally conformal cocontact dynamics should provide a natural framework for
time-dependent dissipative kinetic systems.

\section*{Dedication}

This paper represents one of Miroslav Grmela’s final contributions to the field. We dedicate this publication to his memory as a beloved co-author, mentor, and friend.

\section*{Acknowledgment}

MP was supported by Czech Science Foundation, project 23-05736S. MP is a member of the Nečas centre for mathematical modeling.

\section*{Declaration of generative AI and AI-assisted technologies
in the writing process}

During the preparation of this work, the authors used ChatGPT (OpenAI)
and Claude (Anthropic) for language and editorial suggestions. The authors
also used these tools to assist in cross-checking selected algebraic
computations; every such computation was independently derived and verified
by the authors. After using these tools, the authors reviewed and edited
the content and take full responsibility for the content of the
publication.

\appendix
  
  \section{Appendix} 
\subsection{Schouten--Nijenhuis Bracket and Right Contraction}\label{sec:SN}

We collect here only those multivector identities that are used in the
main text; for the general theory we refer to
\cite{vaisman-book,Marle-Jacobi}.

For decomposable multivector fields
\begin{equation}\label{decompi}
\mathcal X=X_1\wedge\dots\wedge X_m,
\qquad
\mathcal Y=Y_1\wedge\dots\wedge Y_n,
\qquad m,n\geq1,
\end{equation}
the Schouten--Nijenhuis bracket is defined by
\begin{equation}\label{SN-def1}
[\mathcal X,\mathcal Y]_{SN}
=
\sum_{i,j}(-1)^{i+j}[X_i,Y_j]\wedge
X_1\wedge\dots\wedge\widehat X_i\wedge\dots\wedge X_m
\wedge
Y_1\wedge\dots\wedge\widehat Y_j\wedge\dots\wedge Y_n ,
\end{equation}
where a hat denotes omission of the corresponding factor; the result is
a multivector field of order $m+n-1$. For $m=n=1$ the bracket reduces to
the Jacobi--Lie bracket, and for a function $F$ one has
$[X,F]_{SN}=X(F)$. It is supercommutative,
\begin{equation}\label{SN-commute}
[\mathcal Y,\mathcal X]_{SN}
=
-(-1)^{(m-1)(n-1)}[\mathcal X,\mathcal Y]_{SN} ;
\end{equation}
in particular, two bivector fields commute in the Schouten--Nijenhuis
bracket.

The (right) contraction of $\mathcal X$ in \eqref{decompi} with a
one-form $\theta$ is
\begin{equation}\label{right-cont}
\iota_\theta\mathcal X
=
\sum_{i=1}^{m}(-1)^{i+m}\langle\theta,X_i\rangle\,
X_1\wedge\dots\wedge\widehat X_i\wedge\dots\wedge X_m ,
\end{equation}
equivalently $\iota_\theta\mathcal X=\mathcal X^\sharp(\theta)$ in terms
of the musical map induced by $\mathcal X$. Since \eqref{right-cont}
involves no differentiation, it is defined for an arbitrary one-form.

Three identities are used in the main text. First, for a function $F$,
\begin{equation}\label{X-FY-2}
[\mathcal X,F\mathcal Y]_{SN}
=
F[\mathcal X,\mathcal Y]_{SN}
+
\iota_{dF}\mathcal X\wedge\mathcal Y .
\end{equation}
Second, for a vector field $X$, a bivector field $\mathcal Y$ and a
one-form $\theta$,
\begin{equation}\label{iota-SN-2-2-X}
\iota_\theta(X\wedge\mathcal Y)
=
\langle\theta,X\rangle\,\mathcal Y
+
X\wedge\iota_\theta\mathcal Y .
\end{equation}
Third, for a bivector field $\mathcal X$ and a closed one-form $\theta$,
\begin{equation}\label{iota-SN-2-2-X+}
\iota_\theta\bigl([\mathcal X,\mathcal X]_{SN}\bigr)
=
2\,[\mathcal X,\iota_\theta\mathcal X]_{SN} .
\end{equation}
Identity \eqref{iota-SN-2-2-X+} is usually stated for an exact one-form
$\theta=dF$. Both sides are local operators and a closed one-form is
locally exact, so it holds for every closed one-form; this is the form
in which it is applied to the Lee form in
Subsection~\ref{sec:LCCosM}.
 \subsection{(Almost) Poisson Manifolds and Dynamics}\label{sec-Poisson}

\noindent \textbf{Almost Poisson Manifolds.}
An almost Poisson manifold is a pair \((M,\mathbb{L})\), where \(\mathbb{L}\) is a bivector field on \(M\). 
One defines a bracket on \(C^\infty(M)\) by
\begin{equation}\label{almost-Poisson-bracket-app}
\{F,H\}=\mathbb{L}(dF,dH),
\end{equation}
which is skew-symmetric and satisfies the Leibniz rule
\begin{equation}\label{Leibniz-Poisson}
\{F,KH\}=\{F,K\}H+K\{F,H\},
\end{equation}
but in general does not satisfy the Jacobi identity. For a Hamiltonian function \(H\in C^\infty(M)\), the Hamiltonian vector field is defined by
\begin{equation}\label{almost-Poisson-Ham-app}
X_H=\iota_{dH}\mathbb{L}=\sharp_\mathbb{L}(dH)=\mathbb{L}(\bullet,dH).
\end{equation}
The relation between the bracket and the vector field is then
\begin{equation}
X_H(F)=\{F,H\}.
\end{equation}
However, in general, Hamiltonian vector fields are not closed under the Jacobi--Lie bracket, and consequently the so-called characteristic distribution, defined by Hamiltonian vector fields, fails, in general, to be involutive and hence integrable.

\medskip

\noindent \textbf{Poisson Manifolds.}
A Poisson manifold \cite{BhaskaraViswanath,DuZu05,Laurent13,vaisman-book,Weinstein98} is an almost Poisson manifold \((M,\mathbb{L})\) for which the bivector field \(\mathbb{L}\) satisfies
\begin{equation}\label{Poisson-ident-app}
[\mathbb{L},\mathbb{L}]=0,
\end{equation}
with respect to the Schouten--Nijenhuis bracket. This condition implies that the bracket in \eqref{almost-Poisson-bracket-app} also satisfies the Jacobi identity
\begin{equation}\label{Jac-id-Poisson}
\circlearrowright \{F,\{H,K\}\}=0,
\end{equation}
and hence defines a Lie algebra structure on \(C^\infty(M)\). Moreover, the Hamiltonian vector fields defined by \eqref{almost-Poisson-Ham-app} are closed under the Jacobi--Lie bracket:
\begin{equation}\label{Poisson-integ-app}
[X_F,X_H]=-X_{\{F,H\}}.
\end{equation}
Thus, for Poisson manifolds, the characteristic distribution generated by Hamiltonian vector fields is involutive and hence integrable.

\subsection{(Almost) Jacobi Manifolds and Dynamics}\label{sec-Jacobi}

\noindent \textbf{Almost Jacobi Manifolds.}
Let \(M\) be a manifold equipped with a bivector field \(\mathbb{L}^J\) and a vector field \(\mathcal{V}\). In this setting, one defines an almost Jacobi bracket on the space of smooth functions by
\begin{equation}\label{Jacobi-bracket}
\{F,H\}^{J}=\mathbb{L}^J(dF,dH)-F\,\mathcal{V}(H)+H\,\mathcal{V}(F).
\end{equation}
This bracket is skew-symmetric and satisfies the first-order differential operator identity
\begin{equation}\label{Jacobi-bracket-Leibniz}
\{F,KH\}^{J}=\{F,K\}^{J}H+K\{F,H\}^{J}-KH\,\mathcal{V}(F).
\end{equation}
Thus, unlike the Poisson bracket, it does not satisfy the Leibniz rule. In general, the almost Jacobi bracket does not satisfy the Jacobi identity.

For a given Hamiltonian function \(H\), the corresponding almost Jacobi Hamiltonian vector field is defined by
\begin{equation}\label{Ham-v-f-Jac}
X_H^J=\sharp_{\mathbb{L}^J}(dH)+H\mathcal{V}.
\end{equation}
Its action on functions is expressed by
\begin{equation}\label{bra-X-J}
X_H^J(F)=\{F,H\}^{J}+F\,\mathcal{V}(H).
\end{equation}
In particular,
\begin{equation}
X_H^J(H)=H\,\mathcal{V}(H),
\end{equation}
so that the Hamiltonian need not be preserved along the flow. In general, the family of almost Jacobi Hamiltonian vector fields is not closed under the Jacobi--Lie bracket, and hence the distribution they generate is not involutive. In this way, almost Jacobi geometry extends the almost Poisson framework beyond the purely conservative case. The almost Poisson setting is recovered when \(\mathcal{V}=0\), in which case the bracket reduces to the one induced by the bivector field \(\mathbb{L}^J\).

\medskip

\noindent \textbf{Jacobi Manifolds.} 
The necessary and sufficient condition for the bracket \eqref{Jacobi-bracket} to satisfy the Jacobi identity
\begin{equation}\label{Jac-id-Jac}
\circlearrowright \{F,\{H,K\}^{J}\}^{J}=0
\end{equation}
is that the following identities hold:
\begin{equation}\label{ident-Jac}
[\mathbb{L}^J,\mathbb{L}^J]=2\,\mathcal{V}\wedge\mathbb{L}^J,
\qquad
[\mathcal{V},\mathbb{L}^J]=0.
\end{equation}
In this case, the bracket \eqref{Jacobi-bracket} is called a Jacobi bracket, \((\mathbb{L}^J,\mathcal{V})\) is called a Jacobi pair, and the triple \((M,\mathbb{L}^J,\mathcal{V})\) is called a Jacobi manifold \cite{Kirillov,Lich,Lichnerowicz-Jacobi,Marle-Jacobi,vaisman2002jacobi}. The importance of the Jacobi identity is that the Jacobi Hamiltonian vector fields defined in \eqref{Ham-v-f-Jac} are closed under the Jacobi--Lie bracket:
\begin{equation}\label{Jac-integ}
[X_F^J,X_H^J]=-X_{\{F,H\}^{J}}^J.
\end{equation}
Hence the characteristic distribution generated by Jacobi Hamiltonian vector fields is involutive and therefore integrable. In this sense, Jacobi geometry preserves the integrability aspect of Hamiltonian dynamics while allowing a wider class of motions.

\subsection{The Lichnerowicz--de~Rham Differential}
\label{sec-LdR}

Let \(M\) be a smooth manifold endowed with a closed one-form
\(\varphi\). In locally conformal geometry, \(\varphi\) is referred to as
the Lee form. The Lichnerowicz--de~Rham differential associated with
\(\varphi\) is the degree-one operator
\begin{equation}
\label{LdR-Diff-}
d_\varphi:
\Gamma^k(M)
\longrightarrow
\Gamma^{k+1}(M),
\qquad
d_\varphi\beta
=
d\beta-\varphi\wedge\beta.
\end{equation}
A direct calculation gives
\begin{equation}
d_\varphi^2\beta
=
-d\varphi\wedge\beta.
\end{equation}
Since \(d\varphi=0\), it follows that
\begin{equation}
d_\varphi^2=0.
\end{equation}
Thus \(d_\varphi\) determines a cochain complex, whose cohomology is called
the Lichnerowicz--de~Rham cohomology of \(M\)
\cite{GuLi84,HaRy99}.

The local character of this differential becomes transparent after choosing
a sufficiently fine open cover \(\{U_\alpha\}\) of \(M\). Since the Lee form
is closed, it is locally exact, and hence there exist smooth functions
\(\sigma_\alpha\in C^\infty(U_\alpha)\) such that
\begin{equation}
\label{local-Lee-potential}
\varphi|_{U_\alpha}
=
d\sigma_\alpha.
\end{equation}
On each chart, the twisted differential is conjugate to the ordinary
de~Rham differential:
\begin{equation}
\label{LdR-local-conjugation-app}
\left.d_\varphi\beta\right|_{U_\alpha}
=
e^{\sigma_\alpha}
d\left(
e^{-\sigma_\alpha}\beta|_{U_\alpha}
\right).
\end{equation}
Accordingly, \(d_\varphi\)-closedness is locally equivalent to ordinary
closedness after conformal rescaling.

On a nonempty overlap \(U_\alpha\cap U_\beta\), one has
\begin{equation}
d(\sigma_\alpha-\sigma_\beta)
=
0.
\end{equation}
Therefore,
\begin{equation}
\label{Lee-potential-overlap-app}
\sigma_\alpha-\sigma_\beta
=
c_{\alpha\beta},
\end{equation}
where \(c_{\alpha\beta}\) is locally constant. These constants record the
obstruction to replacing the family of local Lee potentials by a single
globally defined potential.

\medskip

\noindent
\textbf{Degree-zero twisted cohomology and conformal normalization.}
At degree zero, the Lichnerowicz--de~Rham differential acts on a function
\(u\in C^\infty(M)\) according to
\begin{equation}
d_\varphi u
=
du-u\varphi.
\end{equation}
In particular, the ordinary constant unit satisfies
\begin{equation}
d_\varphi 1
=
-\varphi.
\end{equation}
Thus, unless the Lee form vanishes, the ordinary unit is not
\(d_\varphi\)-closed.

A conformal normalization is a nonzero solution of
\begin{equation}
\label{twisted-closed-I}
d_\varphi\mathbb I
=
0,
\qquad
\text{equivalently,}
\qquad
d\mathbb I
=
\varphi\,\mathbb I.
\end{equation}
The existence of such a solution as an ordinary global scalar function is a
nontrivial global condition.

Suppose that \(M\) is connected and that
\(\mathbb I\in C^\infty(M)\) is a nonzero solution of
\eqref{twisted-closed-I}. If \(\mathbb I\) vanishes at one point, then the
equation
\[
d\mathbb I
=
\varphi\,\mathbb I
\]
implies, by restriction to any smooth curve issuing from that point, that
\(\mathbb I\) vanishes along the entire curve. Connectedness would then imply
\(\mathbb I\equiv0\). Consequently, every nonzero solution is nowhere
vanishing. Division by \(\mathbb I\) therefore gives
\begin{equation}
\varphi
=
\frac{d\mathbb I}{\mathbb I}
=
d\log|\mathbb I|.
\end{equation}
Hence the existence of a nonzero global scalar solution implies
\begin{equation}
[\varphi]
=
0
\qquad
\text{in }
H^1_{\mathrm{dR}}(M).
\end{equation}

Conversely, suppose that the Lee form is globally exact:
\begin{equation}
\varphi
=
d\sigma
\end{equation}
for some \(\sigma\in C^\infty(M)\). Then
\begin{equation}
\label{exact-conformal-unit-app}
\mathbb I
=
e^\sigma
\end{equation}
is nowhere vanishing and satisfies
\begin{equation}
d_\varphi\mathbb I
=
d(e^\sigma)
-
e^\sigma d\sigma
=
0.
\end{equation}
Moreover, any two nowhere-vanishing global solutions of
\eqref{twisted-closed-I} differ by a nonzero constant on each connected
component of \(M\).

We shall call such a nowhere-vanishing global scalar solution
\(\mathbb I\) a global scalar conformal unit. The preceding argument shows
that a global scalar conformal unit exists if and only if the Lee form is
exact, that is,
\begin{equation}
[\varphi]
=
0
\qquad
\text{in }
H^1_{\mathrm{dR}}(M).
\end{equation}
Thus, global scalar conformal units exist precisely in the exact-Lee sector.

\medskip

\noindent
\textbf{Local conformal units and their overlap behavior.}
Even when the Lee form is not globally exact, conformal units always exist
locally. On \(U_\alpha\), let \(1_\alpha\) be a nowhere-vanishing locally
constant normalization:
\begin{equation}
d1_\alpha
=
0.
\end{equation}
Define
\begin{equation}
\label{local-conformal-unit-app}
\mathbb I_\alpha
=
e^{\sigma_\alpha}1_\alpha.
\end{equation}
Using \(\varphi|_{U_\alpha}=d\sigma_\alpha\), one obtains
\begin{equation}
\begin{split}
d_\varphi\mathbb I_\alpha
&=
d\left(e^{\sigma_\alpha}1_\alpha\right)
-
d\sigma_\alpha\,e^{\sigma_\alpha}1_\alpha
\\
&=
e^{\sigma_\alpha}d1_\alpha
=
0.
\end{split}
\end{equation}
Thus every \(\mathbb I_\alpha\) is a local
\(d_\varphi\)-closed conformal unit.

On an overlap \(U_\alpha\cap U_\beta\), the local representatives satisfy
\begin{equation}
\mathbb I_\alpha
=
e^{c_{\alpha\beta}}
\frac{1_\alpha}{1_\beta}
\mathbb I_\beta.
\end{equation}
They define an ordinary global scalar function precisely when the local
normalizations can be chosen so that
\begin{equation}
\label{conformal-unit-gluing-app}
e^{\sigma_\alpha}1_\alpha
=
e^{\sigma_\beta}1_\beta
\qquad
\text{on }
U_\alpha\cap U_\beta.
\end{equation}
When this compatibility condition holds, the common local expressions glue
to a nowhere-vanishing global function \(\mathbb I\), and the preceding
argument implies that the Lee form is exact.

If \([\varphi]\neq0\), condition
\eqref{conformal-unit-gluing-app} cannot hold for a nowhere-vanishing global
scalar function. In that case, the constants
\begin{equation}
g_{\alpha\beta}
=
e^{c_{\alpha\beta}}
\end{equation}
describe the transition behavior of the local conformal weights. The family
\(\{\mathbb I_\alpha\}\) should then be regarded as a system of local
representatives associated with the flat weight line bundle determined by
the Lee form, rather than as the restrictions of an ordinary global
function.

\noindent 
\textbf{Scope of the Exactness Assumption.} 
Whenever $\mathbb I$, $\mathbb I\,\mathrm{Vol}$ or
$\mathfrak f=\mathbb I f$ are used in this paper as ordinary global
scalar objects, we work in the exact-Lee sector and fix a positive
global conformal unit satisfying \eqref{twisted-closed-I}.

This assumption is not needed for the Hamiltonian formalism itself. The
LCS and LCCos Hamiltonian equations, their Jacobi brackets and the
associated divergence identities are defined for an arbitrary Lee
class; in particular, given a global Hamiltonian $H$, the local
functions $H_\alpha=e^{-\sigma_\alpha}H$ always exist and satisfy the
required conformal compatibility. Exactness enters only through the
weighted volumes and densities, and through the normalization
$\iota_{E_H}\tau=\mathbb I$ of the LCCos evolution vector field.

In the exact sector the weighted formulation admits a direct
interpretation. Setting $\Omega_0=\mathbb I^{-1}\Omega$ and
$H_0=\mathbb I^{-1}H$, one has $d\Omega_0=0$ and
$\iota_{\mathcal E_H}\Omega_0=\mathbb I^{-1}d_\varphi H=dH_0$, so that
$\mathcal E_H$ is the symplectic Hamiltonian vector field of $H_0$ with
respect to $\Omega_0$; the same computation with
$\tau_0=\mathbb I^{-1}\tau$ applies in the LCCos case. Thus the
conformal unit records precisely which rescaling restores the ordinary
symplectic, respectively cosymplectic, formalism. For a non-exact Lee
form no such global rescaling exists, the identities remain valid on
each conformal chart, and their global interpretation requires
densities valued in the flat weight line bundle determined by
$\varphi$; this formulation is not pursued here. 
  \bibliographystyle{abbrv}
\bibliography{references}

\end{document}